\documentclass[12pt,a4paper]{article}

\usepackage{defaultstyle}
\usepackage{algorithmicx}
\usepackage{algpseudocode}
\algnewcommand\algorithmicyield{\textbf{yield}}
\algnewcommand\Yield[1]{\State \algorithmicyield\ #1}
\usetikzlibrary{tikzmark, patterns}

\newcommand{\authoraddress}[2]{%
    \leftskip1cm%
    \par\medskip{\footnotesize\noindent\textsc{#1}
    \par\noindent\textit{Email address: }\href{mailto:#2}{\texttt{#2}}
}}

\makeatletter
\AddToHook{cmd/appendix/before}{%
  \def\cref@section@alias{appendix}%
  \def\cref@subsection@alias{appendix}%
}
\makeatother

\numberwithin{equation}{section}

\AtBeginBibliography{\footnotesize}

\tikzset{
    midarrow/.style={
        decoration={markings, mark=at position 0.6 with {\arrow{Straight Barb[length=2mm,width=1.5mm]}}},
        postaction={decorate}
    }
}

\tikzset{
    Rightarrow/.style={
        double equal sign distance,
        >={Implies},
        ->
    },
    triple/.style={
        -,
        preaction={draw,Rightarrow}
    },
    quadruple/.style={
        preaction={
            draw,
            Rightarrow,
            shorten >=0pt
        },
        shorten >=1pt,
        -,
        double,
        double distance=0.2pt
    },
    vertex/.style={
        minimum width=2cm,
        minimum height=2cm,
        inner sep=0pt,
        outer sep=0pt
    }
}

\AtBeginDocument{%
    \newcommand{\tktrivalcoord}[6]{%
        \begingroup
        \def\mx{#5}\def\my{#6}%
        \def\ax{#1}\def\ay{#2}%
        \def\bx{#3}\def\by{#4}%

        \pgfmathparse{(#1-#3)*(#1-#5)}%
        \ifdim\pgfmathresult pt<0pt
        \def\mx{#1}\def\my{#2}%
        \def\ax{#3}\def\ay{#4}%
        \def\bx{#5}\def\by{#6}%
        \else
        \pgfmathparse{(#3-#1)*(#3-#5)}%
        \ifdim\pgfmathresult pt<0pt
            \def\mx{#3}\def\my{#4}%
            \def\ax{#1}\def\ay{#2}%
            \def\bx{#5}\def\by{#6}%
        \fi
        \fi

        \pgfmathparse{(\ay-\my)*(\by-\my)}%
        \ifdim\pgfmathresult pt>0pt

        \pgfmathparse{%
            (\ax-\mx)*(\ax-\mx)+(\ay-\my)*(\ay-\my)
            -(\bx-\mx)*(\bx-\mx)-(\by-\my)*(\by-\my)}%
        \ifdim\pgfmathresult pt>0pt
            \def\nx{\bx}\def\ny{\by}%
            \def\fx{\ax}\def\fy{\ay}%
        \else
            \def\nx{\ax}\def\ny{\ay}%
            \def\fx{\bx}\def\fy{\by}%
        \fi

        \draw[very thick,color=green!60!black]
            (\nx,\ny)
            .. controls (\nx,{(\ny+2*\my)/3})
                        and (\fx,{(\ny+2*\my)/3})
            .. (\fx,\ny)
            -- (\fx,\fy);

        \draw[very thick,color=green!60!black]
            ({(\ax+\bx)/2},{(\ny+\my)/2})
            .. controls ({(\ax+\bx)/2},\my)
                        and (\mx,{(\ny+\my)/2})
            .. (\mx,\my);

        \else
        \PackageError{tktrival}%
            {The outer points are not on the same side of the middle point}%
            {Both outer endpoints must lie above or both below the middle endpoint.}%
        \fi
        \endgroup
    }%
}

\newcommand{\cat}[1]{\mathscr{C} \mathsf{at}_{\overline{#1}}}

\newcommand{\cptd}[1]{\mathscr{C} \mathsf{ptd}_{\overline{#1}}}

\newcommand{\atm}[2][s]{%
    \ifx b#1\mathsf{Atm}_{#2}\else
    \ifx s#1\mathsf{sAtm}_{#2}\else
    \ifx f#1\mathsf{sAtm}^{\mathrm{fa}}_{#2}\else
    \mathrm{Illegal~option}%
    \fi\fi\fi
}

\newcommand{\mol}[2][s]{%
    \ifx b#1\mathsf{Mol}_{#2}\else
    \ifx s#1\mathsf{sMol}_{#2}\else
    \ifx f#1\mathsf{sMol}^{\mathrm{fa}}_{#2}\else
    \mathrm{Illegal~option}%
    \fi\fi\fi
}

\newcommand{\gset}{\mathsf{g} \mathscr{S} \mathsf{et}}

\newcommand{\set}{\mathscr{S} \mathsf{et}}

\newcommand{\ogpos}{\mathsf{og} \mathscr{P} \mathsf{os}^{\times}}

\newcommand{\bord}{\mathscr{B} \mathsf{ord}^{\mathrm{or}}}

\newcommand{\vect}{\mathscr{V} \mathsf{ect}}

\newcommand{\centstack}[2]{\mathrel{\vcenter{\offinterlineskip\halign{$\hfil##\hfil$\cr{#1}\cr{#2}\cr}}}}

\def\centerarc[#1](#2)(#3:#4:#5)
    {\draw[#1] ($(#2)+({#5*cos(#3)},{#5*sin(#3)})$) arc (#3:#4:#5);}

\begin{document}

\title{\sc Higher Categorical Bistellar Moves}

\author{Aleksandar Ivanov}
\date{}
\maketitle

\begin{abstract}
    We construct $n$-categories $\mathscr{P}_n$ with generators given by simplices and relations given by bistellar moves between them, in any dimension $n$. We show that functors out of $\mathscr{P}_n$ into an $n$-category of theories, give the correct input for the orbifold construction for TQFTs, and subsume all previously known examples, which are all in dimensions $n \leq 3$. In particular, in dimension 3 the construction gives a highly nontrivial categorification of a kind of nonunital, noncounital, Frobenius algebra.
\end{abstract}

\tableofcontents

\section{Introduction}

The goal of this paper is to study the defining relations of orbifold data --- the input of the orbifold construction --- and to capture them into a rigorous, algorithmic framework that provides their definition in any dimension. We start off with a motivational section that situates the work in the context of defect TQFTs and the orbifold construction. This section is not necessary for understanding the technical content of the paper, but it does serve as the main reason to develop it.

\setcounter{subsection}{-1}
\subsection{Motivation}

The orbifold construction is a procedure to produce a new QFT from the data of a known QFT and a collection of some of its symmetries that have to come in a particular form. Intuitively, one thinks of the procedure as an implementation of gauging the chosen discrete symmetries. It was introduced in \cite{ffrs} for two-dimensional, rational CFTs, then subsequently studied more concretely in \cite{orbifold2d} and \cite{orbifoldnd} for oriented TQFTs in arbitrary dimension, and has since been further developed by more authors; see \cite{orbifoldsurvey} for a survey. The construction can also be applied very broadly, not just to topological or rational conformal field theories, as in \cite[]{ising}. The symmetries that go into the construction, which generically mean not only group symmetries, but all noninvertible and higher-categorical symmetries, have to come packaged in a special form for the construction to work. We say that they form an \emph{orbifold datum} for the theory. In low dimension $n$ it is known that:
\begin{itemize}
    \item $(n=1)$ one-dimensional orbifold data are always idempotents, and every idempotent is a framed (or equivalently oriented) orbifold datum, while unoriented, one-dimensional orbifold data are self-adjoint idempotents.
    \item $(n=2)$ unital, counital, symmetric, special\footnote{The condition of being special here means that the multiplication is a left inverse to the comultiplication. Sometimes this goes under the name $\Delta$-separability, which has the benefit of making the connection to separability more apparent, but it has the drawback of including the free variable $\Delta$ in its name, which, if not aesthetically unappealing, is at least logically inconsistent, if one does not use $\Delta$ to denote the comultiplication.} Frobenius algebras give rise to two-dimen\-sional oriented orbifold data \cite{orbifold2d}.
    \item $(n=3)$ a certain categorification of the former algebraic structure (see \cite{orbifoldnd}, \cite{orbifold3d}), with associativity and Frobeniusness up to isomorphism, give rise to some three-dimen\-sional oriented orbifold data, of which spherical fusion categories are a strict subset of examples \cite{sphericalfusion}.
\end{itemize}
Additionally, there is a theory of condensation monads, developed by \cite{condensation}, which plays a similar role to orbifold data, as a state sum construction internal to an ambient category. It is expected that this is the correct definition of `framed orbifold data', i.e. the data needed when working with framed TQFTs. The dependence on tangential structure, however, still needs to be understood more deeply.

To be more specific, what the orbifold construction starts with, in the most studied case --- that of oriented TQFTs --- is the information of the theory, and its symmetries, in the form of a \emph{defect TQFT} (see \cite{dkr} for the two-dimensional case and \cite{orbifoldnd} for arbitrary dimension), which is a symmetric monoidal functor
\begin{equation}
    \mathscr{Z}: \bord_n (\mathbb{D}) \xrightarrow{\qquad} \vect_{\mathbb{k}},
\end{equation}
from a certain category of \emph{stratified} bordisms, labeled by the \emph{defect data} $\mathbb{D}$. An example of a morphism in $\bord_2 (\mathbb{D})$ going from the circle with three point defects to the circle with two points is
\begin{equation}
    \raisebox{-0.5\height}{\includegraphics{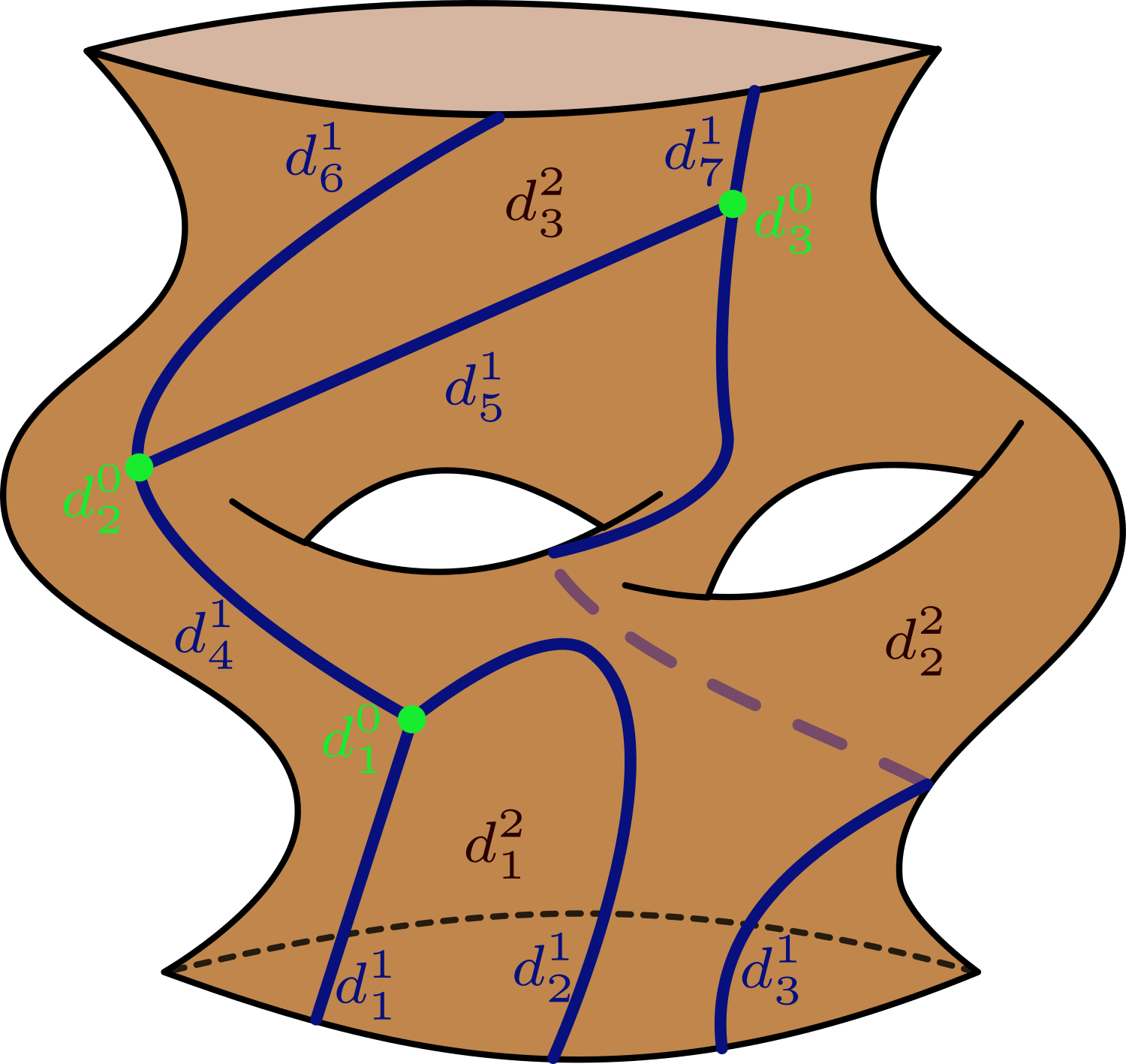}}
\end{equation}
\noindent
where the $d^m_i$ are some elements of $\mathbb{D}$ such that $m$ denotes the dimension in which the element is, and $i$ is just a counter. $\mathbb{D}$ also contains information on which labels are allowed to meet, for example $d^1_1$, $d^1_2$ and $d^1_4$ above are allowed to meet in $d^0_1$. This is an implementation of the paradigm that generalized symmetries should be identified with topological defects. If there is only one possible label in $\mathbb{D}$ for $n$-dimensional (or bulk) pieces, then restriction along the functor
\begin{equation}
    \bord_n \xhookrightarrow{\qquad} \bord_n (\mathbb{D}),
\end{equation}
which picks out the trivially stratified bordisms with the unique label, encodes the data of a TQFT which we think of as the underlying bulk TQFT, while the rest of the data encodes the symmetries of that theory. If there are multiple possible labels for bulk pieces then we have an underlying TQFT for each possible label, and we think of these as different possible phases of the theory. The defects between those encode the symmetries and dualities of the phases.\footnote{In this informal setting, we will use the words phase and theory somewhat interchangeably.}

The way orbifold data $\mathscr{A}$ arise, in this language, is through analyzing triangulations of manifolds. They can be thought of as a form of state sum data internal to the defect TQFT. The defining characteristics of orbifold data are that:
\begin{itemize}
    \item they can label a bordism stratification that is Poincar\'e dual to a triangulation, and
    \item under a change of triangulation through bistellar (also called Pachner) moves, which induce a change of stratification in the Poincar\'e dual picture, the stratified manifolds are evaluated by the theory to the same thing.
\end{itemize}

These properties allow one to define the value of the new theory $\mathscr{Z_{A}}$ on a bordism, by filling the bordism with a condensate of defects that is Poincar\'e dual to a triangulation, and then evaluating with the old defect theory $\mathscr{Z}$. The invariance under bistellar moves provides the well-definition of this construction because of Pachner's theorem:
\begin{theorem}[{\cite{pachner}}]
    Two simplicial complexes which are related by a bistellar move have PL-homeomorphic topological realizations. Furthermore, any two triangulations of a PL-manifold can be transformed into one another by some series of bistellar moves.
\end{theorem}

More specifically, a bistellar move on an $n$-dimensional simplicial complex $K$ is an operation where we first recognize a piece $A \subset \partial \Delta^{n+1}$ of the boundary of the standard simplex in one dimension higher, as a simplicial subcomplex of $K$. Performing the operation then means cutting out $A$ from $K$ and replacing it with its complement in $\partial \Delta^{n+1}$:
\begin{equation}
    K \quad \rightsquigarrow \quad \left( K \setminus A \right) \mathop{\cup}_{\partial A} \left(\partial \Delta^{n+1} \setminus A \right).
\end{equation}
An example in dimension 2 is:
\vspace{-0.2cm}
\begin{equation}
    \raisebox{-0.5\height}{\includegraphics{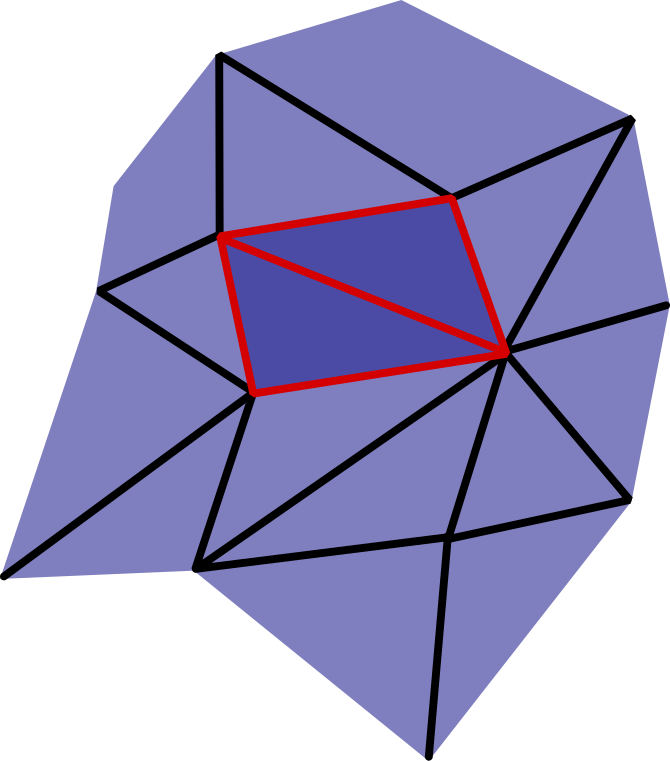}} \qquad \scalebox{2}{$\rightsquigarrow$} \qquad \raisebox{-0.5\height}{\includegraphics{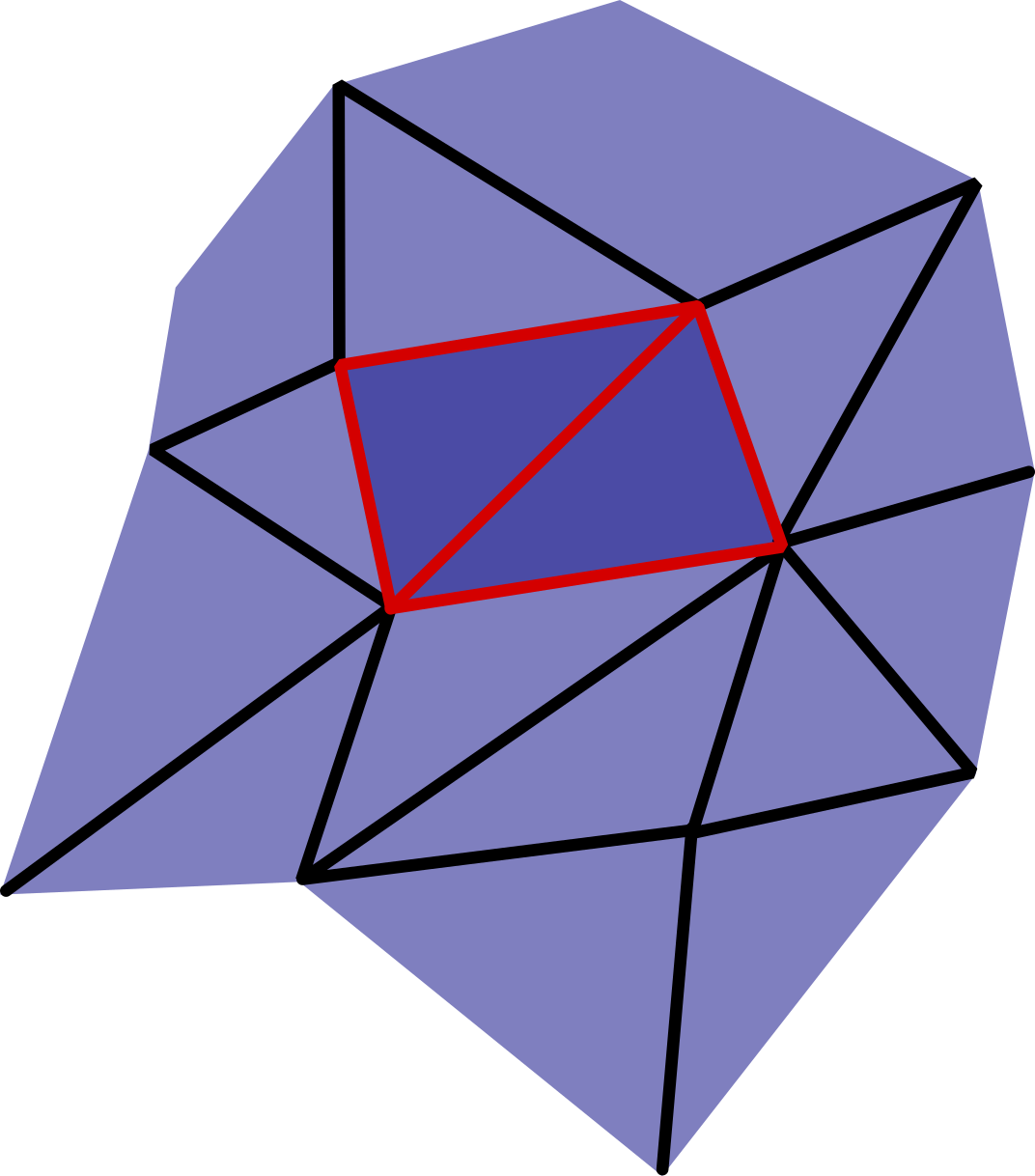}} 
\end{equation}

There is a different, more algebraic or categorical, perspective that one can take on both defects theories and the orbifold construction. In this way of thinking we consider $n$-dimensional theories/phases as objects of some $n$-category (possibly a weak $n$-category, or even an $(\infty, n)$-category), and then we think of the codimension-$k$ defects ($0 < k \leq n$) as $k$-morphisms in this higher category. Topological defects of the same codimension can fuse together in many different directions, which gives the composition operations of the higher category, and the transparent topological defect is always present, and provides the identities of the category. One inspiration for why such an encoding might be helpful can come from the cobordism hypothesis with singularities \cite{cobordismhypothesissingularities}.

Given a defect TQFT we can translate some of its data to this algebraic perspective. That is, for any given $n$-dimensional defect TQFT we can construct an $n$-category in which the objects are the phases of the theory. This was originally done by \cite{dkr} for 2-dimensional, non-extended, oriented theories, \cite{graycat} for 3-dimensional, non-extended oriented theories, and is still an active area of research. The idea is that we generate the lower levels of the higher category from the defect data itself, except for the top level which is constructed differently. For the top level, given two free composites $\underline{\smash x}$ and $\underline{\smash y}$ of codimension-$(n-1)$ defects as source and target, the set of $n$-morphisms between them is the underlying set of
\begin{equation}
    \mathscr{Z}\left(S^{n-1}(\underline{\smash x}; \underline{\smash y})\right),
\end{equation}
where $S^{n-1}(\underline{\smash x}; \underline{\smash y})$ is the $(n-1)$-sphere stratified by the presence of both input and output defects. The compositions are then given by evaluating bordisms. For example in dimension 2, the vertical and horizontal compositions, respectively, are provided by evaluating the bordisms
\vspace{-0.3cm}
\begin{equation}
    \raisebox{-64pt}{\begin{picture}(87,135)
        \put(0,10){\scalebox{.50}{\includegraphics{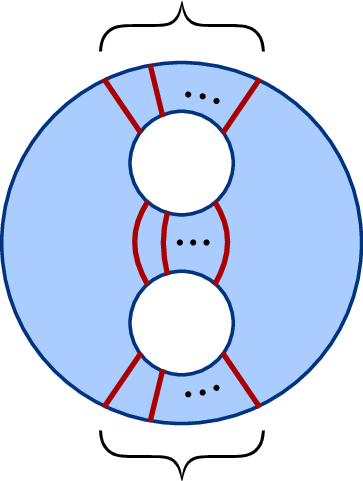}}}
        \put(0,10){
            \setlength{\unitlength}{.50pt}\put(-213,-300){
            \put(233,408)   {\scriptsize $b$ }
            \put(358,408)   {\scriptsize $a$ }
            \put(296,288)   {\scriptsize $\underline{\smash x}^*$ }
            \put(296,539)   {\scriptsize $\underline{\smash z}$ }
            \put(255,352)   {\scriptsize $x_1$ }
            \put(332,352)   {\scriptsize $x_p$ }
            \put(260,410)   {\scriptsize $y_1$ }
            \put(328,412)   {\scriptsize $y_q$ }
            \put(255,472)   {\scriptsize $z_1$ }
            \put(332,473)   {\scriptsize $z_r$ }
        }}
    \end{picture}}
    \hspace{5em}
    \raisebox{-64pt}{\begin{picture}(90,135)
        \put(-5,10){\scalebox{.50}{\includegraphics{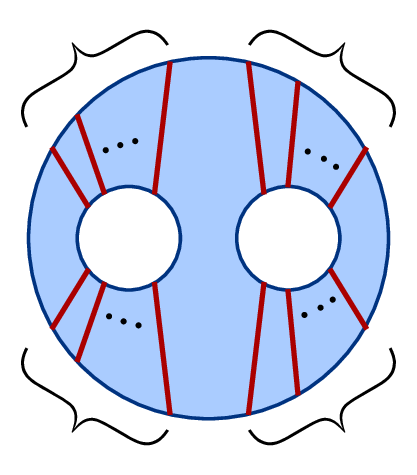}}}
        \put(-5,10){
            \setlength{\unitlength}{.50pt}\put(-200,-301){
            \put(217,410)   {\scriptsize $c$ }
            \put(297,410)   {\scriptsize $b$ }
            \put(374,410)   {\scriptsize $a$ }
            \put(363,307)   {\scriptsize $\underline{\smash x}^*$ }
            \put(227,307)   {\scriptsize $\underline{\smash y}^*$ }
            \put(363,516)   {\scriptsize $\underline{\smash x}'$ }
            \put(227,516)   {\scriptsize $\underline{\smash y}'$ }
            \put(307,375)   {\scriptsize $x_1$ }
            \put(366,392)   {\scriptsize $x_p$ }
            \put(281,375)   {\scriptsize $y_q$ }
            \put(220,392)   {\scriptsize $y_1$ }
            \put(307,452)   {\scriptsize $x_1'$ }
            \put(366,429)   {\scriptsize $x_r'$ }
            \put(281,452)   {\scriptsize $y_s'$ }
            \put(220,429)   {\scriptsize $y_1'$ }
        }}
    \end{picture}}\ ,
\end{equation}
where the bulk pieces have been labeled by $a$, $b$ and $c$, and the various line defects have been labeled by variants of $x$, $y$ and $z$, with a generic composite line defect being labeled by an underlined version of those. The depiction is adapted from \cite[Fig.6]{dkr}.

We can also look at the orbifold construction in this language --- a perspective that was initiated in \cite{orbifoldnd}. For $n$-dimensional theories, an orbifold datum is a choice of a single object and some collection of $k$-morphisms that label simplices, where $0 < k \leq n$, which satisfy the bistellar move relations at level $n$. Moreover, we observe that when we apply the construction to a theory for which we choose the trivial orbifold datum, where all defects are identity defects, then we get the theory back unchanged. Combining these two facts, at the coarsest level, the construction takes a higher category $\mathscr{T}$, whose set of objects is the set of theories we put in, and whose morphisms, and higher morphisms, are the topological defects between them, and returns an expanded set of theories $T'$, one for each choice of orbifold datum, so that
\begin{equation}
    \mathsf{ob}(\mathscr{T}) \ \ \subseteq \ \ T' \ \ := \ \ \left\{ \begin{gathered}
        \text{orbifold data} \\ \text{in } \mathscr{T}
    \end{gathered} \right\},
\end{equation}
where we identify a theory in $\mathscr{T}$ with the trivial orbifold datum on it.

Since we argued that the data of theories and their topological defects should organize into a higher category, the same applies to the theories that the orbifold construction outputs. One expects that they should also organize into a higher category, and that the orbifold construction should generalize to give not just a set of theories, but also all the topological defects between them. In other words, the orbifold construction should refine to a \emph{functor} that takes a certain appropriate kind of higher category and returns a higher category in which the objects are orbifold data.

In this paper we aim to separate the structure of orbifold data from any particular incarnation of it, and capture this structure into an abstract object, so that the representations of this object are the orbifold data.\footnote{This is very analogous to how we study group theory abstractly, and then subsequently study the representation theory of groups.} Specifically, we will define, for each $n \in \mathbb{N}$, a strict $n$-category $\mathscr{P}_n$, called the \emph{walking $n$-bistellar datum}, through a certain kind of generators-and-relations presentation. The generators will be labels for simplices, and the relations will be exactly the bistellar moves. A single $n$-dimensional orbifold datum in a higher category of theories $\mathscr{C}$ will then be a functor\footnote{This extends the analogy with group theory, where representations are functors $\mathsf{B}G \to \vect_{\mathbb{k}}$ from a one object category, which is also the case for orbifold data.}
\begin{equation}
    \mathscr{P}_n \xrightarrow{\qquad \ \ } \mathscr{C}.
\end{equation}
The work of trying to construct an $n$-category where the objects are orbifold data, and the higher morphisms are certain iterated bimodules, will be the topic of a subsequent paper.

An important technical device for the paper will be the concept of a ``frame-acyclic simplicial atom'' (\Cref{def:simplicial_atoms}), which uses concepts developed in \cite{hadzihasanovic}, and which we want to think of as a local piece of a framed triangulation. Examples in low dimensions are

\begin{equation}
    \begin{tikzpicture}[scale=0.5, baseline]
        \node at (0,0.2) [circle, fill, inner sep=1pt] {};
        \node at (2,0.2) [circle, fill, inner sep=1pt] {};
        \draw[midarrow] (0,0.2) -- (2,0.2);
    \end{tikzpicture},
    \qquad \qquad
    \begin{tikzpicture}[scale=0.6, baseline]
        \coordinate (A) at (0,{-sqrt(3)/2*1cm + 0.3cm});
        \coordinate (B) at (2cm,{-sqrt(3)/2*1cm + 0.3cm});
        \coordinate (C) at (1cm,{sqrt(3)/2*1cm + 0.3cm});

        \foreach \n in {A, B, C}
            \node at (\n) [circle, fill, inner sep=1pt] {};

        \draw[midarrow] (A) -- (B);
        \draw[midarrow] (C) -- (B);
        \draw[midarrow] (A) -- (C);

        \draw[double, double distance=1pt, -{Implies[]}] ({1cm}, {-0.25*sqrt(3)*1cm + 0.3cm}) -- ({1cm}, {0.05*sqrt(3)*1cm + 0.3cm});
    \end{tikzpicture},
    \qquad \qquad
    \begin{tikzpicture}[scale=1, baseline=-1ex]
        \coordinate (A) at (0,0-0.5);
        \coordinate (B) at (1,0-0.5);
        \coordinate (C) at (1,1-0.5);
        \coordinate (D) at (0,1-0.5);

        \foreach \n in {A, B, C, D}
            \node at (\n) [circle, fill, inner sep=1pt] {};

        \draw[midarrow] (A) -- (B);
        \draw[midarrow] (B) -- (C);
        \draw[midarrow] (C) -- (D);
        \draw[midarrow] (A) -- (D);
        \draw[midarrow] (A) -- (C);

        \draw[double, double distance=1pt, -{Implies[]}] ({0.25+0.1}, {0.75+0.1-0.5}) -- ({0.25-0.1}, {0.75-0.1-0.5});
        
        \draw[double, double distance=1pt, -{Implies[]}] ({0.75+0.1}, {0.25-0.1-0.5}) -- ({0.75-0.1}, {0.25+0.1-0.5});
    \end{tikzpicture}
    \Rrightarrow
    \begin{tikzpicture}[scale=1, baseline=-1ex]
        \coordinate (A) at (0,0-0.5);
        \coordinate (B) at (1,0-0.5);
        \coordinate (C) at (1,1-0.5);
        \coordinate (D) at (0,1-0.5);

        \foreach \n in {A, B, C, D}
            \node at (\n) [circle, fill, inner sep=1pt] {};

        \draw[midarrow] (A) -- (B);
        \draw[midarrow] (B) -- (C);
        \draw[midarrow] (C) -- (D);
        \draw[midarrow] (A) -- (D);
        \draw[midarrow] (B) -- (D);

        \draw[double, double distance=1pt, -{Implies[]}] ({0.25+0.1}, {0.25-0.1-0.5}) -- ({0.25-0.1}, {0.25+0.1-0.5});

        \draw[double, double distance=1pt, -{Implies[]}] ({0.75+0.1}, {0.75+0.1-0.5}) -- ({0.75-0.1}, {0.75-0.1-0.5});
    \end{tikzpicture}\, ,
\end{equation}
which have an underlying 1-simplex, 2-simplex, and 3-simplex, respectively. When dealing with TQFTs as actual symmetric monoidal functors from a bordism category, having no tangential structure is the easiest case to handle, but when dealing with TQFTs collectively, in the higher category of theories, the framing tangential structure is the easiest to deal with, as exemplified in the cobordism hypothesis \cite{cobordismhypothesis,cobordismhypothesissingularities}. Framed triangulations are one possibility for providing a middle ground between these two extremes by defining a geometrically inspired structure, but one which still has the directedness properties required to fit nicely into higher category theory.

\subsection{Main Results and Layout}

The goal of the paper is to provide and motivate a rigorous description of a family of strict higher categories such that mapping out of them provides the data necessary for the orbifold construction. Thus, the main definition of the paper is:
\begin{specialdefinition}[{\Cref{def:bistellar_datum}}]
    For each $n \in \mathbb{N}$, the strict $n$-category $\mathscr{P}_n$, called the \emph{walking $n$-bistellar datum}, is given by the generators-and-relations presentation
    \begin{equation}
        \mathscr{P}_n := L_n F_\omega (\atm[f]{}),
    \end{equation}
    which is the strict localization of the free strict $\omega$-category generated by the $\omega$-computad of frame-acyclic simplicial atoms.\footnote{The localization functor, the free generation of a strict $\omega$-category, and the concept of an $\omega$-computad will all be explained in \Cref{sec:prerequisites}.}
\end{specialdefinition}

The definition above does not stand alone; there are very good reasons that motivate why this definition is appropriate and captures orbifold data in all the known examples. In the latter part of this paper, we will undertake a comparison to the results of \cite{orbifold2d,orbifoldnd,orbifold3d,condensation}, which is facilitated by the generators-and-relations presentation of the definition. In particular, in low dimensions we can simply unwind this presentation and consider its implications. Thus, in dimension 1, we get the usual idempotent:
\begin{specialtheorem}[{\Cref{cor:1_bistellar_datum}}]
    The walking $1$-bistellar datum $\mathscr{P}_1$ is the 1-category with the presentation
    \begin{equation}
        \left\{ a,\  e:a \to a  \ \mid \ e \circ e = e\right\}.
    \end{equation}
\end{specialtheorem}

Already in dimension 2, we get a slight generalization of the special Frobenius algebras which are expected in the literature:
\begin{specialtheorem}[{\Cref{cor:2_bistellar_datum}}]
    The walking $2$-bistellar datum $\mathscr{P}_2$ is the 2-category with the presentation whose generators are
    \begin{equation}
        a, \quad A: a \to a, \quad \mu: A \otimes A \to A, \quad \delta: A \to A \otimes A,
    \end{equation}
    where $\mu$ and $\delta$ are called multiplication and comultiplication, respectively, and they satisfy the string diagram relations:
    \begin{equation}
        \begin{tikzpicture}[very thick, scale=0.6, color=green!50!black, baseline]
            \tktrival{1}{1}{0}{-1}{1}{1}
            \draw (0.5,-1) -- (0.5,0);
            \tktrival{1}{0}{0}{1}{1}{1}
        \end{tikzpicture} =
        \begin{tikzpicture}[very thick, scale=0.6, color=green!50!black, baseline]
            \draw (-0.5,-1) -- (-0.5,0);
            \tktrival{1}{1}{1}{-1}{1}{1}
            \tktrival{1}{0}{0}{1}{1}{1}
        \end{tikzpicture},\ \ 
        \begin{tikzpicture}[very thick, scale=0.6, color=green!50!black, baseline]
            \tktrival{-1}{-1}{-1}{1}{1}{1}
            \draw (0.5,1) -- (0.5,0);
            \tktrival{-1}{0}{0}{-1}{1}{1}
        \end{tikzpicture} =
        \begin{tikzpicture}[very thick, scale=0.6, color=green!50!black, baseline]
            \tktrival{-1}{-1}{-1}{0}{1}{1}
            \draw (-1,1) -- (-1,0);
            \tktrival{-1}{0}{0}{0}{1}{1}
        \end{tikzpicture},\ \ 
        \begin{tikzpicture}[very thick, scale=0.6, color=green!50!black, baseline]
            \draw (-1,-1) -- (-1,0);
            \tktrival{-1}{0}{0.5}{-1}{1}{1}
            \tktrival{1}{-1}{-1}{0}{1}{1}
            \draw (1,0) -- (1,1);
        \end{tikzpicture} =
        \begin{tikzpicture}[very thick, scale=0.6, color=green!50!black, baseline]
            \tktrival{1}{0}{0}{0}{1}{1}
            \tktrival{-1}{0}{0}{0}{1}{1}
        \end{tikzpicture} =
        \begin{tikzpicture}[very thick, scale=0.6, color=green!50!black, baseline]
            \tktrival{-1}{0}{-1}{-1}{1}{1}
            \draw (0.5,-1) -- (0.5,0);
            \draw (-1.5, 0) -- (-1.5, 1);
            \tktrival{1}{0}{0}{1}{1}{1}
        \end{tikzpicture},\ \ 
        \begin{tikzpicture}[very thick, scale=0.6 * 2/3, color=green!50!black, baseline=0.2cm]
            \tktrival{1}{0}{0}{0}{1.5}{1}
            \tktrival{-1}{0}{0}{0}{1.5}{1}
            \tktrival{1}{0}{0}{2}{1.5}{1}
        \end{tikzpicture} =
        \begin{tikzpicture}[very thick, scale=0.6 * 2/3, color=green!50!black, baseline=0.2cm]
            \tktrival{1}{0}{0}{0}{1.5}{1}
            \draw (0,0) -- (0,2);
        \end{tikzpicture},\ \ 
        \begin{tikzpicture}[very thick, scale=0.6 * 2/3, color=green!50!black, baseline=0.2cm]
            \tktrival{-1}{0}{0}{-1}{1.5}{1}
            \tktrival{1}{0}{0}{1}{1.5}{1}
            \tktrival{-1}{0}{0}{1}{1.5}{1}
        \end{tikzpicture} =
        \begin{tikzpicture}[very thick, scale=0.6 * 2/3, color=green!50!black, baseline=0.2cm]
            \draw (0,-1) -- (0,1);
            \tktrival{-1}{0}{0}{1}{1.5}{1}
        \end{tikzpicture}.
    \end{equation}
\end{specialtheorem}

In dimension 3, we get a genuinely new categorification of the 2-dimensional data, which we also show encompasses the known description of special orbifold data (\Cref{prop:3_bistellar_from_special_orbifold}):
\begin{specialtheorem}[{\Cref{thm:3_bistellar_datum}}]
    The walking 3-bistellar datum $\mathscr{P}_3$ is the 3-category with the presentation whose 0-, 1- and 2-generators are
    \begin{equation}
        a, \quad A: a \to a, \quad \mu: A \boxtimes A \to A, \quad \delta: A \to A \boxtimes A,
    \end{equation}
    and whose 3-generators are:
    \begin{enumerate}
        \item an associator $\alpha$, antiassociator $\hat{\alpha}$, coassociator $\gamma$ and anticoassociator $\hat{\gamma}$:
        \begin{equation}
            \begin{tikzpicture}[very thick, scale=0.6, color=green!50!black, baseline]
                \tktrival{1}{1}{0}{-1}{1}{1}
                \draw (0.5,-1) -- (0.5,0);
                \tktrival{1}{0}{0}{1}{1}{1}
            \end{tikzpicture} \quad \centstack{\stackrel{\alpha}{\Rrightarrow}}{\underset{\hat{\alpha}}{\Lleftarrow}} \quad
            \begin{tikzpicture}[very thick, scale=0.6, color=green!50!black, baseline]
                \draw (-0.5,-1) -- (-0.5,0);
                \tktrival{1}{1}{1}{-1}{1}{1}
                \tktrival{1}{0}{0}{1}{1}{1}
            \end{tikzpicture}\, , \qquad \qquad
            \begin{tikzpicture}[very thick, scale=0.6, color=green!50!black, baseline]
                \tktrival{-1}{-1}{-1}{1}{1}{1}
                \draw (0.5,1) -- (0.5,0);
                \tktrival{-1}{0}{0}{-1}{1}{1}
            \end{tikzpicture} \quad \centstack{\stackrel{\gamma}{\Rrightarrow}}{\underset{\hat{\gamma}}{\Lleftarrow}} \quad
            \begin{tikzpicture}[very thick, scale=0.6, color=green!50!black, baseline]
                \tktrival{-1}{-1}{-1}{0}{1}{1}
                \draw (-1,1) -- (-1,0);
                \tktrival{-1}{0}{0}{0}{1}{1}
            \end{tikzpicture},
        \end{equation}
        \item two Frobeniusators $\phi_L, \phi_R$ and two antiFrobeniusators $\hat{\phi}_L, \hat{\phi}_R$:
        \begin{equation}
            \begin{tikzpicture}[very thick, scale=0.6, color=green!50!black, baseline]
                \tktrival{-1}{0}{-1}{-1}{1}{1}
                \draw (0.5,-1) -- (0.5,0);
                \draw (-1.5, 0) -- (-1.5, 1);
                \tktrival{1}{0}{0}{1}{1}{1}
            \end{tikzpicture} \quad \centstack{\stackrel{\phi_L}{\Rrightarrow}}{\underset{\hat{\phi}_L}{\Lleftarrow}} \quad
            \begin{tikzpicture}[very thick, scale=0.6, color=green!50!black, baseline]
                \tktrival{1}{0}{0}{0}{1}{1}
                \tktrival{-1}{0}{0}{0}{1}{1}
            \end{tikzpicture}\, , \qquad \qquad
            \begin{tikzpicture}[very thick, scale=0.6, color=green!50!black, baseline]
                \tktrival{1}{0}{0}{0}{1}{1}
                \tktrival{-1}{0}{0}{0}{1}{1}
            \end{tikzpicture}
             \quad \centstack{\stackrel{\phi_R}{\Lleftarrow}}{\underset{\hat{\phi}_R}{\Rrightarrow}} \quad
            \begin{tikzpicture}[very thick, scale=0.6, color=green!50!black, baseline]
                \draw (-1,-1) -- (-1,0);
                \tktrival{-1}{0}{0.5}{-1}{1}{1}
                \tktrival{1}{-1}{-1}{0}{1}{1}
                \draw (1,0) -- (1,1);
            \end{tikzpicture}\, ,
        \end{equation}
        \item four bubbleators $\beta_L^U, \beta_R^U, \beta_L^D, \beta_R^D$ and four antibubbleators $\hat{\beta}_L^U, \hat{\beta}_R^U, \hat{\beta}_L^D, \hat{\beta}_R^D$:
        \begin{equation}
            \begin{tikzpicture}[very thick, scale=0.6 * 2/3, color=green!50!black, baseline]
            \draw (-1,-1) -- (-1,0);
            \tktrival{-1}{0}{0.5}{-1}{1}{1}
            \tktrival{1}{-1}{-1}{0}{1}{1}
            \draw (1,0) -- (1,1);
            \tktrival{1}{-1}{-0.5}{1}{1.5}{1}
        \end{tikzpicture} \quad \centstack{\stackrel{\hat{\beta}_R^U}{\Rrightarrow}}{\underset{\beta_R^U}{\Lleftarrow}} \quad
        \begin{tikzpicture}[very thick, scale=0.6 * 2/3, color=green!50!black, baseline]
            \tktrival{1}{0}{0}{0}{1}{1}
            \draw (0,0) -- (0,2);
        \end{tikzpicture}\, , \quad \quad
        \begin{tikzpicture}[very thick, scale=0.6 * 2/3, color=green!50!black, baseline]
            \draw (1,-1) -- (1,0);
            \tktrival{-1}{0}{-0.5}{-1}{1}{1}
            \tktrival{1}{-1}{0}{0}{1}{1}
            \draw (-1,0) -- (-1,1);
            \tktrival{1}{-1}{-1}{1}{1.5}{1}
        \end{tikzpicture} \quad \centstack{\stackrel{\hat{\beta}_L^U}{\Rrightarrow}}{\underset{\beta_L^U}{\Lleftarrow}} \quad
        \begin{tikzpicture}[very thick, scale=0.6 * 2/3, color=green!50!black, baseline]
            \tktrival{1}{0}{0}{0}{1}{1}
            \draw (0,0) -- (0,2);
        \end{tikzpicture}\, , \quad \quad 
        \begin{tikzpicture}[very thick, scale=0.6 * 2/3, color=green!50!black, baseline]
            \tktrival{-1}{-1}{-1.5}{0}{1.5}{1}
            \draw (0,0) -- (0,1);
            \tktrival{-1}{0}{-1.5}{0}{1}{1}
            \tktrival{1}{-1}{-1}{1}{1}{1}
            \draw (-2,1) -- (-2,2);
        \end{tikzpicture} \quad \centstack{\stackrel{\hat{\beta}_R^D}{\Rrightarrow}}{\underset{\beta_R^D}{\Lleftarrow}} \quad
        \begin{tikzpicture}[very thick, scale=0.6 * 2/3, color=green!50!black, baseline]
            \draw (0,-1) -- (0,1);
            \tktrival{-1}{0}{0}{1}{1}{1}
        \end{tikzpicture}\, , \quad \quad
        \begin{tikzpicture}[very thick, scale=0.6 * 2/3, color=green!50!black, baseline]
            \tktrival{-1}{-1}{-1.5}{0}{1.5}{1}
            \draw (-1.5,0) -- (-1.5,1);
            \tktrival{-1}{0}{0}{0}{1}{1}
            \tktrival{1}{-1}{-1.5}{1}{1}{1}
            \draw (0.5,1) -- (0.5,2);
        \end{tikzpicture} \quad \centstack{\stackrel{\hat{\beta}_L^D}{\Rrightarrow}}{\underset{\beta_L^D}{\Lleftarrow}} \quad
        \begin{tikzpicture}[very thick, scale=0.6 * 2/3, color=green!50!black, baseline]
            \draw (0,-1) -- (0,1);
            \tktrival{-1}{0}{0}{1}{1}{1}
        \end{tikzpicture}\, .
        \end{equation}
    \end{enumerate}
    These 3-generators, all of which are noninvertible, satisfy a set of 620 relations which are organized and presented in \Cref{app:4_atoms_list}.
\end{specialtheorem}

In more detail the layout of the paper is as follows: In \Cref{sec:prerequisites} we elaborate on the technical setting in which our results are set. There are two key components that we define. The first is the appropriate concept of a generators-and-relations presentation of a strict $n$-category through strict computads (\Cref{def:gen_and_rel_presentation}), which is how we provide the definition of the walking $n$-bistellar datum. The second is the combinatorics of simplicial atoms, which are certain elementary building blocks coming from the theory of pasting diagrams. These will provide the actual generators and relations themselves.

The first two subsections, \Cref{ssec:higher_categories} and \Cref{ssec:computads}, are a review of the classical concepts of strict higher categories and computads, and some basic useful facts about them. \Cref{ssec:atoms_molecules} consists of the introduction of framed-acyclic atoms from the theory of pasting diagrams as developed by \cite{hadzihasanovic}, as well as the restriction to the simplicial case that is of use for the present purposes. The key technical result is:
\begin{specialtheorem}[{\Cref{prop:extend_atoms_to_molecules}}]
    Let $C$ be a strict $\omega$-category, and let there be a compatible family of maps
    \begin{equation}
        \left\{ \alpha_m: \atm[f]{m} \xrightarrow[\phantom{A}]{\qquad} C_m \right\}_{0 \leq m \leq n}
    \end{equation}
    from the set of frame-acyclic simplicial $n$-atoms into the $n$-cells of $C$. Then this family extends to a map of globular sets
    \begin{equation}
        \tilde{\alpha}_n: \mol[f]{\leq n} \xrightarrow{\qquad} C
    \end{equation}
    which respects composition
    \begin{equation}
        \tilde{\alpha}_n(U \#_k V) = \tilde{\alpha}_n(U) \circ_k \tilde{\alpha}_n(V).
    \end{equation} 
\end{specialtheorem}

In \Cref{sec:bistellar_definition} we will actually see that framed-acyclic simplicial atoms of the relevant dimensions assemble into an $\omega$-computad (\Cref{def:satm_computad}). Consequently, they give the generators and relations of a family consisting of a strict $n$-category $\mathscr{P}_n$, as provided above, for each $n \in \mathbb{N}$. This is where the technical \Cref{prop:extend_atoms_to_molecules} comes into use. The low-dimensional descriptions provided above then follow.

In \Cref{sec:special_cases} we consider special cases and examples. We explain how exactly $\mathscr{P}_2$ is related to the usual special Frobenius algebras that appear in dimension 2 (\Cref{lem:2_bistellar_decomposition} and \Cref{lem:2_bistellar_with_units}). Importantly, in dimension 3, we make a connection to the case of special orbifold data as defined in \cite{orbifoldnd} and \cite{orbifold3d} (\Cref{prop:3_bistellar_from_special_orbifold}). We also comment throughout on the question of tangential structure, since the most explored case in the literature is the oriented case, while the data we provide should be the data playing a role in the framed case.

Finally, in \Cref{app:pseudocode}, we provide pseudocode that finds the set of simplicial $n$-atoms for a given $n$, and which can run very efficiently for $n \leq 4$. This pseudocode is the logic behind the \texttt{python} implementation provided as an auxiliary file \cite{pythoncode}. In \Cref{app:4_atoms_list} we provide some of the output of this code, which in particular generates the list of all simplicial 4-atoms.

\subsection*{Notation:}\label{ssec:notation}

\begin{itemize}
    \item $\mathbb{N}$ denotes the set $\{0, 1, 2, \dots \}$ including 0, and $\omega$ denotes the ordinal that comes after all the natural numbers.
    \item When specifically working in a 1-category, $\circ$ denotes composition.
    \item When specifically working in a 2-category, $\otimes$ denotes the horizontal composition, while $\circ$ denotes the vertical composition.
    \item When specifically working in a 3-category, the three different possible compositions are denoted by $\boxtimes$, at the lowest level, $\otimes$ at the middle level, and then $\circ$ at the highest level.
    \item We read string diagrams in a 2-category right to left for the horizontal composition, and bottom to top for the vertical composition.
    \item In categorical contexts, $=$ denotes a canonical isomorphism between objects.
    \item $\cat{n}$ denotes the category of strict $n$-categories, with the overline emphasizing the fact that the objects are strict. 
\end{itemize}

\subsection*{Acknowledgements:}

We thank Nils Carqueville, Aaron Hofer, Tim L\"uders, Ingo Runkel, \"Od\"ul Tetik and Markus Zetto for insightful conversations and comments. We further thank Nils Carqueville, Tim L\"uders and \"Od\"ul Tetik for a thorough reading of an earlier draft of the manuscript. The author is supported by the Austrian Science Fund FWF under project no.\ P37046.

\section{Prerequisites}\label{sec:prerequisites}

This section provides an overview of the definitions and results that serve as the basis for the main results in this paper. First we will look at globular sets and strict higher categories, which will be the underlying form that our combinatorial data will fit in. Then we will move on to computads, which are a gadget for freely generating higher categories. Most of the results on strict higher categories and computads that we will need can be found in \cite{ara}. Finally, we will introduce the combinatorial notions that are relevant for the construction, which will be atoms and molecules from the theory of pasting diagrams. Our main resource for results about these in the remainder of the article will be \cite{hadzihasanovic}.

While a large part of this section is review, some of the definitions and results are original, especially in \Cref{ssec:computads} and \Cref{ssec:atoms_molecules}. The most important of these are the definitions of simplicial atoms and molecules, and the technical result \Cref{prop:extend_atoms_to_molecules}, which underpins the whole mechanism of the construction.

\setcounter{subsection}{-1}
\subsection{Strict Higher Categories}\label{ssec:higher_categories}

\begin{definition}
    The \emph{globe category} $\mathbb{G}$ is the category with the natural numbers $0, 1, 2, \dots$ as objects and the morphisms
    \begin{equation}
        s^n_k, t^n_k : k \xrightarrow{\quad} n \qquad \forall\, k < n \in \mathbb{N},
    \end{equation}
    called the source and target morphisms, as the nontrivial morphisms. The composition is determined by the globularity relations
    \begin{equation}
        s^n_k s^k_m = s^n_l t^l_m, \qquad \qquad t^n_k s^k_m = t^n_l t^l_m.
    \end{equation}
    This category is clearly generated by the morphisms $s^{n+1}_n$ and $t^{n+1}_n$. Often the number $n$ will be implied from context in which case we will call the source and target maps $s_k, t_k: k \to n$ and furthermore if $k = n - 1$ then we may simply write $s$ and $t$. The full subcategory on integers less than, or equal to, a fixed $m$ will be denoted by $\mathbb{G}_{\leq m}$.
\end{definition}

\begin{definition}
    A \emph{globular set} is a presheaf on the globe category. The category of globular sets is
    \begin{equation}
        \gset := \mathsf{Fun}(\mathbb{G}^{\mathrm{op}}, \set).
    \end{equation}
    For a fixed globular set $X$, the images of the morphisms $s^n_k$ and $t^n_k$ may again be denoted by those same symbols so that
    \begin{equation}
        s_k, t_k: X_{k+1} \xrightarrow{\quad} X_k \qquad \forall\, k < n \in \mathbb{N},
    \end{equation}
    and the conventions surrounding the dropping of indices persist as well. An \emph{$m$-globular set} is a presheaf on the category $\mathbb{G}_{\leq m}$.
\end{definition}

\begin{definition}
    A \emph{strict $\omega$-category} $C$ is a globular set
    \begin{equation}
        \begin{tikzcd}
            \dots \arrow[r, shift left=0.5ex, "\pi^-"] \arrow[r, shift right=0.5ex, "\pi^+"'] & C_n \arrow[r, shift left=0.5ex, "\pi^-"] \arrow[r, shift right=0.5ex, "\pi^+"'] & C_{n-1} \arrow[r, shift left=0.5ex, "\pi^-"] \arrow[r, shift right=0.5ex, "\pi^+"'] & \dots \arrow[r, shift left=0.5ex, "\pi^-"] \arrow[r, shift right=0.5ex, "\pi^+"'] & C_0,
        \end{tikzcd}
    \end{equation}
    with sources and targets denoted by $\pi$, together with identity maps $i_{n}: C_n \to C_{n+1}$ for all $n \in \mathbb{N}$ and composition maps
    \begin{equation}
        C_n \mathop{\times}_{C_k} C_n \xrightarrow{\ \ \circ^n_k \ \ } C_n \qquad \forall\, k < n \in \mathbb{N},
    \end{equation}
    where the pullback is formed with one source and one target map. When $n$ is implicit from the context then we will denote the composition maps simply by $\circ_k$. We call an element of the domain of $\circ^n_k$ a pair of \emph{$k$-composable} $n$-cells. The identity and composition maps have to satisfy the following properties:
    \begin{enumerate}
        \item (Sources and targets of composites) For $k < n \in \mathbb{N}$, $\varepsilon \in \{ -, + \}$, and $x_+$ and $x_-$ $k$-composable $n$-cells, \footnote{Here we adhere to the convention that composition is read from right to left.}
        \begin{equation}
            \pi^\varepsilon (x_+ \circ_k x_-) = \begin{cases}
                \pi^\varepsilon (x_\varepsilon) & k = n-1 \\
                \pi^\varepsilon (x_+) \circ_k \pi^\varepsilon (x_-) & k < n-1
            \end{cases}.
        \end{equation}
        \item (Sources and targets of identities) For $n \in \mathbb{N}$, $\varepsilon \in \{ -, + \}$, and $x$ an $n$-cell, 
        \begin{equation}
            \pi^\varepsilon i_n (x) = x.
        \end{equation}
        \item (Strict associativity for composition) For $k < n \in \mathbb{N}$ and $x_1$, $x_2$ and $x_3$, three $n$-cells, 
        \begin{equation}
            (x_1 \circ_k x_2) \circ_k x_3 =  x_1 \circ_k (x_2 \circ_k x_3).
        \end{equation}
        \item (Strict identities for composition) For $k < n \in \mathbb{N}$, $\varepsilon \in \{ -, + \}$, and $x$ an $n$-cell, 
        \begin{equation}
            x \circ_k i_{n-1} \dots i_{n-k} \pi^-_{n-k}(x) = x = i_{n-1} \dots i_{n-k} \pi^+_{n-k} (x) \circ_k x.
        \end{equation}
        \item (Strict binary interchange) For $k < l < n \in \mathbb{N}$, and $n$-cells $x_1$, $x_2$, $x_3$ and $x_4$, such that the ordered pairs $(x_1, x_2)$ and $(x_3, x_4)$ are both $l$-composable, while the ordered pairs $(x_1, x_3)$ and $(x_2, x_4)$ are both $k$-composable, 
        \begin{equation}
            (x_1 \circ_l x_2) \circ_k (x_3 \circ_l x_4) = (x_1 \circ_k x_3) \circ_l (x_2 \circ_k x_4).
        \end{equation}
        \item (Strict nullary interchange) For $k < n \in \mathbb{N}$, and $x_1$ and $x_2$ $k$-composable $n$-cells,
        \begin{equation}
            i_n (x_1 \circ_k x_2) = i_n(x_1) \circ_k i_n (x_2).
        \end{equation}
    \end{enumerate}
    If the globular set is instead an $m$-globular set, and the composition and identity maps' indices are similarly restricted by $m$, then we have a \emph{strict $m$-category} instead.
\end{definition}

\begin{remark}
    We already used the terminology $k$-composable pairs of $n$-cells for the elements of $C_n \mathop{\times}_{C_k} C_n$. If, for some pair of $n$-cells, we do not want to be specific about what $k$ we need to choose so that the pair is in $C_n \mathop{\times}_{C_k} C_n$, but such a $k$ nonetheless exists, then we say that the pair is \emph{composable}.
\end{remark}

\begin{remark}
    This definition of strict higher categories sees them as globular sets with extra structure. There is a different definition, sometimes called the `single-set' definition, in which all the cells live in a single set, where they stand-in both for the cells themselves and for the identities on the cells. The source, targets and compositions are then added as extra structure.
\end{remark}

\begin{definition}
    A map $F\colon C \to D$ between strict $m$-categories, for $m \in \mathbb{N} \cup \{ \omega\}$, otherwise known as a \emph{strict $m$-functor}, is a map of the underlying $m$-globular sets such that
    \begin{equation}
        \begin{tikzcd}
            C_n \mathop{\times}_{C_k} C_n \arrow[r, "\circ_k"] \arrow[d, "{(F_n, F_n)}"'] & C_n \arrow[d, "F_n"] \\
            D_n \mathop{\times}_{D_k} D_n \arrow[r, "\circ_k"] & D_n
        \end{tikzcd}\, , \qquad\qquad
        \begin{tikzcd}
            C_n \arrow[r, "i_n"] \arrow[d, "F_n"'] & C_{n+1} \arrow[d, "F_{n+1}"] \\
            D_n \arrow[r, "i_k"] & D_{n+1}
        \end{tikzcd}
    \end{equation}
    commute for all $k < n$, and $n$ a valid index for a set of cells of $C$ and $D$.
\end{definition}

\begin{lemma}
    Strict $n$-categories and strict $n$-functors, for each $n \in \mathbb{N} \cup \{ \omega\}$, assemble into a category, which we will denote by $\cat{n}$.
\end{lemma}

The main operations that we will need between higher categories of different dimensions famously come in the form of two adjunctions. We provide this result, and an outline of its proof, mostly to fix what the functors do concretely. For more details see e.g. \cite[\S 14.4.6]{ara}.

\begin{proposition}\label{prop:cat_adjunctions}
    The inclusion $I_n$ of strict $n$-categories into strict $(n+1)$-categories, by adding only identity cells at level $n+1$, has both a left and a right adjoint
    \begin{equation}
        L_n, R_n: \cat{n+1} \xrightarrow{\quad} \cat{n},
    \end{equation}
    where $R_n$ is the truncation functor which forgets the level $n+1$, while $L_n$ is the localization functor that identifies all $n$-cells that are related by a zigzag of $(n+1)$-cells.
\end{proposition}

\begin{proof}
    The inclusion functor $I_n$ sends a strict $n$-category $C$ to the strict $(n+1)$-category with underlying $(n+1)$-globular set
    \begin{equation}
        \begin{tikzcd}
            C_n \arrow[r, shift left=0.6ex, Rightarrow, no head] \arrow[r, shift right=0.6ex, Rightarrow, no head] & C_n \arrow[r, shift left=0.5ex, "\pi^-"] \arrow[r, shift right=0.5ex, "\pi^+"'] & C_{n-1} \arrow[r, shift left=0.5ex, "\pi^-"] \arrow[r, shift right=0.5ex, "\pi^+"'] & \dots \arrow[r, shift left=0.5ex, "\pi^-"] \arrow[r, shift right=0.5ex, "\pi^+"'] & C_0,
        \end{tikzcd}
    \end{equation}
    At the top level $n+1$ we add the identity map $i_{n} = \mathrm{id}$ and the composition map
    \begin{equation}
        C_n \mathop{\times}_{C_n} C_n = C_n \xrightarrow{\ \ \mathrm{id} \ \ } C_n,
    \end{equation}
    with everything else being inherited from level $n$. It is easy to see that this collection satisfies the axioms of a strict $(n+1)$-category.

    To see the adjunction $I_n \dashv R_n$ observe the bijection of hom-sets
    \begin{equation}
        \cat{n+1}(I_n (C), D) \xrightarrow{\ \ \cong \ \ } \cat{n} (C, R_n (D)).
    \end{equation}
    Indeed, since strict $(n+1)$-functors preserve identities we have that a functor $F: I_n (C) \to D$ satisfies $F_{n+1} = i_n F_n$, so that it is fully determined by its components up to level $n$. 

    The functor $L_n$ sends a strict $(n+1)$-category $C$ to the strict $n$-category with underlying globular set
    \begin{equation}
        \begin{tikzcd}
            C_n/_{\approx}\arrow[r, shift left=0.5ex, "\pi^-/_\approx"] \arrow[r, shift right=0.5ex, "\pi^+/_\approx"'] & C_{n-1} \arrow[r, shift left=0.5ex, "\pi^-"] \arrow[r, shift right=0.5ex, "\pi^+"'] & \dots \arrow[r, shift left=0.5ex, "\pi^-"] \arrow[r, shift right=0.5ex, "\pi^+"'] & C_0,
        \end{tikzcd}
    \end{equation}
    where $\approx$ is the congruence relation generated by the relation
    \begin{equation}
        x \approx_* y :\Leftrightarrow \left(\exists z \in C_{n+1}: \pi^- (z) = x \ \mathrm{and} \ \pi^+ (z) = y \right),
    \end{equation}
    and the sources and target maps descend to equivalence classes by globularity. The new identity map is $\tilde{i}_n = q i_n$, where $q: C_n \to C_n/_\approx$ is the quotient map, and the new composition maps are given by
    \begin{equation}
        q(x) \ \tilde{\circ}_k \ q(y) = q (x \circ_k y),
    \end{equation}
    and are well-defined because $\approx$ is a congruence relation. It is easily seen that the axioms of a strict $n$-category are satisfied.
    
    To get an adjunction $L_n \dashv I_n$, we will construct a bijection
    \begin{equation}
        \cat{n} (L_n (C), D) \xrightarrow{\ \ \cong \ \ } \cat{n+1} (C, I_n (D))
    \end{equation}
    On the right-hand side the data defining a strict $(n+1)$-functor $F: C \to I_n (D)$ reduce to the data for a strict $n$-functor together with
    \begin{equation}
        F_n \pi^- = F_{n+1} = F_n \pi^+
    \end{equation}
    On the left-hand side, the data defining a strict $n$-functor $G: L_n (C) \to D$ is exactly the data of a strict $n$-functor with components $(G_n q, G_{n-1}, \dots, G_0)$. Because of the universal property of the quotient, maps like $F_n$ (with $ F_n \pi^- = F_n \pi^+$) are bijectively related maps like $G_n q$.
\end{proof}

\begin{remark}
    The inclusion of strict $n$-categories into strict $m$-categories for $n < m$ is just the composition of the stepwise inclusion maps. Since adjunctions paste together the compositions of the stepwise left (or right) adjoint is the full left (or right) adjoint. Properly we should denote these maps with $I^m_n$, $L^m_n$ and $R^m_n$, but we will mostly abuse notation and drop $m$.
\end{remark}

\begin{remark}
    For the functors $I^\omega_n$, $L^\omega_n$ and $R^\omega_n$, we cannot get them automatically from composing stepwise adjoints, but we can define them hands-on, analogously to what we did for the stepwise versions, and we can also provide the adjunctions $L^{\omega}_n \dashv I^{\omega}_n \dashv R^{\omega}_n$ manually.
\end{remark}

\begin{lemma}\label{lem:loc_and_rest_interaction}
    Let $k < n < m \in \mathbb{N} \cup \{\omega\}$. Then
    \begin{equation}
        R_k L_n = R_k : \cat{m} \xrightarrow{\quad} \cat{k} \quad \mathrm{and} \quad L_k R_n = L_k: \cat{m} \xrightarrow{\quad} \cat{k}.
    \end{equation}
\end{lemma}

\begin{proof}
    The first statement is obvious from the definitions. The second statement reduces to comparing the congruence relation $\approx_{k+1}$ defining $C_k$ modulo $C_{k+1}$ as in the proof of \Cref{prop:cat_adjunctions}, and the congruence relation defining $C_k$ modulo $C_{k+1}/_{\approx_{k+2}}$. These two relations can easily be seen to be the same because the existence of an element $a \in C_{k+1}$ implies the existence of an equivalence class $[a]_{\approx_{k+2}}$, and vice versa.
\end{proof}

\subsection{Computads}\label{ssec:computads}

In the entirety of this work the word computad will be used for strict computads (as opposed to computads over an arbitrary globular operad). This is a concept, whose purpose is to freely generate strict higher categories. Roughly, an $n$-computad $A$ consists of a set $A_i$ of generating cells in each dimension $i$ up to $n$, together with source and target maps $b^{\pm}_i$, such that $b^-(x)$ and $b^+(x)$ are parallel in the free category that is generated by the computad. The important point is that asking the source and target of a generating cell to be generating cells is too restrictive. Instead, we want to allow the sources and targets to be composites. This means, though, that we have to define things inductively, because to be able to provide the codomain of the source and target maps we need to have defined what it means to freely generate a strict category with the data up to one level lower.

\begin{definition}
    For any $n \in \mathbb{N}$ there is a category $\cat{n}^+$ given by the pullback
    \begin{equation}
        \begin{tikzcd}
            \cat{n}^+ \arrow[rd, phantom, very near start, "\ulcorner"] \arrow[r] \arrow[d] & \gset_{n+1} \arrow[d] \\
            \cat{n} \arrow[r] & \gset_n
        \end{tikzcd}.
    \end{equation}
    An object of this category is given by an $(n+1)$-globular set whose restriction to an $n$-globular set carries the structure of an $n$-category. When $n=-1$ we define $\cat{-1}^+ := \cat{0} = \set$.
\end{definition}

\begin{proposition}[{\cite[Prop.\,15.1.3]{ara}, \cite[\S.4.1]{metayer}}]
    The functor that forgets the categorical structure on the top level $\cat{n+1} \xrightarrow{\qquad} \cat{n}^+$ has a free left adjoint
    \begin{equation}
        \tilde{F}_{n+1}: \cat{n}^+ \xrightarrow{\qquad} \cat{n+1},
    \end{equation}
    which acts as the identity on every level except the top level, where it has all well-typed composites of generator $(n+1)$-cells and identities of $n$-cells modulo relations that impose the strict $(n+1)$-category axioms.
\end{proposition}

\begin{definition}
    The category of $0$-computads is just $\cptd{0} = \set$, and it comes with the functor $F_0: \cptd{0} \to \cat{0} = \set$, which is the identity functor of $\set$. Having defined the category of $n$-computads $\cptd{n}$ together with a functor $F_n:\cptd{n} \to \cat{n}$, we define the category of $(n+1)$-computads $\cptd{n+1}$, together with the functor $F_{n+1}: \cptd{n+1} \to \cat{n+1}$ as
    \vspace{-0.3cm}
    \begin{equation}
        \begin{tikzcd}
            \cptd{n+1} \arrow[rd, phantom, very near start, "\ulcorner"] \arrow[r] \arrow[rr, bend left=3ex, "F_{n+1}"] \arrow[d] & \cat{n}^+ \arrow[r, "\tilde{F}_{n+1}"']\arrow[d] & \cat{n+1} \\
            \cptd{n} \arrow[r, "F_n"] & \cat{n} &
        \end{tikzcd}.
    \end{equation}
\end{definition}

\begin{remark}
    The above definition works in two steps. In one step we define $\cat{n}^+$, which does the job of adding new generator cells in dimension $n+1$ together with source and target maps that make the whole collection globular --- or, in other words, so that the source and target of a given $(n+1)$-cell are parallel. However, for $n>0$, not all strict $n$-categories are free, so when defining an $n$-computad we have to recursively make sure that all lower levels are free. This is the second step.
\end{remark}

\begin{remark}
    Unwinding the above definition, to define an $(n+1)$-computad $A$ we need the data of an $n$-computad $A_{\leq n}$, a set of $(n+1)$-cells $A_{n+1}$ and two maps $b^{\pm}_n$ providing the source ($-$) and target ($+$), so that
    \begin{equation}
        \begin{tikzcd}
            A_{n+1} \arrow[r, shift left=0.5ex, "b_n^-"] \arrow[r, shift right=0.5ex, "b_n^+"'] & F_n(A_{\leq n})_n \arrow[r, shift left=0.5ex, "\pi^-"] \arrow[r, shift right=0.5ex, "\pi^+"'] & F_n(A_{\leq n})_{n-1} 
        \end{tikzcd}
    \end{equation}
    is a $2$-globular set. In practice, we will do this inductively, starting with a set of 0-cells $A_0$ and then adding sets of cells with appropriate sources and targets dimension by dimension. With this approach, the most technical part of constructing a computad is providing sources and targets that satisfy the requirements.
\end{remark}

\begin{remark}
    When it is necessary to be explicit we will denote the functor that forgets an $n$-computad down to a $k$-computad by
    \begin{equation}
        V^n_k: \cptd{n} \xrightarrow{\qquad} \cptd{k}.
    \end{equation}
    We continue to adopt the abuse of dropping $n$ when it is implied from context.
\end{remark}

\begin{proposition}\label{prop:F_is_free}
    For any $n \in \mathbb{N}$ the free functor $F_n$ is a left adjoint to the forgetful functor
    \begin{equation}
        U_n: \cat{n} \xrightarrow{\qquad} \cptd{n}.
    \end{equation}
\end{proposition}

\begin{proof}
    Since we know that $\tilde{F}_n$ has a forgetful right adjoint, we just need to prove that the functor $\hat{F}_n: \cptd{n} \to \cat{n-1}^+$ has one as well. The composition of these will then be the desired right adjoint functor.

    There is a natural candidate for this right adjoint $\cat{n-1}^+ \to \cptd{n}$, which we call $\hat{U}_n$, and which is defined inductively. Given an object $C$ of $\cat{n-1}^+$, $\hat{U}_n$ returns the $n$-computad which consists of:
    \begin{enumerate}
        \item the set of $n$-cells $C_\star$, where $C_\star$ is the top level of cells of $C$,
        \item the $(n-1)$-computad $U_{n-1} (R_{n-1} C)$, where $U_{n-1}$ is defined by the inductive hypothesis,
        \item the source and target maps are given by
        \begin{equation}
            \begin{tikzcd}
                b^\varepsilon: C_\star \arrow[r, "\pi^\varepsilon"] & C_{n-1} \arrow[r, Rightarrow, no head] & (R_{n-1} C)_{n-1} \arrow[r, "\eta_{n-1}"] & (F_{n-1} U_{n-1} R_{n-1} C)_{n-1}
            \end{tikzcd},
        \end{equation}
        where $\eta$ is the unit for the adjunction $F_{n-1} \dashv U_{n-1}$. These source and target maps are globular because
        \begin{equation}
            \pi^\varepsilon b^- (x) = \pi^\varepsilon \eta_{n-1} \pi^- (x) = \eta_{n-2} \pi^\varepsilon \pi^- (x)= \eta_{n-2} \pi^\varepsilon \pi^+ (x) = \pi^\varepsilon b^+ (x).
        \end{equation}
    \end{enumerate}

    We will show the adjunction $\hat{F}_n \dashv \hat{U}_n$ by an isomorphism of hom-sets. A map in $\cat{n-1}^+$ from $\hat{F}_n (A)$ to $C$ is given by a strict $(n-1)$-functor $F_{n-1} V_{n-1} A \to R_{n-1} C$ to, together with a map $A_n \to C_\star$ such that both squares
    \begin{equation}\label{eq:cptd_catplus_cond_on_left}
        \begin{tikzcd}[cramped]
            A_n \arrow[r, shift left=0.5ex] \arrow[r, shift right=0.5ex] \arrow[d] & (F_{n-1}V_{n-1} A)_{n-1} \arrow[d] \\
            C_\star \arrow[r, shift left=0.5ex, "b^-"] \arrow[r, shift right=0.5ex, "b^+"'] & C_{n-1}
        \end{tikzcd}
    \end{equation}
    commute. On the other hand, a map of $n$-computads from $A$ to $\hat{U}_n C$ is given by a map of $(n-1)$-computads $p: V_{n-1} A \to U_{n-1} R_{n-1} C$, together with a map $A_n \to C_\star$ such that both squares
    \begin{equation}\label{eq:cptd_catplus_cond_on_right}
        \begin{tikzcd}[cramped]
            A_n \arrow[r, shift left=0.5ex] \arrow[r, shift right=0.5ex] \arrow[d] & (F_{n-1}V_{n-1} A)_{n-1} \arrow[d, "(F_{n-1}p)_{n-1}"] \\
            C_\star \arrow[r, shift left=0.5ex, "b^-"] \arrow[r, shift right=0.5ex, "b^+"'] & (F_{n-1} U_{n-1} R_{n-1} C)_{n-1}
        \end{tikzcd}
    \end{equation}
    commute. Inductively we know that there is a bijection $F_{n-1} \dashv U_{n-1}$, so that the data that makes up the two maps is the same. Moreover, the definition of the source and target maps of $\hat{U}_n C$ is exactly such that the conditions of \Cref{eq:cptd_catplus_cond_on_left} and the ones of \Cref{eq:cptd_catplus_cond_on_right} are equivalent. 

    In the base case, $F_0$ is just the identity functor on $\set$ so it clearly has a right adjoint $U_0$, which is also the identity functor. This completes the induction.
\end{proof}

\begin{lemma}\label{lem:free_cat_same_below}
    Let $k < n$ be in $\mathbb{N} \cup \{ \omega \}$. Then
    \begin{equation}
        R_k F_n = F_k V_k : \cptd{n} \xrightarrow{\qquad} \cat{k}.
    \end{equation}
\end{lemma}

\begin{proof}
    Without loss of generality, we will only show this for $k=n-1$. Observe that, because $\tilde{F}_n$ acts as the identity on levels below $n$, there is a commutative diagram
    \begin{equation}
        \begin{tikzcd}
            \cat{n-1}^+ \arrow[d] \arrow[r, "\tilde{F}_n"] & \cat{n} \arrow[dl, "R_{n-1}"] \\
            \cat{n-1}
        \end{tikzcd}.
    \end{equation}
    We get the claim by pasting this diagram with the defining diagram of $\cptd{n}$.
\end{proof}

The technology of computads allows us to formalize what we mean by a free strict $n$-category, with $n \in \mathbb{N} \cup \{ \omega \}$, but we are also interested in describing certain strict $n$-categories that depart from freeness in a controlled way, namely, quotients of free $n$-categories by sets of relations.

We can put this in more formal language by working in $\cat{n}^+$. We think of an object
\begin{equation}\label{eq:object_of_catplus}
    \begin{tikzcd}
        C_\star \arrow[r, shift left=0.5ex, "r^-"] \arrow[r, shift right=0.5ex, "r^+"'] & C_{n} \arrow[r, shift left=0.5ex, "\pi^-"] \arrow[r, shift right=0.5ex, "\pi^+"'] & C_{n-1} \arrow[r, shift left=0.5ex, "\pi^-"] \arrow[r, shift right=0.5ex, "\pi^+"'] & \dots \arrow[r, shift left=0.5ex, "\pi^-"] \arrow[r, shift right=0.5ex, "\pi^+"'] & C_0
    \end{tikzcd}
\end{equation}
of $\cat{n}^+$ as an $n$-category $C$ together with a set of relations $C_\star$, where the maps $r^-$ and $r^+$ tell us which two $n$-cells a given relation identifies.

\begin{definition}\label{def:cat_quotient_by_relations}
    There is a quotient functor
    \begin{equation}
        Q_n: \cat{n}^+ \xrightarrow{\qquad} \cat{n},
    \end{equation}
    which takes an object like in \Cref{eq:object_of_catplus} to
    \begin{equation}
        \begin{tikzcd}
            C_n /_{\sim} \arrow[r, shift left=0.5ex, "\pi^-/_{\sim}"] \arrow[r, shift right=0.5ex, "\pi^+/_{\sim}"'] & C_{n-1} \arrow[r, shift left=0.5ex, "\pi^-"] \arrow[r, shift right=0.5ex, "\pi^+"'] & \dots \arrow[r, shift left=0.5ex, "\pi^-"] \arrow[r, shift right=0.5ex, "\pi^+"'] & C_0
        \end{tikzcd},
    \end{equation}
    where the equivalence relation $\sim$ on $n$-cells $x,y \in C_n$ is the congruence relation generated by
    \begin{equation}
        x \sim_* y :\Leftrightarrow \left( \exists a \in C_\star: r^- (a) = x \text{ and } r^+(a) = y \right).
    \end{equation}
\end{definition}

\begin{proof}
    Explicitly, the relation is given by $x \sim y$ if there exists a finite sequence
    \begin{equation}\label{eq:sequence_for_relation}
        \left( (x = x^0), (x^1_{1}, \dots, x^1_{i_1}), \dots, (x^{m-1}_{1}, \dots, x^{m-1}_{i_{m-1}}), (x^m = y) \right),
    \end{equation}
    where each subsequence $(x^p_{1}, \dots, x^p_{i_p})$ for $p=1, \dots, m-1$ is a sequence of composable $n$-cells, and where for each consecutive pair of subsequences, $(x^p_{1}, \dots, x^p_{i_p})$ and $(x^{p+1}_{1}, \dots, x^{p+1}_{i_{p+1}})$, either
    \begin{enumerate}
        \item $x^p_1 \circ \dots \circ x^p_{i_p} = x^{p+1}_{1} \circ \dots \circ x^{p+1}_{i_{p+1}}$, where $\circ$ denotes the relevant composition operations determined by the composability criterion, or
        \item $i_p = i_{p+1}$ and $\forall i \in \{1, \dots, i_p\}: x^p_i \sim_* x^{p+1}_i$ or
        \item $i_p = i_{p+1}$ and $\forall i \in \{1, \dots, i_p\}: x^{p+1}_i \sim_* x^p_i$.
    \end{enumerate}

    It is clear that this relation is an equivalence relation. What we need to show is that the maps $\pi^{\varepsilon}$ for $\varepsilon \in \{ +, - \}$ descend to equivalence classes making $\pi^{\varepsilon}/_{\sim}$ well-defined, and that the top two levels still form a category.

    The first check is guaranteed by globularity. Indeed, if $x \sim y$ and the relevant sequence is as in \Cref{eq:sequence_for_relation}, where we further denote $x^p = x^p_1 \circ \dots \circ x^p_{i_p}$,  then there exists some $a^p \in C_\star$ such that
    \begin{equation}
        \pi^{\varepsilon} (x^p) = \pi^{\varepsilon} (r^- (a^p)) = \pi^{\varepsilon} (r^+ (a^p)) = \pi^{\varepsilon} (x^{p+1}),    
    \end{equation}
    or the same equation with $r^-$ and $r^+$ switched. In any case, $\pi^{\varepsilon} (x) = \pi^{\varepsilon} (y)$.

    Composition descends because we ask for a congruence relation. To be explicit, though, observe that if $(x,y)$ is a pair of composable $n$-cells such that a sequence $(x^i)_{i=0,\dots,k}$ exhibits the equivalence $x \sim x'$ and a sequence $(y^i)_{i=0,\dots,l}$ exhibits the equivalence $y \sim y'$ then the sequence
    \begin{equation}
        \left( (x \circ y), (x, y), (x^1, y), \dots, (x', y), (x', y^1), \dots, (x', y'), (x' \circ y')\right)
    \end{equation}
    exhibits the equivalence $x \circ y \sim x' \circ y'$. It is clear that the equivalence class of an identity is again an identity for this new composition, and that all axioms of an $n$-category are implied by their counterparts for $C$.

    Given any morphism $F: C \to D$ in $\cat{n}^+$, if $x \sim_* y$ with $a \in C_\star$ witnessing the equivalence then $F_n(x) \sim_* F_n(y)$ with $F_{n+1}(a)$ witnessing the equivalence, and because of functoriality this also extends to the full congruence relation $\sim$. Thus, $F$ descends to the quotient, and its functoriality is also preserved.
\end{proof}

\begin{definition}\label{def:gen_and_rel_presentation}
    A strict $n$-category has a \emph{generators-and-relations presentation} if it is in the image of the functor
    \begin{equation}
        \cptd{n+1} \xrightarrow{\quad} \cat{n}^+ \xrightarrow{\ \ Q_n \ \ } \cat{n}.
    \end{equation}
    The $k$-cells, for $k \leq n$, of the $(n+1)$-computad are the generators, while the $(n+1)$-cells are the relations.
\end{definition}

\begin{proposition}\label{prop:localization_is_presentation}
    For all $k < n \in \mathbb{N} \cup \{\omega\}$, there is a commutative diagram 
    \begin{equation}
        \begin{tikzcd}
            \cptd{n} \arrow[r, "F_{n}"] \arrow[d] & \cat{n} \arrow[rd, "L_k"] \\
            \cptd{k+1} \arrow[r] & \cat{k}^+ \arrow[r, "Q_k"] & \cat{k}
        \end{tikzcd}.
    \end{equation}
\end{proposition}

\begin{proof}
    By virtue of the second part of \Cref{lem:loc_and_rest_interaction} and \Cref{lem:free_cat_same_below} we can take $n=k+1$ without loss of generality, for the case of finite $n$. In this case we need to show that
    \begin{equation}
        L_k F_{k+1} = Q_k \hat{F}_{k+1},
    \end{equation}
    where we use the notation $\hat{F}$ as in the proof of \Cref{prop:F_is_free}.

    Given any $(k+1)$-computad $A$ the action of both $L_k F_{k+1}$ and $Q_k \hat{F}_{k+1}$ on levels less than $k$ is identical, by the construction of the functors involved. We only need to check that the sets of cells at the top level are isomorphic. 
    
    Denote the generator $(k+1)$-cells of $A$ by $A_{k+1}$ and let $B_{k+1} := F_{k+1}(A)_{k+1}$ and $B_k := F_{k+1} (A)_{k} = F_k V_k (A)_k$. At the top level each set is defined as a quotient by some equivalence relation on $n$-cells. What we need to show is then that these equivalence relations mutually imply each other. The left-hand-side equivalence relation is the one defining the localization, which is generated by 
    \begin{equation}
        (x \approx_* y) :\Leftrightarrow \left( \exists z \in B_{k+1}: \pi^-(z) = x \ \mathrm{and} \ \pi^+(z) = y \right),
    \end{equation}
    while the right-hand-side equivalence relation is the congruence relation for the quotient functor generated by
    \begin{equation}
        (x \sim_* y) :\Leftrightarrow \left( \exists z \in A_{k+1}: b^-(z) = x \ \mathrm{and} \ b^+(z) = y \right).
    \end{equation}
    
    It is clear that the generating relation $\sim_*$ implies $\approx_*$ because of the inclusion $A_{k+1} \subset B_{k+1}$, thus this is also true of the full congruence relations $\sim$ and $\approx$. What is left to be shown then is that $\approx_*$ implies $\sim$. Since $\sim$ is already a congruence relation this would mean that the full congruence relation $\approx$ also implies $\sim$.
    
    Up to identifications that implement the axioms of a $(k+1)$-category, the elements of $B_{k+1}$ are lists whose elements are either generating $(k+1)$-cells from $A_k$ or identities of elements of $B_k$, such that neighboring pairs are composable. Suppose $x \approx_* y$ and let $z = (z_1, \dots, z_p)$ be a list which implements the equivalence. We have that
    \begin{equation}
        \pi^- (z) = \pi^- (z_p) = x \qquad \mathrm{and} \qquad \pi^+ (z) = \pi^+ (z_1) = y.
    \end{equation}
    Since neighboring pairs are composable we have that if, for a given $1 \leq i \leq p$, $z_i$ is an identity then $\pi^- (z_i) = \pi^+(z_i)$, and if it is not then $\pi^-(z_i) = b^-(z_i) \sim_* b^+(z_i) = \pi^+(z_i)$. In any case,
    \begin{equation}
        x = \pi^-{z_p} \sim \pi^+ (z_p) = \pi^-(z_{p-1}) \sim \dots \sim \pi^+(z_2) = \pi^-(z_1) \sim \pi^+(z_1) = y.
    \end{equation}
    
    In the case $n=\omega$, define $C^n := I_\omega F_n V_n A$ for any $n \in \mathbb{N}$, and $C^\omega := F_\omega A$. Then
    \begin{align}
        L_k C^{\omega} &= L_k \mathsf{colim} \left(C^{k+1} \to C^{k+2} \to C^{k+3} \to \dots \right) \notag\\ 
        &= \mathsf{colim} \left( L_k F_{k+1} V_{k+1} A \to L_k F_{k+2} V_{k+2} A  \to \dots \right) \notag\\
        &= \mathsf{colim} \left( Q_k \hat{F}_{k+1} V_{k+1} A \to Q_k \hat{F}_{k+1} V_{k+1} A  \to \dots \right) \notag\\
        &= Q_k \hat{F}_{k+1} V_{k+1} A,
    \end{align}
    where the second equality is because $L_k$ is a left adjoint, the third equality is given by the finite-$n$ proof from before, and the fourth equality is because the colimit is of a constant.
\end{proof}

\begin{remark}
    The preceding proposition implies that an $\omega$-computad provides a sequence of generators-and-relations presentations, where the $n$-th term in the sequence is an $n$-category, and such that we can descend down the sequence by using the localization functors $L$. This is the exact form that we want for bistellar data, which is why it will come from an $\omega$-computad. 
\end{remark}

\subsection{Atoms and Molecules}\label{ssec:atoms_molecules}

Below we summarize certain definitions and results from \cite{hadzihasanovic}, as well as introduce small modifications of those, which serve as the basic combinatorial building blocks in our constructions. We also prove \Cref{prop:extend_atoms_to_molecules}, which is the cornerstone technical result that underpins the rest of the constructions of the paper. The presentation is not the same as in \textit{loc.~cit.}, because we do not need all of the technology developed there, and we can use some of the later results proven there to simplify some of the original definitions in the context of the present purposes.

\begin{definition}
    An \emph{oriented graded poset}\footnote{The adjective `oriented' is meant to convey a sense of direction. It is somewhat of a misnomer, because it does not have anything to do with the $\mathrm{SO}$ family of Lie groups. An adjective that would be more appropriate is `directed', as in `directed graph'.} is a graded poset $P \to \mathbb{N}$ together with, for each $x \in P$, a partition of the faces of $x$ into two disjoint sets
    \begin{equation}
        \Delta x = \Delta^- x \sqcup \Delta^+ x,
    \end{equation}
    where we say that the faces in $\Delta^- x$ are the \emph{source} faces, and the ones in $\Delta^+ x$ are the \emph{target} faces. A strict map of oriented graded posets $f: P \to Q$ is a map of graded posets such that
    \begin{equation}
        f(\Delta^\varepsilon x) = \Delta^\varepsilon f(x),
    \end{equation}
    for all $x \in P$ and $\varepsilon \in \{ -, + \}$. These assemble into a category $\ogpos$.
\end{definition}

\begin{lemma}[{\cite[Prop.\,2.2.21]{hadzihasanovic}}]
    The category of oriented graded posets $\ogpos$ has pushouts of inclusions along inclusions.
\end{lemma}

\begin{definition}\label{def:dimension}
    The \emph{dimension} of a graded poset $P$ is the largest degree of an element
    \begin{equation}
        \dim(P) := \max (P \to \mathbb{N}).
    \end{equation}
    By convention, the dimension of the empty poset is $-1$. 
\end{definition}

\begin{remark}
    We inherit this notion of dimension to oriented graded posets as well. Furthermore, we will denote the set of maximal elements of an oriented graded poset $P$ by $\mathsf{M}P$, and those of grade $n$ by $\mathsf{M}P_n$.
\end{remark}

\begin{definition}
    An oriented graded poset $P$ is \emph{pure} if the restriction of the grading map $P \to \mathbb{N}$ to the set $\mathsf{M}P$ of maximal elements of $P$ is constant. When it is, the constant is clearly equal to the dimension.
\end{definition}

\begin{definition}[{\cite[Def.\,2.1.13]{hadzihasanovic}}]\label{def:boundary}
    Oriented graded posets have a negative and a positive boundary, where for each $k \in \mathbb{N}$ and each $\varepsilon \in \{ -,+ \}$, there is an oriented graded subposet
    \begin{equation}
        \partial^\varepsilon_k P = \downarrow\! \{ x \in P_k \mid \nabla^{-\varepsilon}x = \emptyset \} \cup \bigcup_{l < k} \downarrow\! \mathsf{M}P_l,
    \end{equation}
    where $\downarrow$ is the downward closure, and $\nabla^\varepsilon x$ denotes the $\varepsilon$-cofaces of $x$. When we omit the subscript the implicit subscript is $\dim P - 1$. When we omit the superscript we mean the union of both boundaries.
\end{definition}

\begin{definition}[{\cite[Def.\,3.3.2]{hadzihasanovic}}]
    The categories of \emph{atoms}, \emph{round molecules} and \emph{molecules} are the smallest full subcategories of $\ogpos$ that fit in a chain of full subcategories
    \begin{equation}
        \atm[b]{} \ \ \subseteq \ \ \mol[b]{}^{\mathrm{rd}} \ \ \subseteq \ \ \mol[b]{} \ \ \subseteq \ \ \ogpos
    \end{equation}
    and satisfy the following properties:
    \begin{enumerate}
        \item The oriented graded poset $* \xrightarrow{\{ 0 \}} \mathbb{N}$ is in the category of atoms.
        \item A molecule $U$ is a round molecule if for all $k < \dim U$
        \begin{equation}
            \partial^-_k U \cap \partial^+_k U = \partial_{k - 1} U.
        \end{equation}
        \item Given two round molecules $U$ and $V$, such that $\dim U = \dim V$ and $\partial^{\varepsilon} U = \partial^{\varepsilon} V$ for $\varepsilon \in \{ -,+ \}$, then there is an atom $U \Rightarrow V$ whose oriented graded poset is 
        \begin{equation}
            (U \mathop{\amalg}_{\partial U} V)^{\vartriangleright} \ \ \text{ with } \ \ \Delta^- \top = \Delta \top \cap U \text{ and } \Delta^+ \top = \Delta \top \cap V.
        \end{equation}
        \item Given two molecules $U$ and $V$, such that there is an integer $k$ for which $\partial^-_k U = \partial^+_k V$, then there is a molecule $U \#_k V$ whose oriented graded poset is the pushout
        \begin{equation}
             \begin{tikzcd}
                 \partial^+_k U  \arrow[r, Rightarrow, no head] \arrow[d, hook] & \partial^-_k V \arrow[rd, phantom, very near end, "\lrcorner"] \arrow[r, hook] & V \arrow[d] \\
                 U \arrow[rr] & & U \#_k V
             \end{tikzcd}.
         \end{equation}
         In this situation we sometimes say that the two molecules $U$ and $V$ are $k$-gluable.
    \end{enumerate}
\end{definition}

\begin{example}
    We will mostly draw molecules with the possibly degenerate polytopes that they represent, and decorate the cells with arrows to show directionality. In this reading the poset relation is always the face relation of the polytope, the grading is given by the dimension grading, and the arrows go from the source boundary to the target boundary. Some examples are:
    \begin{equation}
        \begin{tikzpicture}[scale=1.5, baseline=0.5]
            \coordinate (A) at (0,0);
            \coordinate (B) at (1,0);
            \coordinate (C) at (2,0);
            \coordinate (D) at (0.25,0.25);

            \foreach \n in {A, B, C}
                \node at (\n) [circle, fill, inner sep=1pt] {};

            \draw[midarrow] (A) .. controls ($(A) + (0.25, 0.35)$) and ($(B) + (-0.25, 0.35)$) .. (B);
            \draw[midarrow] (A) .. controls ($(A) + (0.25, -0.35)$) and ($(B) + (-0.25, -0.35)$) .. (B);
            \draw[midarrow] (B) -- (C);

            \draw[double, double distance=1pt, -{Implies[]}] ({0.51}, {0-0.1}) -- ({0.51}, {0+0.1});
        \end{tikzpicture},
        \qquad
        \begin{tikzpicture}[scale=1, baseline=-0.2cm]

            \def\R{1cm}
            \coordinate (P1) at (162:\R);
            \coordinate (P2) at (90:\R);
            \coordinate (P3) at (18:\R);
            \coordinate (P4) at (-126:\R);
            \coordinate (P5) at (-54:\R);

            \foreach \n in {P1, P2, P3, P4, P5}
                \node at (\n) [circle, fill, inner sep=1pt] {};
            
            \draw[midarrow] (P1) -- (P2);
            \draw[midarrow] (P2) -- (P3);
            \draw[midarrow] (P1) -- (P4);
            \draw[midarrow] (P4) -- (P5);
            \draw[midarrow] (P5) -- (P3);
            
            \draw[double, double distance=1pt, -{Implies[]}] (0,0.3) -- (0,-0.3);
        \end{tikzpicture},
        \qquad
        \begin{tikzpicture}[scale=1, baseline]
            \coordinate (A) at (-1,{-1/(sqrt(3))});
            \coordinate (B) at (1,{-1/(sqrt(3))});
            \coordinate (C) at (0,{2/sqrt(3)});
            \coordinate (D) at (0, 0);

            \foreach \n in {A, B, C, D}
                \node at (\n) [circle, fill, inner sep=1pt] {};

            \draw[midarrow] (A) -- (B);
            \draw[midarrow] (C) -- (B);
            \draw[midarrow] (A) -- (C);
            \draw[midarrow] (A) -- (D);
            \draw[midarrow] (D) -- (C);
            \draw[midarrow] (D) -- (B);

            \draw[double, double distance=1pt, -{Implies[]}] ({-0.25+0.1}, {0.25-0.1}) -- ({-0.25-0.1}, {0.25+0.1});
            \draw[double, double distance=1pt, -{Implies[]}] ({0.25-0.1}, {0.25-0.1}) -- ({0.25+0.1}, {0.25+0.1});
            \draw[double, double distance=1pt, -{Implies[]}] ({0}, {-0.25-0.2}) -- ({0}, {-0.25+0.15});
        \end{tikzpicture}
        \Rrightarrow
        \begin{tikzpicture}[scale=1, baseline]
            \coordinate (A) at (-1,{-1/(sqrt(3))});
            \coordinate (B) at (1,{-1/(sqrt(3))});
            \coordinate (C) at (0,{2/sqrt(3)});
            \coordinate (D) at (0, 0);

            \foreach \n in {A, B, C}
                \node at (\n) [circle, fill, inner sep=1pt] {};

            \draw[midarrow] (A) -- (B);
            \draw[midarrow] (C) -- (B);
            \draw[midarrow] (A) -- (C);

            \draw[double, double distance=1pt, -{Implies[]}] ({0}, {0-0.3}) -- ({0}, {0+0.3});
        \end{tikzpicture}.
    \end{equation}
    Of the above, the pentagon and the tetrahedron are atoms, because they have a maximum, which also determines their dimensions, which, of course, are 2 and 3, respectively. Since they are atoms they are also pure. The first example, which is a globe with an attached line on one end, is a molecule of dimension 2, but not an atom, since it has two maximal elements. It is also not pure, because these maximal elements have different dimensions.
    
    Furthermore, consider the molecule that is glued as
    \begin{equation}
        \begin{tikzpicture}[scale=1, baseline]
            \coordinate (A) at (0,0-0.5);
            \coordinate (B) at (1,0-0.5);
            \coordinate (C) at (1,1-0.5);
            \coordinate (D) at (0,1-0.5);

            \foreach \n in {A, B, C, D}
                \node at (\n) [circle, fill, inner sep=1pt] {};

            \draw[midarrow] (B) -- (A);
            \draw[midarrow] (B) -- (C);
            \draw[midarrow] (C) -- (D);
            \draw[midarrow] (A) -- (D);
            \draw[midarrow] (A) -- (C);

            \draw[double, double distance=1pt, -{Implies[]}] ({0.25-0.1}, {0.75-0.1-0.5}) -- ({0.25+0.1}, {0.75+0.1-0.5});

            \draw[double, double distance=1pt, -{Implies[]}] ({0.75-0.1}, {0.25-0.1-0.5}) -- ({0.75+0.1}, {0.25+0.1-0.5});
        \end{tikzpicture} \ = \
        \begin{tikzpicture}[scale=1, baseline]
            \coordinate (A) at (0,0-0.5);
            \coordinate (B) at (1,0-0.5);
            \coordinate (C) at (1,1-0.5);
            \coordinate (D) at (0,1-0.5);

            \foreach \n in {A, B, C, D}
                \node at (\n) [circle, fill, inner sep=1pt] {};

            \draw[midarrow] (B) -- (A);
            \draw[midarrow] (B) -- (C);
            \draw[midarrow] (C) -- (D);
            \draw[midarrow] (A) -- (C);

            \draw[double, double distance=1pt, -{Implies[]}] ({0.75-0.1}, {0.25-0.1-0.5}) -- ({0.75+0.1}, {0.25+0.1-0.5});
        \end{tikzpicture} \ \#_1 \ \begin{tikzpicture}[scale=1, baseline]
            \coordinate (A) at (0,0-0.5);
            \coordinate (B) at (1,0-0.5);
            \coordinate (C) at (1,1-0.5);
            \coordinate (D) at (0,1-0.5);

            \foreach \n in {A, B, C, D}
                \node at (\n) [circle, fill, inner sep=1pt] {};

            \draw[midarrow] (B) -- (A);
            \draw[midarrow] (C) -- (D);
            \draw[midarrow] (A) -- (D);
            \draw[midarrow] (A) -- (C);

            \draw[double, double distance=1pt, -{Implies[]}] ({0.25-0.1}, {0.75-0.1-0.5}) -- ({0.25+0.1}, {0.75+0.1-0.5});
        \end{tikzpicture} \ = \ \left( \begin{tikzpicture}[scale=1, baseline]
            \coordinate (C) at (1,1-0.5);
            \coordinate (D) at (0,1-0.5);

            \foreach \n in {C, D}
                \node at (\n) [circle, fill, inner sep=1pt] {};

            \draw[midarrow] (C) -- (D);
        \end{tikzpicture} \ \#_0 \ \begin{tikzpicture}[scale=1, baseline]
            \coordinate (A) at (0,0-0.5);
            \coordinate (B) at (1,0-0.5);
            \coordinate (C) at (1,1-0.5);

            \foreach \n in {A, B, C}
                \node at (\n) [circle, fill, inner sep=1pt] {};

            \draw[midarrow] (B) -- (A);
            \draw[midarrow] (B) -- (C);
            \draw[midarrow] (A) -- (C);

            \draw[double, double distance=1pt, -{Implies[]}] ({0.75-0.1}, {0.25-0.1-0.5}) -- ({0.75+0.1}, {0.25+0.1-0.5});
        \end{tikzpicture} \right) \ \#_1 \ \left( \begin{tikzpicture}[scale=1, baseline]
            \coordinate (A) at (0,0-0.5);
            \coordinate (C) at (1,1-0.5);
            \coordinate (D) at (0,1-0.5);

            \foreach \n in {A, C, D}
                \node at (\n) [circle, fill, inner sep=1pt] {};

            \draw[midarrow] (C) -- (D);
            \draw[midarrow] (A) -- (D);
            \draw[midarrow] (A) -- (C);

            \draw[double, double distance=1pt, -{Implies[]}] ({0.25-0.1}, {0.75-0.1-0.5}) -- ({0.25+0.1}, {0.75+0.1-0.5});
        \end{tikzpicture} \ \#_0 \ \begin{tikzpicture}[scale=1, baseline]
            \coordinate (A) at (0,0-0.5);
            \coordinate (B) at (1,0-0.5);

            \foreach \n in {A, B}
                \node at (\n) [circle, fill, inner sep=1pt] {};

            \draw[midarrow] (B) -- (A);
        \end{tikzpicture} \right).
    \end{equation}
    This is, in fact, the only possible way to glue the left-hand side molecule starting from atomic pieces. This shows us that gluing molecules that are not pure can still result in a pure molecule, and for certain molecules this is the only way to get them.

    Subsequently, it also makes clear why we want to glue molecules over boundary pieces of arbitrary codimension relative to the interior of a gluing piece; in situations where the gluing pieces are of different dimensions the concept of codimension is not even well-defined.
\end{example}

\begin{lemma}[{\cite[Lem.\,3.3.8]{hadzihasanovic}}]\label{lem:mol_is_globular}
    If $U$ is an $n$-molecule, for each $k \leq n$ and $\varepsilon \in \{ -,+ \}$, $\partial_k^\varepsilon U$ is a $k$-molecule. Moreover, the construction $\partial^\varepsilon$ is globular for molecules, so that
    \begin{equation}
        \partial^\varepsilon_{n-2} \partial_{n-1}^+ U = \partial^\varepsilon_{n-2} \partial_{n-1}^- U = \partial^\varepsilon_{n-2} U.
    \end{equation}
\end{lemma}

Molecules are such that they have an internal concept of directionality. To make the concept precise we have the \emph{flow graph} of a molecule \cite[Def.\,4.3.1]{hadzihasanovic}. For our purposes we will only need the flow graph of the maximal elements.

\begin{definition}[{\cite[Def.\,4.3.2]{hadzihasanovic}}]
    Let $k$ be an integer and $U$ be a molecule. The \emph{maximal $k$-flow graph} $\mathscr{M}_k U$ of $U$ is the graph whose vertices are the maximal elements of $U$ of degree bigger than $k$, and in which there is an edge from $x$ to $y$ if $x \neq y$ and there exists an element $z$ of degree $k$ which is both in the positive $k$-boundary of $\downarrow\! x$ and in the negative $k$-boundary of $\downarrow\! y$.
\end{definition}

\begin{example}
    Consider the molecule
    \begin{equation}
        U := \begin{tikzpicture}[scale=1, baseline]
            \coordinate (A) at (-1,{-1/(sqrt(3))});
            \coordinate (B) at (1,{-1/(sqrt(3))});
            \coordinate (C) at (0,{2/sqrt(3)});
            \coordinate (D) at (0, 0);
            \coordinate (E) at (3, {-1/(sqrt(3))});

            \foreach \n in {A, B, C, D, E}
                \node at (\n) [circle, fill, inner sep=1pt] {};

            \draw[midarrow] (A) -- (B);
            \draw[midarrow] (C) -- (B);
            \draw[midarrow] (A) -- (C);
            \draw[midarrow] (A) -- (D);
            \draw[midarrow] (C) -- (D);
            \draw[midarrow] (D) -- (B);
            \draw[midarrow] (B) -- node[above]{$d$} (E);

            \draw[double, double distance=1pt, -{Implies[]}] ({-0.25-0.1}, {0.25+0.1}) -- ({-0.25+0.1}, {0.25-0.1});
            \node at ({-0.25-0.1}, {0.25-0.21}) {$a$};
            \draw[double, double distance=1pt, -{Implies[]}] ({0.25+0.1}, {0.25+0.1}) -- ({0.25-0.1}, {0.25-0.1});
            \node at ({0.25+0.11}, {0.25-0.21}) {$b$};
            \draw[double, double distance=1pt, -{Implies[]}] ({0}, {-0.25+0.14}) -- ({0}, {-0.25-0.2});
            \node at ({0-0.22}, {-0.25-0.09}) {$c$};
        \end{tikzpicture},
    \end{equation}
    with the maximal elements called $a$, $b$, $c$ and $d$, as in the picture. The two nonempty maximal flow graphs for this molecule are
    \begin{equation}
        \mathscr{M}_0 U = \begin{tikzcd}
            a & b \arrow[d] \\
            c \arrow[r] & d
        \end{tikzcd} \qquad \qquad
        \mathscr{M}_1 U = \begin{tikzcd}
            a \arrow[d] & b \arrow[l] \arrow[ld] \\
            c
        \end{tikzcd}.
    \end{equation}
    In the maximal 0-flow graph both directed edges are provided by the existence of the vertex where the line $\downarrow d$ attaches to the rest of the molecule. In the maximal 1-flow graph the arrows are given by the edges that are `internal' in the big triangular piece.
\end{example}

We now introduce two different integers that can be assigned to any molecule, that are very helpful in technically dealing with molecules. We will see their usage most prominently in \Cref{prop:layering_existence}.

\begin{definition}[{\cite[Def.\,4.4.1]{hadzihasanovic}}]
    Let $U$ be a molecule. The \emph{frame dimension} $\mathrm{frdim} U$ of $U$ is
    \begin{equation}
        \mathrm{frdim} U := \dim \bigcup \left\{ \downarrow\! x \, \cap \downarrow\! y \mid x, y \in \mathsf{M}U, x \neq y \right\}.
    \end{equation}
\end{definition}

\begin{definition}[{\cite[Def.\,4.2.7]{hadzihasanovic}}]
    Let $U$ be a molecule. The \emph{layering dimension $\mathrm{lydim}U$} of $U$ is the minimal degree $k$, for which there is at most one maximal element of $U$ with degree bigger than $k+1$. By convention, when $U$ is an atom then $\mathrm{lydim} U = -1$. 
\end{definition}

Up to this point we have not imposed any acyclicity conditions on molecules. However, when specifically considering atoms, we are dealing with objects that topologically represent discs, and, since we are motivated by a relation to framing, these should satisfy some form of acyclicity. Here we will work with the weakest possible form of acyclicity that is known to be enough for proving key propositions, but it is, \emph{a priori}, an open question whether geometric considerations impose even more stringent acyclicity conditions.

\begin{definition}[{\cite[Def.\,8.1.1]{hadzihasanovic}}]
    A molecule $U$ is \emph{frame-acyclic} if for any submolecule\footnote{See \cite[]{hadzihasanovic} for more intricate details of the definition of submolecule.} $V$ of $U$, the directed graph $\mathscr{M}_{\mathrm{frdim}V} V$ is acyclic.
\end{definition}

In low dimensions the condition of frame-acyclicity has no bearing because it is always satisfied. More precisely all molecules of dimension $\leq 3$ are automatically frame-acyclic \cite[Thm.\,8.4.11]{hadzihasanovic}, even more though:

\begin{lemma}\label{lem:fa_less5_automatic}
    All atoms of dimension less than or equal to 4 are frame-acyclic.
\end{lemma}

\begin{proof}
    The only dimension to check is dimension 4. Let $U$ is a 4-atom and $W$ be a submolecule of $U$. If $W$ contains the top element of $U$ then by the closure property $W = U$. In that case $\mathrm{frdim} W = -1$ and $\mathscr{M}_{-1} W$ is the directed graph with one element and no edges, which is clearly acyclic. If $W$ does not contain the top element of $U$ then $\dim W < \dim U = 4$, and so $\mathscr{M}_k W$ is acyclic for dimensional reasons.
\end{proof}

\begin{definition}[{\cite[Def.\,4.2.1]{hadzihasanovic}}]
    Let $U$ be a molecule, $-1 \leq k < \dim U$, and $m$ be the number of maximal elements of $U$ of degree more than $k$. A \emph{$k$-layering} of $U$ is a sequence of molecules $(U_i)_{i=1}^m$ such that
    \begin{equation}
        U = U_1 \#_k U_2 \#_k \dots \#_k U_m,
    \end{equation}
    and $\dim U_i > k$ for all $i \in \{1, \dots, m\}$.
\end{definition}

\begin{remark}
    The idea of a $k$-layering is that it is a decomposition of the molecule, in which there is always one maximal element per constituent, as long as that maximal element is of degree larger than $k$. This does not mean that the constituents are atoms, because they can still contain more than one maximal element, as is the case for certain layerings of whisker pasting diagrams in a higher category.
\end{remark}

\begin{proposition}\label{prop:layering_existence}
    Every molecule $U$ has a $k$-layering for some $k$ such that
    \begin{equation}
        \mathrm{frdim} U \leq k \leq \mathrm{lydim} U < \dim U.
    \end{equation}
    If $U$ has a $k$-layering then it also has a $(k+1)$-layering, as long as $k+1 < \dim U$. Moreover, if $U$ is frame-acyclic then it has a $k$-layering for every $k$ in the range $\mathrm{frdim} U \leq k < \dim U$.
\end{proposition}

\begin{proof}
    These are just the results \cite[Lem.\,4.2.6, Thm.\,4.2.12, Cor.\,4.4.6, Thm.\,8.1.4]{hadzihasanovic}.
\end{proof}

\begin{definition}\label{def:simplicial_atoms}
    A \emph{simplicial molecule} is a molecule whose underlying poset is the order poset of a simplicial complex. Similarly, we get \emph{simplicial atoms} and \emph{framed-acyclic simplicial atoms} so that there is a collection of pullback diagrams
    \begin{equation}
        \begin{tikzcd}
            \atm[f]{} \arrow[rd, phantom, very near start, "\ulcorner"]\arrow[r, hook] \arrow[d] & \atm{} \arrow[rd, phantom, very near start, "\ulcorner"] \arrow[r, hook] \arrow[d] & \mol{} \arrow[rd, phantom, very near start, "\ulcorner"] \arrow[rr] \arrow[d] & {} & \mathsf{s} \mathscr{C}\mathsf{px}^{\hookrightarrow} \arrow[d] \\
            \atm[b]{}^{\mathrm{fa}} \arrow[r, hook] & \atm[b]{} \arrow[r, hook] & \mol[b]{} \arrow[r] & \ogpos \arrow[r] & \mathsf{g}\mathscr{P}\mathsf{os}
        \end{tikzcd}.
    \end{equation}
\end{definition}

\begin{remark}
    The restriction to only the injective maps of simplicial complexes $\mathsf{s} \mathscr{C}\mathsf{px}^{\hookrightarrow}$ is necessary to get a forgetful functor to graded posets, since noninjective maps can decrease the grade of an element. This will not make a big difference to us, because we will mostly care about the objects of the categories being defined, which is also reflected in our choice of strict maps between oriented graded posets.
\end{remark}

\begin{remark}
    We will often want to restrict our attention to atoms or molecules of a given dimension, or up to a given dimension. We denote these by adding the subscript $n$ or $\leq n$ to the names of the categories.
\end{remark}

\begin{lemma}\label{lem:boundary_map_restricts}
    Let $U$ be an $n$-molecule, $k \leq n$, and $\varepsilon \in \{ -,+ \}$. If $U$ is simplicial, then $\partial_k^\varepsilon U$ is as well, and if $U$ is frame-acyclic then $\partial_k^\varepsilon U$ is as well.
\end{lemma}

\begin{proof}
    For the first part of the lemma observe that the boundary is a finite union of downward-closures inside a simplicial complex. In other words it is a finite union of simplicial subcomplexes and this is again a simplicial complex. For the second part, $\partial_k^\varepsilon U$ is a submolecule of $U$, so it is also frame-acyclic.
\end{proof}

\begin{lemma}
    If $W$ is a simplicial molecule and $W = U \#_k V$ for some integer $k$ and $k$-gluable molecules $U$ and $V$, then both $U$ and $V$ are also simplicial.
\end{lemma}

\begin{proof}
    Both  $U$ and $V$ are downward-closed subsets in the underlying poset of $W$. Downward-closed subsets of the underlying poset of a simplicial complex are simplicial subcomplexes.
\end{proof}

\begin{remark}
    The previous two lemmata, together with the fact that frame-acyclicity is defined iteratively on submolecules, allows us to make arguments involving boundaries and submolecules of frame-acyclic simplicial molecules, while remaining in the same collection; all of them are guaranteed to be frame-acyclic and simplicial again.
\end{remark}

\begin{lemma}\label{lem:simp_atoms_are_simplices}
    Every simplicial $n$-atom's underlying simplicial complex is isomorphic to $\Delta^n$.
\end{lemma}

\begin{proof}
    By construction a simplicial $n$-atom has an underlying finite poset with a maximum $m$ in it. Let the underlying simplicial complex be $(V, S)$, where $V$ is the set of vertices and $S \subseteq \mathscr{P}(V)$ the subposet of simplices ordered by inclusion. Since $S$ has a maximum
    \begin{equation}
        m = \bigcup_{s \in S} s = V.
    \end{equation}
    But since a simplicial complex must have all subsets of $m$ in its poset $S$
    \begin{equation}
        S \supseteq \downarrow_{\mathscr{P}(V)} m = \mathscr{P}(V) \quad \Rightarrow \quad S = \mathscr{P}(V).
    \end{equation}
    Since $S$ is finite it must be the case that $V \cong [k]$ for some integer $k$, and since the dimension of $m$ in the poset grading must be $n$ we have that $k=n$.
\end{proof}

\begin{remark}
    In particular, the above lemma implies that, for each $n \in \mathbb{N}$, there are only finitely many (isomorphism classes of) simplicial $n$-atoms. This is in contrast to the case of all atoms, where, starting from dimension $n=2$, there are countably infinitely many (isomorphism classes of) $n$-atoms.
\end{remark}

\begin{lemma}[{\cite[Cor.\,3.4.12]{hadzihasanovic}}]
    Given two molecules $U$ and $V$ in $\mol[b]{}$, if there exists an isomorphism $\phi: U \to V$ between them, then it is unique.
\end{lemma}

\begin{remark}
    Because of the above lemma we are somewhat justified in abusing notation and letting the symbols $\atm[b]{n}$, $\mol[b]{n}$, $\atm{n}$ or $\mol{n}$, and all related ones, also denote the sets of isomorphism classes of objects. We will do this consistently in any statement from here on out.
\end{remark}

The intuition we want to have for molecules is that they are, in some sense, free composites of their atomic building blocks. This intuition remains even when we restrict to simplicial molecules with the caveat that now not all composites are allowed. The following proposition makes this statement concrete, and is the technical backbone of our main constructions.

\begin{proposition}\label{prop:extend_atoms_to_molecules}
    Let $n \geq 1$ be an integer, and $C$ be a strict $n$-category. If there exists a commutative diagram of sets
    \begin{equation}
        \begin{tikzcd}
            \atm[f]{n} \arrow[r, "\partial^\varepsilon_{n-1}"] \arrow[d, "\alpha_n"'] & \mol[f]{\leq n - 1} \arrow[d, "\tilde{\alpha}_{n-1}"] \\
            C_n \arrow[r, "\pi^\varepsilon"] & C_{n-1}
        \end{tikzcd},
    \end{equation}
    where for every frame-acyclic simplicial molecule $W$ of dimension less than $n$, and a splitting of that molecule $W = U \#_l V$ with $l < n - 1$, the map $\tilde{\alpha}_{n-1}$ satisfies
    \begin{equation}\label{eq:compositionality_of_extension}
        \tilde{\alpha}_{n-1}(W) = \tilde{\alpha}_{n-1}(U) \circ_l \tilde{\alpha}_{n-1}(V),
    \end{equation}
    then the map $\alpha_n$ can be extended to a map $\tilde{\alpha}_n$, which fits in the commutative diagram of sets
    \begin{equation}
        \begin{tikzcd}
            \mol[f]{\leq n} \arrow[r, "\partial^\varepsilon_{n-1}"] \arrow[d, "\tilde{\alpha}_n"'] & \mol[f]{\leq n - 1} \arrow[d, "\tilde{\alpha}_{n-1}"] \\
            C_n \arrow[r, "\pi^\varepsilon"] & C_{n-1}
        \end{tikzcd},
    \end{equation}
    and also satisfies \Cref{eq:compositionality_of_extension} with the appropriately extended range of splittings.
\end{proposition}

\begin{proof}
    From the given data we can construct another commutative diagram
    \begin{equation}
        \begin{tikzcd}
            \atm[f]{n} \sqcup \mol[f]{\leq n -1} \arrow[r, "\partial^\varepsilon_{n-1}"] \arrow[d, "{(\alpha_n, i_{n-1} \tilde{\alpha}_{n-1})}"'] & \mol[f]{\leq n - 1} \arrow[d, "\tilde{\alpha}_{n-1}"] \\
            C_n \arrow[r, "\pi^\varepsilon"] & C_{n-1}
        \end{tikzcd},
    \end{equation}
    where $\hat{\alpha}_n := (\alpha_n, i_{n-1} \tilde{\alpha}_{n-1})$ continues to satisfy \Cref{eq:compositionality_of_extension}.

    Take the set $\mol[f]{n}$ and make it into a poset $M_n$ by equipping it with the partial order of submolecule inclusion. In this poset the minimal elements are exactly the frame-acyclic simplicial atoms. Indeed, if $U$ is minimal then pick one of its maximal elements $x$, and since the atom $\downarrow x$ is a submolecule of $U$, $U$ must be atomic by minimality. Alternatively if $U$ is atomic then by \cite[Prop.\,5.3.11]{hadzihasanovic} any submolecule inclusion is an isomorphism, which shows that $U$ is minimal in $M_n$.
    
    Furthermore, $M_n$ is well-founded. We can see this because it has a strictly monotone map to $\mathbb{N}$ given by the cardinality of the underlying poset of a molecule. If there were an infinite descending chain in $M_n$ then this would produce a map $\mathbb{N}^{\mathrm{op}} \hookrightarrow M_n \to \mathbb{N}$, which would be a strictly monotone map from the total order $\mathbb{N}^{\mathrm{op}}$ to its opposite, which is a contradiction.

    In defining $\tilde{\alpha}_n$, if we need to evaluate a molecule of dimension less than $n$, or an atom of dimension $n$ then we use $\hat{\alpha}_n$. For the rest we proceed by well-founded induction on the poset $M_n$, for which the base case is already covered. If $U$ is an $n$-dimensional, frame-acyclic simplicial molecule that is not atomic then we set $k := \mathrm{lydim} U \geq 0$, and we pick a $k$-layering $(U_i)_{i=1}^m$ of $U$ so we can define
    \begin{equation}
        \tilde{\alpha}_n (U) := \tilde{\alpha}_n (U_1) \circ_k \tilde{\alpha}_n (U_2) \circ_k \dots \circ_k \tilde{\alpha}_n (U_m),
    \end{equation}
    We're allowed to do this because such a layering must exist by \Cref{prop:layering_existence}, and we can proceed inductively because each constituent of a layering must be a proper submolecule since $U$ is not atomic.

    We have to show that:
    \begin{enumerate}
        \item\label{it:well-defined} the assignment $\tilde{\alpha}_n (U)$ is well-defined,
        \item\label{it:boundaries} the assignment $\tilde{\alpha}_n (U)$ is such that it fits into the required commutative diagram, and
        \item\label{it:decomposition} the assignment $\tilde{\alpha}_n (U)$ satisfies \Cref{eq:compositionality_of_extension} under decomposition.
    \end{enumerate}

    The easiest of these is \Cref{it:boundaries}, and it follows immediately from the observation that the properties of the boundary of a gluing composite are exactly the properties of the source/target of a composite in a category. Formally, this is \cite[Lem.\,3.1.5]{hadzihasanovic} and \cite[Lem.\,3.1.15]{hadzihasanovic}.

    For \Cref{it:decomposition} we argue as follows. Given $W$, with $\mathrm{lydim} W = k'$, and a decomposition $U \#_l V$, as well as, layerings $(U_i)$ and $(V_j)$, both at the respective layering dimensions, the proof of \cite[Thm.\,4.2.12]{hadzihasanovic} explicitly constructs, for certain integers $n_U$ and $n_V$, a $k'$-layering $(W_1, \dots W_{n_U}, W_{n_U+1}, \dots, W_{n_U + n_V})$ such that
    \begin{equation}
        W_1 \#_{k'} \dots \#_{k'} W_{n_U} = \begin{cases}
            U & l = {k'} \\
            U \#_l \partial_{k'}^- V & l < k'
        \end{cases}, \quad
        W_{n_U + 1} \#_{k'} \dots \#_{k'} W_{n_U + n_V} = \begin{cases}
            V & l = {k'} \\
            \partial_{k'}^+ U \#_l V & l < k'
        \end{cases}
    \end{equation}
    In the first case, $l = k'$, we are then done because
    \begin{equation}
        \tilde{\alpha}_n (W) = \tilde{\alpha}_n (W_1) \circ_{l} \dots \circ_{l} \tilde{\alpha}_n (W_{n_U}) \circ_{l} \tilde{\alpha}_n (W_{n_U + 1}) \circ_{l} \dots \circ_{l} \tilde{\alpha}_n (W_{n_U + n_V}) = \tilde{\alpha}_n (U) \circ_{l} \tilde{\alpha}_n (V),
    \end{equation}
    where we used the inductive hypothesis to recombine the pieces for the submolecules $U$ and $V$. In the second case, $l < k'$, the same computation gives
    \begin{equation}
        \tilde{\alpha}_n (W) = \tilde{\alpha}_n (U \#_l \partial_{k'}^- V) \circ_{k'} \tilde{\alpha}_n (\partial_{k'}^+ U \#_l V) = \left( \tilde{\alpha}_n (U) \circ_l \tilde{\alpha}_n (\partial_{k'}^- V) \right) \circ_{k'} \left( \tilde{\alpha}_n (\partial_{k'}^+ U) \circ_l \tilde{\alpha}_n (V) \right),
    \end{equation}
    where the inductive hypothesis again allowed us to decompose the submolecules. By using binary interchange in the category $C$ we have
    \begin{equation}
        \tilde{\alpha}_n (W) = \left( \tilde{\alpha}_n (U) \circ_{k'} \tilde{\alpha}_n (\partial_{k'}^+ U) \right) \circ_{l} \left( \tilde{\alpha}_n (\partial_{k'}^- V) \circ_{k'} \tilde{\alpha}_n (V) \right) = \tilde{\alpha}_n (U) \circ_{l} \tilde{\alpha}_n (V),
    \end{equation}
    where we also used that gluing a boundary onto the same molecule, e.g. $\partial_{k'}^+ V \#_{k'} V $, gives back the molecule $V$, and the inductive hypothesis to allow us to recompose $U$ and $V$ back together.

    To show well-definition --- \Cref{it:well-defined} on the list --- we partially borrow an idea from the proof of \cite[Lem.\,8.2.13]{hadzihasanovic}, and pick two different $k$-layerings of $U$
    \begin{equation}
        U = U_1 \#_k \dots \#_k U_m,  \qquad \qquad U = U''_1 \#_k \dots \#_k U''_m.
    \end{equation}
    A $k$-layering induces a total order on $\mathsf{M}U_{>k}$, the maximal elements of $U$ which have degree $>k$, namely, by seeing which constituent the element is in. Let $\sigma$ then be the unique permutation acting on $\mathsf{M}U_{>k}$ that takes the order induced by $(U_i)$ to the order induced by $(U''_j)$,\footnote{We emphasize that it is not necessarily the case that $U_{i} = U''_{\sigma(i)}$, but they do have the same single maximal element of degree $> k$.} and let $d$ be the number of adjacent transpositions of $\sigma$
    \begin{equation}
        d = \left| \left\{ (i, i') \mid i < i', \sigma(i') < \sigma(i) \right\} \right|
    \end{equation}
    We proceed by induction on $d$. If $d=0$ then the two layerings are the same by \cite[Prop.\,4.3.11]{hadzihasanovic}. When $d > 0$ we will construct a $k$-layering $(U'_l)_{l=0}^m$ such that the unique permutation $\sigma'$, that sends the order induced by $(U'_l)$ to the one induced by $(U''_j)$, is such that
    \begin{equation}
        \sigma = \sigma' \tau_{p, p+1},
    \end{equation}
    where $\tau_{p, p+1}$ is the adjacent transposition on the elements $p$ and $p+1$. If the two choices of $k$-layering $(U_i)$ and $(U'_l)$ are such that they provide the same definition of $\tilde{\alpha}_n (U)$ then, because $d' < d$, well-definition is proved by the inductive hypothesis in $d$.

    It is clear that for $l \neq p, p+1$, we can take $U'_l = U_l$. Now consider $U_p \#_k U_{p+1}$. It has two maximal elements of degree $>k$, but no edge between them in the maximal $k$-flow graph $\mathscr{M}_k (U_p \#_k U_{p+1})$, because we know there exist two $k$-layerings, $(U_i)$ and $(U''_j)$, where the maximal elements are differently ordered. This implies that $r:= \mathrm{frdim} (U_p \#_k U_{p+1}) < k$. Because of the frame-acyclicity assumption on $U$, the submolecules $U_p$, $U_{p+1}$, and $U_p \#_k U_{p+1}$ are all also frame-acyclic, and, since $\mathrm{frdim} U_p, \mathrm{frdim} U_{p+1} \leq r$, there are $r$-layerings
    \begin{equation}
        U_p = W_1 \#_r \dots \#_r W_{n_p - 1}, \qquad \qquad U_{p+1} = \tilde{W}_2 \#_r \dots \#_r \tilde{W}_{n_{p+1}},
    \end{equation}
    for some integers $n_p, n_{p+1} > 1$, where we've chosen the ranges in an unorthodox way to indicate that. While not an $r$-layering, we can write
    \begin{equation}
        U_p = W_1 \#_r \dots \#_r W_{n_p - 1} \#_r W_{n_p}, \qquad \qquad U_{p+1} = \tilde{W_1} \#_r \tilde{W}_2 \#_r \dots \#_r \tilde{W}_{n_{p+1}},
    \end{equation}
    with $W_{n_p} := \partial_r^+ W_{n_p - 1}$, and $\tilde{W}_1 := \partial_r^- \tilde{W}_2$. Observe that these are decompositions of $U_p$ and $U_{p+1}$ such that the maximal elements of degree $> k$ are not in either $W_{n_p}$ or $\tilde{W}_1$, and then, using binary interchange for molecular gluing \cite[Prop.\,3.4.22]{hadzihasanovic}, define
    \begin{equation}
        U_p \#_k U_{p+1} = \left( (W_1 \#_r \dots \#_r W_{n_p - 1}) \#_k \tilde{W}_1 \right) \#_r \left( W_{n_p} \#_k (\tilde{W}_2 \#_r \dots \#_r \tilde{W}_{n_{p+1}}) \right) = : W \#_r \tilde{W},
    \end{equation}
    where $W$ has the maximal element of degree $>k$ of $U_p$ and $\tilde{W}$ the one of $U_{p+1}$. Using this we can finally construct
    \begin{equation}
        U'_p := \partial_k^- W \#_r \tilde{W}, \qquad \qquad U'_{p+1} := W \#_r \partial_k^+ \tilde{W},
    \end{equation}
    where we emphasize the swap of which factor the maximal element is in. $(U'_l)$ is clearly a $k$-layering because there are enough constituents, and the $k$-boundaries agree in the appropriate way for gluing. The only thing left to show is that the two $k$-gluings $(U_i)$ and $(U'_l)$ provide agreeing definitions of $\tilde{\alpha}_n(U)$. This is seen from the following calculation
    \begin{align}
        \tilde{\alpha}_n (U'_p) \circ_k& \tilde{\alpha}_n (U'_{p+1}) = \tilde{\alpha}_n (\partial_k^- W \#_r \tilde{W}) \circ_k \tilde{\alpha}_n (W \#_k \partial_k^+ \tilde{W}) \\
        &= \left( \tilde{\alpha}_n (\partial_k^- W) \circ_r \tilde{\alpha}_n (\tilde{W}) \right) \circ_k \left( \tilde{\alpha}_n (W) \circ_r \tilde{\alpha}_n (\partial_k^+ \tilde{W}) \right) \notag\\
        &= \left(\vphantom{\tilde{W}} \tilde{\alpha}_n (\partial_k^- W) \circ_k \tilde{\alpha}_n (W) \right) \circ_r \left( \tilde{\alpha}_n (\tilde{W}) \circ_k \tilde{\alpha}_n (\partial_k^+ \tilde{W}) \right) \notag\\
        &= \tilde{\alpha}_n (W) \circ_r \tilde{\alpha}_n (\tilde{W}) \notag\\
        &= \tilde{\alpha}_n \left( (W_1 \#_r \dots \#_r W_{n_p - 1}) \#_k \tilde{W}_1 \right) \circ_r \tilde{\alpha}_n \left( W_{n_p} \#_k (\tilde{W}_2 \#_r \dots \#_r \tilde{W}_{n_{p+1}}) \right) \notag\\
        &= \left( \tilde{\alpha}_n (W_1 \#_r \dots \#_r W_{n_p - 1}) \circ_k \tilde{\alpha}_n (\tilde{W}_1) \right) \circ_r \left( \tilde{\alpha}_n (W_{n_p}) \circ_k \tilde{\alpha}_n (\tilde{W}_2 \#_r \dots \#_r \tilde{W}_{n_{p+1}}) \right) \notag\\
        &= \left(\vphantom{\tilde{W}} \tilde{\alpha}_n (W_1 \#_r \dots \#_r W_{n_p - 1}) \circ_r \tilde{\alpha}_n (W_{n_p}) \right) \circ_k \left( \tilde{\alpha}_n (\tilde{W}_1) \circ_r \tilde{\alpha}_n (\tilde{W}_2 \#_r \dots \#_r \tilde{W}_{n_{p+1}}) \right) \notag\\
        &= \tilde{\alpha}_n (U_p) \circ_k \tilde{\alpha}_n (U_{p+1}), \notag
    \end{align}
    where we repeatedly use binary interchange in $C$, and the inductive hypothesis on the proper submolecules $U_p'$, $U_{p+1}'$, $W$ and $\tilde{W}$ and their submolecules.
\end{proof}

\begin{remark}
    It is obvious from the proof that the veracity of the proposition does not actually depend on the fact that the atoms are simplicial; the same proof would go through even for all frame-acyclic molecules. The only effect of being simplicial is that we cannot assume that gluing two simplicial molecules results in a simplicial molecule again, so we have to state the condition of \Cref{eq:compositionality_of_extension} as a decomposition condition rather than as a gluing condition.
\end{remark}

\begin{remark}
    As mentioned before we have the intuition that molecules are built from atoms freely. This result falls short of actually showing such a result, which in fact is also not a result that is known to be true. The most one can currently say is that atoms generate molecules (or, equivalently, that the set of atoms is a basis for the strict $\omega$-category of molecules \cite[Prop.\,5.2.7]{hadzihasanovic}), in the sense that any molecule is a finite gluing of atoms. For our purposes though, the result proven above will be sufficient.
\end{remark}

\section{Bistellar Data}\label{sec:bistellar_definition}

We want to construct the walking bistellar datum in any dimension using frame-acyclic simplicial atoms as generators and relations. This means that functors out of it will be representations of the bistellar datum in some higher category of theories, which could be strict or weak. The reason for using simplicial atoms as the generators and relations for a strict higher category is because we concern ourselves with triangulations and bistellar moves. After constructing the relevant higher categories we will see how they relate to known examples in the literature.

Using the setup from \Cref{ssec:computads} we want to define an $\omega$-computad, which will provide a family of strict $n$-categories, with $n \in \mathbb{N}$, with generators-and-relations presentations for them. This is the correct strategy because a bistellar move is a relation that identifies any two pieces of a triangulation if together they form the boundary of a simplex in one dimension higher. In our language of molecules, the strict localization identifies any two molecules which are the source and target boundary of an atom in one dimension higher. Since \Cref{lem:simp_atoms_are_simplices} tells us that atoms are simplices, this is the exact implementation of bistellar moves that we're looking for.

The construction proceeds inductively in the dimension $n$. At each step we define an $n$-computad together with a map of sets that interprets molecules as certain freely generated composite morphisms in the free category at level $n$. This has the drawback that the definition has a lot of moving parts and is not immediately easy to parse, but the motivation behind it remains the same; we want to make the frame-acyclic simplicial atoms act as the generators of a strict $\omega$-category which we will then strictly localize.

\begin{definition}\label{def:satm_computad}
    The 0-computad $\atm[f]{\leq 0}$ is just the set $\atm[f]{0}$. Given the data of the $n$-computad $\atm[f]{\leq n}$, we can recognize its top cells as\footnote{The map $\alpha_n$ is the component on $n$-cells of the unit of the free-forgetful adjunction $F_n \dashv U_n$ between $n$-computads and strict $n$-categories evaluated at $\atm[f]{n}$.}
    \begin{equation}
        \alpha_n: \atm[f]{n} \xrightarrow{\qquad} F_n (\atm[f]{\leq n})_n,
    \end{equation}
    and, by using \Cref{prop:extend_atoms_to_molecules}, this inductively produces a map
    \begin{equation}
        \mol[f]{\leq n} \xrightarrow{\ \ \tilde{\alpha}_n \ \ } F_n (\atm[f]{\leq n})_n.
    \end{equation}
    The $(n+1)$-computad $\atm[f]{n+1}$ is then defined as the collection consisting of
    \begin{enumerate}
        \item the $n$-computad $\atm[f]{\leq n}$
        \item the set of $(n+1)$-cells $\atm[f]{n+1}$
        \item the source ($\varepsilon=-$) and target ($\varepsilon=+$) maps
        \begin{equation}
            b^{\varepsilon}_n: \atm[f]{n+1} \xrightarrow[\qquad]{\partial^{\varepsilon}_n} \mol[f]{n} \xhookrightarrow{\qquad} \mol[f]{\leq n} \xrightarrow[\qquad]{\tilde{\alpha}_n} F_n (\atm[f]{\leq n})_n.
        \end{equation}
    \end{enumerate}
    The family of computads $\{\atm[f]{\leq n}\}_{n \in \mathbb{N}}$ altogether define the $\omega$-computad $\atm[f]{}$.
\end{definition}

\begin{proof}
    For the preceding construction to be well-defined we need to guarantee that:
    \begin{enumerate}
        \item the assumptions of \Cref{prop:extend_atoms_to_molecules} are satisfied and
        \item for any atom $A$, the morphisms $b^-(A)$ and $b^+(A)$ are parallel.
    \end{enumerate}

    For the first point, since the map $\tilde{\alpha}_{n}$ is defined by induction on $n$, we proceed in this way. When $n=0$ we have that $\tilde{\alpha}_0 = \alpha_0 = \mathrm{id}_{\atm[f]{0}}$, and it vacuously satisfies the gluing property. Having defined $\tilde{\alpha}_{n-1}$ inductively, so that it satisfies the gluing condition of \Cref{prop:extend_atoms_to_molecules}, we have the diagram
    \begin{equation}
        \begin{tikzcd}
            \atm[f]{n} \arrow[r, "\partial^\varepsilon_{n-1}"] \arrow[d, "\alpha_n"'] & \mol[f]{\leq n - 1} \arrow[rd, "\tilde{\alpha}_{n-1}"] \\
            F_n (\atm[f]{\leq n})_n \arrow[r, "\pi^\varepsilon"] & F_n (\atm[f]{\leq n})_{n-1} \arrow[r, Rightarrow, no head] & F_{n-1} (\atm[f]{\leq n-1})_{n-1}
        \end{tikzcd},
    \end{equation}
    where we have already used \Cref{lem:free_cat_same_below} to recognize that the sets of cells in degrees below $n$ are the same as the ones generated by the lower computads. This diagram is commutative because given $U \in \atm[f]{n}$ and $\varepsilon \in \{-,+\}$, we check
    \begin{equation}
        \pi^{\varepsilon} (\alpha_n (U)) := b^{\varepsilon}_{n-1} (U) = \tilde{\alpha}_{n-1} (\partial^{\varepsilon}_{n-1} U),
    \end{equation}
    where the first equality is by the definition of what the source and target maps do to the generators. This satisfies all the assumptions and thus produces $\tilde{\alpha}_n$.

    For the second point, given any $(n+1)$-atom $U \in \atm[f]{n+1}$ and $\varepsilon \in \{-,+\}$, we have
    \begin{equation}
        \pi^- b^{\varepsilon}_n (U) = \pi^- \tilde{\alpha}_n (\partial^{\varepsilon} U) = \tilde{\alpha}_{n-1} (\partial^- \partial^{\varepsilon} U) = \tilde{\alpha}_{n-1} (\partial^+ \partial^{\varepsilon} U) = \pi^+ \tilde{\alpha}_n (\partial^{\varepsilon} U) = \pi^+ b^{\varepsilon}_n (U),
    \end{equation}
    where we use the properties of $\tilde{\alpha}_n$ as guaranteed by \Cref{prop:extend_atoms_to_molecules}, and the globularity for the maps $\partial^{\varepsilon}$ from \Cref{lem:mol_is_globular}.
\end{proof}

\begin{remark}
    For any given low enough dimension one does not need to use \Cref{prop:extend_atoms_to_molecules} to observe that the sources and targets of simplicial atoms can be identified with certain composites in the free category one level lower. In principle, one can just assign the composites by hand by using a graphical calculus. But such a perspective becomes more and more cumbersome and cannot extend to arbitrary dimensions, where an organizational principle becomes necessary.
\end{remark}

\begin{definition}\label{def:bistellar_datum}
    The $\omega$-category
    \begin{equation}
        \mathscr{P} := F_{\omega} (\atm[f]{})
    \end{equation}
    is the walking \emph{$\omega$-bistellar datum}. For any $n \in \mathbb{N}$, the walking \emph{$n$-bistellar datum} is
    \begin{equation}
        \mathscr{P}_n := L_n (\mathscr{P}).
    \end{equation}
\end{definition}

\begin{remark}
    Naively, one might want to circumvent the free generation procedure we used to define bistellar data by simply using molecules as all $n$-cells and building a category out of that. This is not the correct thing to do because of two reasons in particular:
    \begin{enumerate}
        \item While $\mol[b]{}$ can be made into a strict $\omega$-category, the same is not true for $\mol{}$ or $\mol[f]{}$, because not all gluing is allowed; if the result of gluing is not simplicial, for example, then we throw it out.
        \item As remarked on before, it is not known whether even $\mol[b]{}$ is a free strict $\omega$-category; all that we know is that if it is, then $\atm[b]{}$ is an $\omega$-computad, and in that case, it is exactly the computad that generates $\mol[b]{}$.
    \end{enumerate}
\end{remark}

Apart from the fact that the generators are simplices, using the localization also implements bistellar moves. Namely, two local pieces of a simplicial complex are identified when the two pieces together form the boundary of a simplex in one dimension higher. Combinations of bistellar moves are then implemented by the higher gluing of simplices. The following theorems go some way into motivating this further by showing that in low dimensional cases we reproduce data that is expected in the literature.

\begin{theorem}\label{thm:satm_in_low_d}
    The sets of simplicial atoms dimensions 0, 1, 2, and 3 are
    \begin{equation}
        \atm{0} = \left\{
            \begin{tikzpicture}
                \node at (0,0) [circle, fill, inner sep=1pt] {};
            \end{tikzpicture}
        \right\}, \quad
        \atm{1} = \left\{
            \begin{tikzpicture}[scale=0.5, baseline]
                \node at (0,0.2) [circle, fill, inner sep=1pt] {};
                \node at (2,0.2) [circle, fill, inner sep=1pt] {};
                \draw[midarrow] (0,0.2) -- (2,0.2);
            \end{tikzpicture}
        \right\}, \quad
        \atm{2} = \left\{
            \begin{tikzpicture}[scale=0.6, baseline]
                \coordinate (A) at (0,{-sqrt(3)/2*1cm + 0.3cm});
                \coordinate (B) at (2cm,{-sqrt(3)/2*1cm + 0.3cm});
                \coordinate (C) at (1cm,{sqrt(3)/2*1cm + 0.3cm});

                \foreach \n in {A, B, C}
                    \node at (\n) [circle, fill, inner sep=1pt] {};

                \draw[midarrow] (A) -- (B);
                \draw[midarrow] (C) -- (B);
                \draw[midarrow] (A) -- (C);

                \draw[double, double distance=1pt, -{Implies[]}] ({1cm}, {0.05*sqrt(3)*1cm + 0.3cm}) -- ({1cm}, {-0.25*sqrt(3)*1cm + 0.3cm});
            \end{tikzpicture},
            \begin{tikzpicture}[scale=0.6, baseline]
                \coordinate (A) at (0,{-sqrt(3)/2*1cm + 0.3cm});
                \coordinate (B) at (2cm,{-sqrt(3)/2*1cm + 0.3cm});
                \coordinate (C) at (1cm,{sqrt(3)/2*1cm + 0.3cm});

                \foreach \n in {A, B, C}
                    \node at (\n) [circle, fill, inner sep=1pt] {};

                \draw[midarrow] (A) -- (B);
                \draw[midarrow] (C) -- (B);
                \draw[midarrow] (A) -- (C);

                \draw[double, double distance=1pt, -{Implies[]}] ({1cm}, {-0.25*sqrt(3)*1cm + 0.3cm}) -- ({1cm}, {0.05*sqrt(3)*1cm + 0.3cm});
            \end{tikzpicture}
        \right\},
    \end{equation}
    \vspace{-1ex}
    \begin{equation}\notag
        \atm{3} = \left\{
            \begin{tikzpicture}[scale=1, baseline]
                \coordinate (A) at (0,0-0.5);
                \coordinate (B) at (1,0-0.5);
                \coordinate (C) at (1,1-0.5);
                \coordinate (D) at (0,1-0.5);

                \foreach \n in {A, B, C, D}
                    \node at (\n) [circle, fill, inner sep=1pt] {};

                \draw[midarrow] (A) -- (B);
                \draw[midarrow] (B) -- (C);
                \draw[midarrow] (C) -- (D);
                \draw[midarrow] (A) -- (D);
                \draw[midarrow] (A) -- (C);

                \draw[double, double distance=1pt, -{Implies[]}] ({0.25+0.1}, {0.75+0.1-0.5}) -- ({0.25-0.1}, {0.75-0.1-0.5});
                
                \draw[double, double distance=1pt, -{Implies[]}] ({0.75+0.1}, {0.25-0.1-0.5}) -- ({0.75-0.1}, {0.25+0.1-0.5});
            \end{tikzpicture}
            \mathrel{\vcenter{\offinterlineskip\halign{$\hfil#\hfil$\cr\Rrightarrow\cr\Lleftarrow\cr}}}
            \begin{tikzpicture}[scale=1, baseline]
                \coordinate (A) at (0,0-0.5);
                \coordinate (B) at (1,0-0.5);
                \coordinate (C) at (1,1-0.5);
                \coordinate (D) at (0,1-0.5);

                \foreach \n in {A, B, C, D}
                    \node at (\n) [circle, fill, inner sep=1pt] {};

                \draw[midarrow] (A) -- (B);
                \draw[midarrow] (B) -- (C);
                \draw[midarrow] (C) -- (D);
                \draw[midarrow] (A) -- (D);
                \draw[midarrow] (B) -- (D);

                \draw[double, double distance=1pt, -{Implies[]}] ({0.25+0.1}, {0.25-0.1-0.5}) -- ({0.25-0.1}, {0.25+0.1-0.5});

                \draw[double, double distance=1pt, -{Implies[]}] ({0.75+0.1}, {0.75+0.1-0.5}) -- ({0.75-0.1}, {0.75-0.1-0.5});
            \end{tikzpicture},
            \begin{tikzpicture}[scale=1, baseline]
                \coordinate (A) at (0,0-0.5);
                \coordinate (B) at (1,0-0.5);
                \coordinate (C) at (1,1-0.5);
                \coordinate (D) at (0,1-0.5);

                \foreach \n in {A, B, C, D}
                    \node at (\n) [circle, fill, inner sep=1pt] {};

                \draw[midarrow] (A) -- (B);
                \draw[midarrow] (B) -- (C);
                \draw[midarrow] (C) -- (D);
                \draw[midarrow] (A) -- (D);
                \draw[midarrow] (A) -- (C);

                \draw[double, double distance=1pt, -{Implies[]}] ({0.25-0.1}, {0.75-0.1-0.5}) -- ({0.25+0.1}, {0.75+0.1-0.5});
                
                \draw[double, double distance=1pt, -{Implies[]}] ({0.75-0.1}, {0.25+0.1-0.5}) -- ({0.75+0.1}, {0.25-0.1-0.5});
            \end{tikzpicture}
            \mathrel{\vcenter{\offinterlineskip\halign{$\hfil#\hfil$\cr\Rrightarrow\cr\Lleftarrow\cr}}}
            \begin{tikzpicture}[scale=1, baseline]
                \coordinate (A) at (0,0-0.5);
                \coordinate (B) at (1,0-0.5);
                \coordinate (C) at (1,1-0.5);
                \coordinate (D) at (0,1-0.5);

                \foreach \n in {A, B, C, D}
                    \node at (\n) [circle, fill, inner sep=1pt] {};

                \draw[midarrow] (A) -- (B);
                \draw[midarrow] (B) -- (C);
                \draw[midarrow] (C) -- (D);
                \draw[midarrow] (A) -- (D);
                \draw[midarrow] (B) -- (D);

                \draw[double, double distance=1pt, -{Implies[]}] ({0.25-0.1}, {0.25+0.1-0.5}) -- ({0.25+0.1}, {0.25-0.1-0.5});

                \draw[double, double distance=1pt, -{Implies[]}] ({0.75-0.1}, {0.75-0.1-0.5}) -- ({0.75+0.1}, {0.75+0.1-0.5});
            \end{tikzpicture},
            \begin{tikzpicture}[scale=1, baseline]
                \coordinate (A) at (0,0-0.5);
                \coordinate (B) at (1,0-0.5);
                \coordinate (C) at (1,1-0.5);
                \coordinate (D) at (0,1-0.5);

                \foreach \n in {A, B, C, D}
                    \node at (\n) [circle, fill, inner sep=1pt] {};

                \draw[midarrow] (B) -- (A);
                \draw[midarrow] (B) -- (C);
                \draw[midarrow] (C) -- (D);
                \draw[midarrow] (A) -- (D);
                \draw[midarrow] (A) -- (C);

                \draw[double, double distance=1pt, -{Implies[]}] ({0.25-0.1}, {0.75-0.1-0.5}) -- ({0.25+0.1}, {0.75+0.1-0.5});

                \draw[double, double distance=1pt, -{Implies[]}] ({0.75-0.1}, {0.25-0.1-0.5}) -- ({0.75+0.1}, {0.25+0.1-0.5});
            \end{tikzpicture}
            \mathrel{\vcenter{\offinterlineskip\halign{$\hfil#\hfil$\cr\Rrightarrow\cr\Lleftarrow\cr}}}
            \begin{tikzpicture}[scale=1, baseline]
                \coordinate (A) at (0,0-0.5);
                \coordinate (B) at (1,0-0.5);
                \coordinate (C) at (1,1-0.5);
                \coordinate (D) at (0,1-0.5);

                \foreach \n in {A, B, C, D}
                    \node at (\n) [circle, fill, inner sep=1pt] {};

                \draw[midarrow] (B) -- (A);
                \draw[midarrow] (B) -- (C);
                \draw[midarrow] (C) -- (D);
                \draw[midarrow] (A) -- (D);
                \draw[midarrow] (B) -- (D);

                \draw[double, double distance=1pt, -{Implies[]}] ({0.25-0.1}, {0.25-0.1-0.5}) -- ({0.25+0.1}, {0.25+0.1-0.5});
                
                \draw[double, double distance=1pt, -{Implies[]}] ({0.75-0.1}, {0.75-0.1-0.5}) -- ({0.75+0.1}, {0.75+0.1-0.5});
            \end{tikzpicture},
            \begin{tikzpicture}[scale=1, baseline]
                \coordinate (A) at (0,0-0.5);
                \coordinate (B) at (1,0-0.5);
                \coordinate (C) at (1,1-0.5);
                \coordinate (D) at (0,1-0.5);

                \foreach \n in {A, B, C, D}
                    \node at (\n) [circle, fill, inner sep=1pt] {};

                \draw[midarrow] (B) -- (A);
                \draw[midarrow] (B) -- (C);
                \draw[midarrow] (C) -- (D);
                \draw[midarrow] (A) -- (D);
                \draw[midarrow] (B) -- (D);

                \draw[double, double distance=1pt, -{Implies[]}] ({0.25-0.1}, {0.25-0.1-0.5}) -- ({0.25+0.1}, {0.25+0.1-0.5});
                
                \draw[double, double distance=1pt, -{Implies[]}] ({0.75-0.1}, {0.75-0.1-0.5}) -- ({0.75+0.1}, {0.75+0.1-0.5});
            \end{tikzpicture}
            \mathrel{\vcenter{\offinterlineskip\halign{$\hfil#\hfil$\cr\Rrightarrow\cr\Lleftarrow\cr}}}
            \begin{tikzpicture}[scale=1, baseline]
                \coordinate (A) at (0,0-0.5);
                \coordinate (B) at (1,0-0.5);
                \coordinate (C) at (1,1-0.5);
                \coordinate (D) at (0,1-0.5);

                \foreach \n in {A, B, C, D}
                    \node at (\n) [circle, fill, inner sep=1pt] {};

                \draw[midarrow] (B) -- (A);
                \draw[midarrow] (B) -- (C);
                \draw[midarrow] (C) -- (D);
                \draw[midarrow] (A) -- (D);
                \draw[midarrow] (C) -- (A);

                \draw[double, double distance=1pt, -{Implies[]}] ({0.25-0.1}, {0.75-0.1-0.5}) -- ({0.25+0.1}, {0.75+0.1-0.5});
                
                \draw[double, double distance=1pt, -{Implies[]}] ({0.75-0.1}, {0.25-0.1-0.5}) -- ({0.75+0.1}, {0.25+0.1-0.5});
            \end{tikzpicture},
        \right.
    \end{equation}
    \begin{equation}\notag
        \qquad \qquad \qquad \left.
            \begin{tikzpicture}[scale=0.6, baseline]
                \coordinate (A) at (-1,{-1/(sqrt(3))});
                \coordinate (B) at (1,{-1/(sqrt(3))});
                \coordinate (C) at (0,{2/sqrt(3)});
                \coordinate (D) at (0, 0);

                \foreach \n in {A, B, C, D}
                    \node at (\n) [circle, fill, inner sep=1pt] {};

                \draw[midarrow] (A) -- (B);
                \draw[midarrow] (C) -- (B);
                \draw[midarrow] (A) -- (C);
                \draw[midarrow] (A) -- (D);
                \draw[midarrow] (C) -- (D);
                \draw[midarrow] (D) -- (B);

                \draw[double, double distance=1pt, -{Implies[]}] ({-0.25-0.1}, {0.25+0.1}) -- ({-0.25+0.1}, {0.25-0.1});
                
                \draw[double, double distance=1pt, -{Implies[]}] ({0.25+0.1}, {0.25+0.1}) -- ({0.25-0.1}, {0.25-0.1});
                \draw[double, double distance=1pt, -{Implies[]}] ({0}, {-0.25+0.14}) -- ({0}, {-0.25-0.2});
            \end{tikzpicture}
            \mathrel{\vcenter{\offinterlineskip\halign{$\hfil#\hfil$\cr\Rrightarrow\cr\Lleftarrow\cr}}}
            \begin{tikzpicture}[scale=0.6, baseline]
                \coordinate (A) at (-1,{-1/(sqrt(3))});
                \coordinate (B) at (1,{-1/(sqrt(3))});
                \coordinate (C) at (0,{2/sqrt(3)});
                \coordinate (D) at (0, 0);

                \foreach \n in {A, B, C}
                    \node at (\n) [circle, fill, inner sep=1pt] {};

                \draw[midarrow] (A) -- (B);
                \draw[midarrow] (C) -- (B);
                \draw[midarrow] (A) -- (C);

                \draw[double, double distance=1pt, -{Implies[]}] ({0}, {0+0.3}) -- ({0}, {0-0.3});
            \end{tikzpicture},
            \begin{tikzpicture}[scale=0.6, baseline]
                \coordinate (A) at (-1,{-1/(sqrt(3))});
                \coordinate (B) at (1,{-1/(sqrt(3))});
                \coordinate (C) at (0,{2/sqrt(3)});
                \coordinate (D) at (0, 0);

                \foreach \n in {A, B, C, D}
                    \node at (\n) [circle, fill, inner sep=1pt] {};

                \draw[midarrow] (A) -- (B);
                \draw[midarrow] (C) -- (B);
                \draw[midarrow] (A) -- (C);
                \draw[midarrow] (A) -- (D);
                \draw[midarrow] (D) -- (C);
                \draw[midarrow] (D) -- (B);

                \draw[double, double distance=1pt, -{Implies[]}] ({-0.25-0.1}, {0.25+0.1}) -- ({-0.25+0.1}, {0.25-0.1});
                \draw[double, double distance=1pt, -{Implies[]}] ({0.25+0.1}, {0.25+0.1}) -- ({0.25-0.1}, {0.25-0.1});
                \draw[double, double distance=1pt, -{Implies[]}] ({0}, {-0.25+0.14}) -- ({0}, {-0.25-0.2});
            \end{tikzpicture}
            \mathrel{\vcenter{\offinterlineskip\halign{$\hfil#\hfil$\cr\Rrightarrow\cr\Lleftarrow\cr}}}
            \begin{tikzpicture}[scale=0.6, baseline]
                \coordinate (A) at (-1,{-1/(sqrt(3))});
                \coordinate (B) at (1,{-1/(sqrt(3))});
                \coordinate (C) at (0,{2/sqrt(3)});
                \coordinate (D) at (0, 0);

                \foreach \n in {A, B, C}
                    \node at (\n) [circle, fill, inner sep=1pt] {};

                \draw[midarrow] (A) -- (B);
                \draw[midarrow] (C) -- (B);
                \draw[midarrow] (A) -- (C);

                \draw[double, double distance=1pt, -{Implies[]}] ({0}, {0+0.3}) -- ({0}, {0-0.3});
            \end{tikzpicture},
            \begin{tikzpicture}[scale=0.6, baseline]
                \coordinate (A) at (-1,{-1/(sqrt(3))});
                \coordinate (B) at (1,{-1/(sqrt(3))});
                \coordinate (C) at (0,{2/sqrt(3)});
                \coordinate (D) at (0, 0);

                \foreach \n in {A, B, C, D}
                    \node at (\n) [circle, fill, inner sep=1pt] {};

                \draw[midarrow] (A) -- (B);
                \draw[midarrow] (C) -- (B);
                \draw[midarrow] (A) -- (C);
                \draw[midarrow] (A) -- (D);
                \draw[midarrow] (C) -- (D);
                \draw[midarrow] (D) -- (B);

                \draw[double, double distance=1pt, -{Implies[]}] ({-0.25+0.1}, {0.25-0.1}) -- ({-0.25-0.1}, {0.25+0.1});
                \draw[double, double distance=1pt, -{Implies[]}] ({0.25-0.1}, {0.25-0.1}) -- ({0.25+0.1}, {0.25+0.1});
                \draw[double, double distance=1pt, -{Implies[]}] ({0}, {-0.25-0.14}) -- ({0}, {-0.25+0.2});
            \end{tikzpicture}
            \mathrel{\vcenter{\offinterlineskip\halign{$\hfil#\hfil$\cr\Rrightarrow\cr\Lleftarrow\cr}}}
            \begin{tikzpicture}[scale=0.6, baseline]
                \coordinate (A) at (-1,{-1/(sqrt(3))});
                \coordinate (B) at (1,{-1/(sqrt(3))});
                \coordinate (C) at (0,{2/sqrt(3)});
                \coordinate (D) at (0, 0);

                \foreach \n in {A, B, C}
                    \node at (\n) [circle, fill, inner sep=1pt] {};

                \draw[midarrow] (A) -- (B);
                \draw[midarrow] (C) -- (B);
                \draw[midarrow] (A) -- (C);

                \draw[double, double distance=1pt, -{Implies[]}] ({0}, {0-0.3}) -- ({0}, {0+0.3});
            \end{tikzpicture},
            \begin{tikzpicture}[scale=0.6, baseline]
                \coordinate (A) at (-1,{-1/(sqrt(3))});
                \coordinate (B) at (1,{-1/(sqrt(3))});
                \coordinate (C) at (0,{2/sqrt(3)});
                \coordinate (D) at (0, 0);

                \foreach \n in {A, B, C, D}
                    \node at (\n) [circle, fill, inner sep=1pt] {};

                \draw[midarrow] (A) -- (B);
                \draw[midarrow] (C) -- (B);
                \draw[midarrow] (A) -- (C);
                \draw[midarrow] (A) -- (D);
                \draw[midarrow] (D) -- (C);
                \draw[midarrow] (D) -- (B);

                \draw[double, double distance=1pt, -{Implies[]}] ({-0.25+0.1}, {0.25-0.1}) -- ({-0.25-0.1}, {0.25+0.1});
                \draw[double, double distance=1pt, -{Implies[]}] ({0.25-0.1}, {0.25-0.1}) -- ({0.25+0.1}, {0.25+0.1});
                \draw[double, double distance=1pt, -{Implies[]}] ({0}, {-0.25-0.14}) -- ({0}, {-0.25+0.2});
            \end{tikzpicture}
            \mathrel{\vcenter{\offinterlineskip\halign{$\hfil#\hfil$\cr\Rrightarrow\cr\Lleftarrow\cr}}}
            \begin{tikzpicture}[scale=0.6, baseline]
                \coordinate (A) at (-1,{-1/(sqrt(3))});
                \coordinate (B) at (1,{-1/(sqrt(3))});
                \coordinate (C) at (0,{2/sqrt(3)});
                \coordinate (D) at (0, 0);

                \foreach \n in {A, B, C}
                    \node at (\n) [circle, fill, inner sep=1pt] {};

                \draw[midarrow] (A) -- (B);
                \draw[midarrow] (C) -- (B);
                \draw[midarrow] (A) -- (C);

                \draw[double, double distance=1pt, -{Implies[]}] ({0}, {0-0.3}) -- ({0}, {0+0.3});
            \end{tikzpicture} 
        \right\}
      \end{equation}
\end{theorem}

\begin{proof}
    As is already apparent from the way we wrote the set of atoms in 3 dimensions, the strategy is to recognize the boundary pieces of the atom and fill in the top cell afterward. By \Cref{lem:simp_atoms_are_simplices} the shape of an atom is always a simplex, so we only need to try finitely many combinations.
    
    The $n=0$ case is definitional. In $n=1$ there is only one choice for the top cell, because the source and target boundaries have to be the same. For $n=2$, the $2$-simplex has three $1$-faces, so there is only one possible split into two non-empty pieces, namely one boundary gets two $1$-simplices and the other gets one. There is a unique molecule structure on both boundary pieces so that we get two simplicial $2$-atoms, depending on which boundary piece is source and which one is target.

    For $n=3$ there are two possible ways to split the boundary of the $3$-simplex into two non-empty parts, the 2-2 and 1-3 bistellar moves, as they are known in the literature. The 8 possible molecular structures on the 2-2 pieces are displayed in the first row of the list of simplicial atoms in 3 dimensions.  For the 1-3 moves, given a description of the 3 side, the 1 side is fully fixed by the requirement that both sides share both their source and target boundary. The 8 possibilities are the ones displayed in the second row of the list of simplicial atoms in 3 dimensions.
\end{proof}

The above \Cref{thm:satm_in_low_d} and \Cref{prop:localization_is_presentation} do the lion's share of the work towards writing down a generators-and-relations presentation of bistellar data in low dimensions.

\begin{corollary}\label{cor:1_bistellar_datum}
    The walking $1$-bistellar datum $\mathscr{P}_1$ has generators
    \begin{equation}
        a, \qquad e:a \to a
    \end{equation}
    and satisfies the relation
    \begin{equation}
        e \circ e = e.
    \end{equation}
\end{corollary}

\begin{proof}
    Since $\atm{0}{} = \{ \begin{tikzpicture}[scale=0.5, baseline=-0.1cm] \node at (0,0) [circle, fill, inner sep=1pt] {}; \end{tikzpicture} \}$ and $\atm{1}{} = \{ \begin{tikzpicture}[scale=0.5, baseline] \node at (0,0.2) [circle, fill, inner sep=1pt] {}; \node at (2,0.2) [circle, fill, inner sep=1pt] {};\draw[midarrow] (0,0.2) -- (2,0.2); \end{tikzpicture}\}$ these are the two mentioned generators. For the relations we have $\atm{2}{}$. At this level we have two generators, as seen in \Cref{thm:satm_in_low_d}, but both impose the same relation,
    \begin{equation}
        \tilde{\alpha}_1 \left( \begin{tikzpicture}[scale=0.5, baseline]
            \node at (0,0.2) [circle, fill, inner sep=1pt] {};
            \node at (2,0.2) [circle, fill, inner sep=1pt] {};
            \node at (4,0.2) [circle, fill, inner sep=1pt] {};

            \draw[midarrow] (0,0.2) -- (2,0.2);
            \draw[midarrow] (2,0.2) -- (4,0.2);
        \end{tikzpicture} \right) \quad = \quad
        \alpha_1 \left( \begin{tikzpicture}[scale=0.5, baseline]
            \node at (0,0.2) [circle, fill, inner sep=1pt] {};
            \node at (2,0.2) [circle, fill, inner sep=1pt] {};

            \draw[midarrow] (0,0.2) -- (2,0.2);
        \end{tikzpicture} \right),
    \end{equation}
    because the boundary maps $\partial^\varepsilon$, which provide the objects being identified, swap sign $\varepsilon \mapsto - \varepsilon$ under the permutation of the two generators. The fact that the gluing of two copies of the one-dimensional generator is their composition in the category $\mathscr{P}_1$ is exactly recognized by the construction in \Cref{prop:extend_atoms_to_molecules}. More explicitly, this is because
    \begin{equation}
        \tilde{\alpha}_1 \left( \begin{tikzpicture}[scale=0.5, baseline]
            \node at (0,0.2) [circle, fill, inner sep=1pt] {};
            \node at (2,0.2) [circle, fill, inner sep=1pt] {};
            \node at (4,0.2) [circle, fill, inner sep=1pt] {};

            \draw[midarrow] (0,0.2) -- (2,0.2);
            \draw[midarrow] (2,0.2) -- (4,0.2);
        \end{tikzpicture} \right) \quad = \quad
        \tilde{\alpha}_1 \left( \begin{tikzpicture}[scale=0.5, baseline]
            \node at (0,0.2) [circle, fill, inner sep=1pt] {};
            \node at (2,0.2) [circle, fill, inner sep=1pt] {};

            \draw[midarrow] (0,0.2) -- (2,0.2);
        \end{tikzpicture} \ \ \#_0 \ \ \begin{tikzpicture}[scale=0.5, baseline]
            \node at (0,0.2) [circle, fill, inner sep=1pt] {};
            \node at (2,0.2) [circle, fill, inner sep=1pt] {};

            \draw[midarrow] (0,0.2) -- (2,0.2);
        \end{tikzpicture} \right) \quad = \quad
        \alpha_1 \left( \begin{tikzpicture}[scale=0.5, baseline]
            \node at (0,0.2) [circle, fill, inner sep=1pt] {};
            \node at (2,0.2) [circle, fill, inner sep=1pt] {};

            \draw[midarrow] (0,0.2) -- (2,0.2);
        \end{tikzpicture} \right) \ \circ \ 
        \alpha_1 \left( \begin{tikzpicture}[scale=0.5, baseline]
            \node at (0,0.2) [circle, fill, inner sep=1pt] {};
            \node at (2,0.2) [circle, fill, inner sep=1pt] {};

            \draw[midarrow] (0,0.2) -- (2,0.2);
        \end{tikzpicture} \right).
    \end{equation}
\end{proof}

\begin{corollary}\label{cor:2_bistellar_datum}
    The walking $2$-bistellar datum $\mathscr{P}_2$ has generators
    \begin{equation}
        a, \quad A: a \to a, \quad \mu: A \otimes A \to A, \quad \delta: A \to A \otimes A,
    \end{equation}
    where $\mu$ and $\delta$ are called multiplication and comultiplication, respectively, and they satisfy the relations, which we provide with string diagrams:\footnote{Although it will not play a big role because of symmetry, we remind the reader of our conventions for string diagrams, as seen in the notation subsection of the introduction, \Cref{ssec:notation}.}
    \begin{enumerate}
        \item associativity for $\mu$
        \begin{equation}
            \begin{tikzpicture}[very thick, scale=0.6, color=green!50!black, baseline]
                \tktrival{1}{1}{0}{-1}{1}{1}
                \draw (0.5,-1) -- (0.5,0);
                \tktrival{1}{0}{0}{1}{1}{1}
            \end{tikzpicture} \quad = \quad
            \begin{tikzpicture}[very thick, scale=0.6, color=green!50!black, baseline]
                \draw (-0.5,-1) -- (-0.5,0);
                \tktrival{1}{1}{1}{-1}{1}{1}
                \tktrival{1}{0}{0}{1}{1}{1}
            \end{tikzpicture},
        \end{equation}
        \item coassociativity for $\delta$
        \begin{equation}
            \begin{tikzpicture}[very thick, scale=0.6, color=green!50!black, baseline]
                \tktrival{-1}{-1}{-1}{1}{1}{1}
                \draw (0.5,1) -- (0.5,0);
                \tktrival{-1}{0}{0}{-1}{1}{1}
            \end{tikzpicture} \quad = \quad
            \begin{tikzpicture}[very thick, scale=0.6, color=green!50!black, baseline]
                \tktrival{-1}{-1}{-1}{0}{1}{1}
                \draw (-1,1) -- (-1,0);
                \tktrival{-1}{0}{0}{0}{1}{1}
            \end{tikzpicture},
        \end{equation}
        \item Frobenius relations relating $\mu$ and $\delta$
        \begin{equation}
            \begin{tikzpicture}[very thick, scale=0.6, color=green!50!black, baseline]
                \draw (-1,-1) -- (-1,0);
                \tktrival{-1}{0}{0.5}{-1}{1}{1}
                \tktrival{1}{-1}{-1}{0}{1}{1}
                \draw (1,0) -- (1,1);
            \end{tikzpicture} \quad = \quad
            \begin{tikzpicture}[very thick, scale=0.6, color=green!50!black, baseline]
                \tktrival{1}{0}{0}{0}{1}{1}
                \tktrival{-1}{0}{0}{0}{1}{1}
            \end{tikzpicture} \quad = \quad
            \begin{tikzpicture}[very thick, scale=0.6, color=green!50!black, baseline]
                \tktrival{-1}{0}{-1}{-1}{1}{1}
                \draw (0.5,-1) -- (0.5,0);
                \draw (-1.5, 0) -- (-1.5, 1);
                \tktrival{1}{0}{0}{1}{1}{1}
            \end{tikzpicture},
        \end{equation}
        and
        \item bubble omissions
        \begin{equation}
            \begin{tikzpicture}[very thick, scale=0.7 * 2/3, color=green!50!black, baseline]
                \tktrival{1}{0}{0}{0}{1}{1}
                \tktrival{-1}{0}{0}{0}{1}{1}
                \tktrival{1}{0}{0}{2}{1}{1}
            \end{tikzpicture} \quad = \quad
            \begin{tikzpicture}[very thick, scale=0.7 * 2/3, color=green!50!black, baseline]
                \tktrival{1}{0}{0}{0}{1}{1}
                \draw (0,0) -- (0,2);
            \end{tikzpicture}, \qquad \qquad
            \begin{tikzpicture}[very thick, scale=0.7 * 2/3, color=green!50!black, baseline]
                \tktrival{-1}{0}{0}{-1}{1}{1}
                \tktrival{1}{0}{0}{1}{1}{1}
                \tktrival{-1}{0}{0}{1}{1}{1}
            \end{tikzpicture} \quad = \quad
            \begin{tikzpicture}[very thick, scale=0.7 * 2/3, color=green!50!black, baseline]
                \draw (0,-1) -- (0,1);
                \tktrival{-1}{0}{0}{1}{1}{1}
            \end{tikzpicture}.
        \end{equation}
    \end{enumerate}
\end{corollary}

\begin{proof}
    The generators are read off from the sets $\atm{i}{}$ for $i \in \{0, 1, 2\}$ of \Cref{thm:satm_in_low_d}. The relations come from the set $\atm{3}{}$, from the same theorem. To translate the simplicial diagrams into string diagrams we need to take the Poincar\'e dual. To induce sources and targets on the dual there is in principle a choice, but having fixed our conventions for how we write string diagrams this choice is automatically fixed. For example, for the negative boundary of the first 3-cell on the list we have
    \begin{equation}
        \begin{tikzpicture}[scale=2, baseline]
            \coordinate (A) at (0,0-0.5);
            \coordinate (B) at (1,0-0.5);
            \coordinate (C) at (1,1-0.5);
            \coordinate (D) at (0,1-0.5);

            \foreach \n in {A, B, C, D}
                \node at (\n) [circle, fill, inner sep=1pt] {};

            \draw[midarrow] (A) -- (B);
            \draw[midarrow] (B) -- (C);
            \draw[midarrow] (C) -- (D);
            \draw[midarrow] (A) -- (D);
            \draw[midarrow] (A) -- (C);

            \draw[double, double distance=1pt, -{Implies[]}] ({0.25+0.1}, {0.75+0.1-0.5}) -- ({0.25-0.1}, {0.75-0.1-0.5});
            
            \draw[double, double distance=1pt, -{Implies[]}] ({0.75+0.1}, {0.25-0.1-0.5}) -- ({0.75-0.1}, {0.25+0.1-0.5});
        \end{tikzpicture}
        \quad \rightsquigarrow \quad
        \begin{tikzpicture}[scale=2, baseline=1cm]
            \coordinate (A) at (0,0);
            \coordinate (B) at (1,0);
            \coordinate (C) at (1,1);
            \coordinate (D) at (0,1);

            \foreach \n in {A, B, C, D}
                \node at (\n) [circle, fill, inner sep=1pt] {};

            \draw[midarrow] (A) -- (B);
            \draw[midarrow] (B) -- (C);
            \draw[midarrow] (C) -- (D);
            \draw[midarrow] (A) -- (D);
            \draw[midarrow] (A) -- (C);

            \draw[double, double distance=1pt, -{Implies[]}] ({0.25+0.1}, {0.75+0.1}) -- ({0.25-0.1}, {0.75-0.1});
            
            \draw[double, double distance=1pt, -{Implies[]}] ({0.75+0.1}, {0.25-0.1}) -- ({0.75-0.1}, {0.25+0.1});

            \draw[very thick, color=green!50!black, midarrow] (0.25, 1.25) -- (0.25, 0.75);
            \draw[very thick, color=green!50!black, midarrow] (0.25, 0.75) -- (-0.25, 0.75);
            \draw[very thick, color=green!50!black, midarrow] (0.75, 0.25) -- (0.25, 0.75);
            \draw[very thick, color=green!50!black, midarrow] (0.75, -0.25) -- (0.75, 0.25);
            \draw[very thick, color=green!50!black, midarrow] (1.25, 0.25) -- (0.75, 0.25);
        \end{tikzpicture}
        \quad \rightsquigarrow \quad
        \begin{tikzpicture}[very thick, scale=1.2, color=green!50!black, baseline]
            \tktrival{1}{1}{0}{-1}{1}{1}
            \draw (0.5,-1) -- (0.5,0);
            \tktrival{1}{0}{0}{1}{1}{1}
        \end{tikzpicture}
    \end{equation}

    Translating the 2-dimensional source and target boundaries of a simplicial $3$-atom that are related to a 2-2 bistellar move into string diagrams immediately gives the relations
    \begin{equation}
        \begin{tikzpicture}[very thick, scale=0.5, color=green!50!black, baseline]
            \tktrival{1}{1}{0}{-1}{1}{1}
            \draw (0.5,-1) -- (0.5,0);
            \tktrival{1}{0}{0}{1}{1}{1}
        \end{tikzpicture} = 
        \begin{tikzpicture}[very thick, scale=0.5, color=green!50!black, baseline]
            \draw (-0.5,-1) -- (-0.5,0);
            \tktrival{1}{1}{1}{-1}{1}{1}
            \tktrival{1}{0}{0}{1}{1}{1}
        \end{tikzpicture}, \qquad
        \begin{tikzpicture}[very thick, scale=0.5, color=green!50!black, baseline]
            \tktrival{-1}{-1}{-1}{1}{1}{1}
            \draw (0.5,1) -- (0.5,0);
            \tktrival{-1}{0}{0}{-1}{1}{1}
        \end{tikzpicture} =
        \begin{tikzpicture}[very thick, scale=0.5, color=green!50!black, baseline]
            \tktrival{-1}{-1}{-1}{0}{1}{1}
            \draw (-1,1) -- (-1,0);
            \tktrival{-1}{0}{0}{0}{1}{1}
        \end{tikzpicture}, \qquad
        \begin{tikzpicture}[very thick, scale=0.5, color=green!50!black, baseline]
            \tktrival{-1}{0}{-1}{-1}{1}{1}
            \draw (0.5,-1) -- (0.5,0);
            \draw (-1.5, 0) -- (-1.5, 1);
            \tktrival{1}{0}{0}{1}{1}{1}
        \end{tikzpicture} = 
        \begin{tikzpicture}[very thick, scale=0.5, color=green!50!black, baseline]
            \tktrival{1}{0}{0}{0}{1}{1}
            \tktrival{-1}{0}{0}{0}{1}{1}
        \end{tikzpicture}, \qquad
        \begin{tikzpicture}[very thick, scale=0.5, color=green!50!black, baseline]
            \tktrival{1}{0}{0}{0}{1}{1}
            \tktrival{-1}{0}{0}{0}{1}{1}
        \end{tikzpicture} = 
        \begin{tikzpicture}[very thick, scale=0.5, color=green!50!black, baseline]
            \draw (-1,-1) -- (-1,0);
            \tktrival{-1}{0}{0.5}{-1}{1}{1}
            \tktrival{1}{-1}{-1}{0}{1}{1}
            \draw (1,0) -- (1,1);
        \end{tikzpicture}\, ,
    \end{equation}
    with simplicial $3$-atoms whose positive and negative boundaries are swapped contributing the same relation, as in the 1-dimensional case. From the $3$-atoms related to 3-1 bistellar moves we get exactly the relations
    \begin{equation}
        \begin{tikzpicture}[very thick, scale=0.6 * 2/3, color=green!50!black, baseline]
            \draw (-1,-1) -- (-1,0);
            \tktrival{-1}{0}{0.5}{-1}{1}{1}
            \tktrival{1}{-1}{-1}{0}{1}{1}
            \draw (1,0) -- (1,1);
            \tktrival{1}{-1}{-0.5}{1}{1.5}{1}
        \end{tikzpicture} =
        \begin{tikzpicture}[very thick, scale=0.6 * 2/3, color=green!50!black, baseline]
            \tktrival{1}{0}{0}{0}{1}{1}
            \draw (0,0) -- (0,2);
        \end{tikzpicture}, \qquad
        \begin{tikzpicture}[very thick, scale=0.6 * 2/3, color=green!50!black, baseline]
            \draw (1,-1) -- (1,0);
            \tktrival{-1}{0}{-0.5}{-1}{1}{1}
            \tktrival{1}{-1}{0}{0}{1}{1}
            \draw (-1,0) -- (-1,1);
            \tktrival{1}{-1}{-1}{1}{1.5}{1}
        \end{tikzpicture} =
        \begin{tikzpicture}[very thick, scale=0.6 * 2/3, color=green!50!black, baseline]
            \tktrival{1}{0}{0}{0}{1}{1}
            \draw (0,0) -- (0,2);
        \end{tikzpicture}, \qquad
        \begin{tikzpicture}[very thick, scale=0.6 * 2/3, color=green!50!black, baseline]
            \tktrival{-1}{-1}{-1.5}{0}{1.5}{1}
            \draw (0,0) -- (0,1);
            \tktrival{-1}{0}{-1.5}{0}{1}{1}
            \tktrival{1}{-1}{-1}{1}{1}{1}
            \draw (-2,1) -- (-2,2);
        \end{tikzpicture} =
        \begin{tikzpicture}[very thick, scale=0.6 * 2/3, color=green!50!black, baseline]
            \draw (0,-1) -- (0,1);
            \tktrival{-1}{0}{0}{1}{1}{1}
        \end{tikzpicture}, \qquad
        \begin{tikzpicture}[very thick, scale=0.6 * 2/3, color=green!50!black, baseline]
            \tktrival{-1}{-1}{-1.5}{0}{1.5}{1}
            \draw (-1.5,0) -- (-1.5,1);
            \tktrival{-1}{0}{0}{0}{1}{1}
            \tktrival{1}{-1}{-1.5}{1}{1}{1}
            \draw (0.5,1) -- (0.5,2);
        \end{tikzpicture} =
        \begin{tikzpicture}[very thick, scale=0.6 * 2/3, color=green!50!black, baseline]
            \draw (0,-1) -- (0,1);
            \tktrival{-1}{0}{0}{1}{1}{1}
        \end{tikzpicture},
    \end{equation}
    but these can be simplified by using the Frobenius relations to get the desired
    \begin{equation}
        \begin{tikzpicture}[very thick, scale=0.6 * 2/3, color=green!50!black, baseline]
            \tktrival{1}{0}{0}{0}{1}{1}
            \tktrival{-1}{0}{0}{0}{1}{1}
            \tktrival{1}{0}{0}{2}{1}{1}
        \end{tikzpicture} =
        \begin{tikzpicture}[very thick, scale=0.6 * 2/3, color=green!50!black, baseline]
            \tktrival{1}{0}{0}{0}{1}{1}
            \draw (0,0) -- (0,2);
        \end{tikzpicture}, \qquad 
        \begin{tikzpicture}[very thick, scale=0.6 * 2/3, color=green!50!black, baseline]
            \tktrival{-1}{0}{0}{-1}{1}{1}
            \tktrival{1}{0}{0}{1}{1}{1}
            \tktrival{-1}{0}{0}{1}{1}{1}
        \end{tikzpicture} =
        \begin{tikzpicture}[very thick, scale=0.6 * 2/3, color=green!50!black, baseline]
            \draw (0,-1) -- (0,1);
            \tktrival{-1}{0}{0}{1}{1}{1}
        \end{tikzpicture}.
    \end{equation}
\end{proof}

\begin{theorem}\label{thm:3_bistellar_datum}
    The walking 3-bistellar datum $\mathscr{P}_3$ has generators
    \begin{equation}
        a, \quad A: a \to a, \quad \mu: A \otimes A \to A, \quad \delta: A \to A \otimes A,
    \end{equation}
    as well as 3-generators:
    \begin{enumerate}
        \item an associator $\alpha$ and antiassociator $\hat{\alpha}$
        \begin{equation}
            \begin{tikzpicture}[very thick, scale=0.6, color=green!50!black, baseline]
                \tktrival{1}{1}{0}{-1}{1}{1}
                \draw (0.5,-1) -- (0.5,0);
                \tktrival{1}{0}{0}{1}{1}{1}
            \end{tikzpicture} \quad \stackrel{\alpha}{\Rrightarrow} \quad
            \begin{tikzpicture}[very thick, scale=0.6, color=green!50!black, baseline]
                \draw (-0.5,-1) -- (-0.5,0);
                \tktrival{1}{1}{1}{-1}{1}{1}
                \tktrival{1}{0}{0}{1}{1}{1}
            \end{tikzpicture}\, , \qquad \qquad
            \begin{tikzpicture}[very thick, scale=0.6, color=green!50!black, baseline]
                \tktrival{1}{1}{0}{-1}{1}{1}
                \draw (0.5,-1) -- (0.5,0);
                \tktrival{1}{0}{0}{1}{1}{1}
            \end{tikzpicture} \quad \stackrel{\hat{\alpha}}{\Lleftarrow} \quad
            \begin{tikzpicture}[very thick, scale=0.6, color=green!50!black, baseline]
                \draw (-0.5,-1) -- (-0.5,0);
                \tktrival{1}{1}{1}{-1}{1}{1}
                \tktrival{1}{0}{0}{1}{1}{1}
            \end{tikzpicture}\, ,
        \end{equation}
        \item a coassociator $\gamma$ and anticoassociator $\hat{\gamma}$
        \begin{equation}
            \begin{tikzpicture}[very thick, scale=0.6, color=green!50!black, baseline]
                \tktrival{-1}{-1}{-1}{1}{1}{1}
                \draw (0.5,1) -- (0.5,0);
                \tktrival{-1}{0}{0}{-1}{1}{1}
            \end{tikzpicture} \quad \stackrel{\gamma}{\Rrightarrow} \quad
            \begin{tikzpicture}[very thick, scale=0.6, color=green!50!black, baseline]
                \tktrival{-1}{-1}{-1}{0}{1}{1}
                \draw (-1,1) -- (-1,0);
                \tktrival{-1}{0}{0}{0}{1}{1}
            \end{tikzpicture}\, , \qquad \qquad
            \begin{tikzpicture}[very thick, scale=0.6, color=green!50!black, baseline]
                \tktrival{-1}{-1}{-1}{1}{1}{1}
                \draw (0.5,1) -- (0.5,0);
                \tktrival{-1}{0}{0}{-1}{1}{1}
            \end{tikzpicture} \quad \stackrel{\hat{\gamma}}{\Lleftarrow} \quad
            \begin{tikzpicture}[very thick, scale=0.6, color=green!50!black, baseline]
                \tktrival{-1}{-1}{-1}{0}{1}{1}
                \draw (-1,1) -- (-1,0);
                \tktrival{-1}{0}{0}{0}{1}{1}
            \end{tikzpicture}\, ,
        \end{equation}
        \item two Frobeniusators $\phi_L, \phi_R$ and two antiFrobeniusators $\hat{\phi}_L, \hat{\phi}_R$
        \begin{equation}
            \begin{tikzpicture}[very thick, scale=0.6, color=green!50!black, baseline]
                \tktrival{-1}{0}{-1}{-1}{1}{1}
                \draw (0.5,-1) -- (0.5,0);
                \draw (-1.5, 0) -- (-1.5, 1);
                \tktrival{1}{0}{0}{1}{1}{1}
            \end{tikzpicture} \quad \centstack{\stackrel{\phi_L}{\Rrightarrow}}{\underset{\hat{\phi}_L}{\Lleftarrow}} \quad
            \begin{tikzpicture}[very thick, scale=0.6, color=green!50!black, baseline]
                \tktrival{1}{0}{0}{0}{1}{1}
                \tktrival{-1}{0}{0}{0}{1}{1}
            \end{tikzpicture}\, , \qquad \qquad
            \begin{tikzpicture}[very thick, scale=0.6, color=green!50!black, baseline]
                \tktrival{1}{0}{0}{0}{1}{1}
                \tktrival{-1}{0}{0}{0}{1}{1}
            \end{tikzpicture}
             \quad \centstack{\stackrel{\phi_R}{\Lleftarrow}}{\underset{\hat{\phi}_R}{\Rrightarrow}} \quad
            \begin{tikzpicture}[very thick, scale=0.6, color=green!50!black, baseline]
                \draw (-1,-1) -- (-1,0);
                \tktrival{-1}{0}{0.5}{-1}{1}{1}
                \tktrival{1}{-1}{-1}{0}{1}{1}
                \draw (1,0) -- (1,1);
            \end{tikzpicture}\, ,
        \end{equation}
        \item four bubbleators $\beta_L^U, \beta_R^U, \beta_L^D, \beta_R^D$ and four antibubbleators $\hat{\beta}_L^U, \hat{\beta}_R^U, \hat{\beta}_L^D, \hat{\beta}_R^D$
        \begin{equation}
            \begin{tikzpicture}[very thick, scale=0.6 * 2/3, color=green!50!black, baseline]
            \draw (-1,-1) -- (-1,0);
            \tktrival{-1}{0}{0.5}{-1}{1}{1}
            \tktrival{1}{-1}{-1}{0}{1}{1}
            \draw (1,0) -- (1,1);
            \tktrival{1}{-1}{-0.5}{1}{1.5}{1}
        \end{tikzpicture} \quad \centstack{\stackrel{\hat{\beta}_R^U}{\Rrightarrow}}{\underset{\beta_R^U}{\Lleftarrow}} \quad
        \begin{tikzpicture}[very thick, scale=0.6 * 2/3, color=green!50!black, baseline]
            \tktrival{1}{0}{0}{0}{1}{1}
            \draw (0,0) -- (0,2);
        \end{tikzpicture}\, , \quad \ \ 
        \begin{tikzpicture}[very thick, scale=0.6 * 2/3, color=green!50!black, baseline]
            \draw (1,-1) -- (1,0);
            \tktrival{-1}{0}{-0.5}{-1}{1}{1}
            \tktrival{1}{-1}{0}{0}{1}{1}
            \draw (-1,0) -- (-1,1);
            \tktrival{1}{-1}{-1}{1}{1.5}{1}
        \end{tikzpicture} \quad \centstack{\stackrel{\hat{\beta}_L^U}{\Rrightarrow}}{\underset{\beta_L^U}{\Lleftarrow}} \quad
        \begin{tikzpicture}[very thick, scale=0.6 * 2/3, color=green!50!black, baseline]
            \tktrival{1}{0}{0}{0}{1}{1}
            \draw (0,0) -- (0,2);
        \end{tikzpicture}\, , \quad \ \ 
        \begin{tikzpicture}[very thick, scale=0.6 * 2/3, color=green!50!black, baseline]
            \tktrival{-1}{-1}{-1.5}{0}{1.5}{1}
            \draw (0,0) -- (0,1);
            \tktrival{-1}{0}{-1.5}{0}{1}{1}
            \tktrival{1}{-1}{-1}{1}{1}{1}
            \draw (-2,1) -- (-2,2);
        \end{tikzpicture} \quad \centstack{\stackrel{\hat{\beta}_R^D}{\Rrightarrow}}{\underset{\beta_R^D}{\Lleftarrow}} \quad
        \begin{tikzpicture}[very thick, scale=0.6 * 2/3, color=green!50!black, baseline]
            \draw (0,-1) -- (0,1);
            \tktrival{-1}{0}{0}{1}{1}{1}
        \end{tikzpicture}\, , \quad \ \ 
        \begin{tikzpicture}[very thick, scale=0.6 * 2/3, color=green!50!black, baseline]
            \tktrival{-1}{-1}{-1.5}{0}{1.5}{1}
            \draw (-1.5,0) -- (-1.5,1);
            \tktrival{-1}{0}{0}{0}{1}{1}
            \tktrival{1}{-1}{-1.5}{1}{1}{1}
            \draw (0.5,1) -- (0.5,2);
        \end{tikzpicture} \quad \centstack{\stackrel{\hat{\beta}_L^D}{\Rrightarrow}}{\underset{\beta_L^D}{\Lleftarrow}} \quad
        \begin{tikzpicture}[very thick, scale=0.6 * 2/3, color=green!50!black, baseline]
            \draw (0,-1) -- (0,1);
            \tktrival{-1}{0}{0}{1}{1}{1}
        \end{tikzpicture}\, .
        \end{equation}
    \end{enumerate}
    These 3-generators, all of which are noninvertible, satisfy a set of 620 relations. These are organized and presented in \Cref{thm:satm_4}.
\end{theorem}

\begin{proof}
    All the generators can be recognized from \Cref{thm:satm_in_low_d}. The nontrivial part is the relations, which are given by 4-dimensional simplicial atoms. There are 1240 of these which can be found by a computer. A description of one possible search algorithm is described in \Cref{app:pseudocode}.
    Under the flip of the direction of the 4-cell these group into pairs such that each induces the same relation giving the most naive number of relations as 620. If we take into account the flip symmetry at all levels then relations like for example the pentagon axioms for all of the associator, antiassociator, coassociator and anticoassociator will all be related by flips and the number of independent classes like that is 80, which are the ones presented in \Cref{thm:satm_4}.
\end{proof}

\begin{remark}
    The set of 620 relations is not guaranteed to be minimal. In fact, it is not expected to be, because it is a presentation that preserves a lot of the symmetry of the underlying diagrams, as we saw in the two-dimensional case. To get an idea for the organizational principle, if we had organized the relations in the two-dimensional case so that they are generated under flips by some subset then that subset would have 3 elements: (co)associativity, Frobenius relations and bubble relations.
\end{remark}

\begin{remark}
    Some of the relations that appear on the list are familiar in the literature, like the pentagon identity for the associator --- which also comes along with pentagon identities for the antiassociator, coassociator and anticoassociator due to the flip symmetry. Others are generalizations of what the Frobenius property should be in one dimension higher, and yet others describe how exactly three-dimensional bubbles are allowed to appear and disappear.
\end{remark}

\section{Examples and Special Cases}\label{sec:special_cases}

\setcounter{subsection}{-1}
\subsection{Dimension 2}

Even though we specify the walking $n$-bistellar datum as a strict $n$-category, we can still map out of it, and into a weak 2-category. Bistellar data will then consist of not just the literal images of the generators under the mapping, but also coherence data. We write this out in the 2-dimensional case.

\begin{example}[Weak 2-bistellar data]
    Let $F: \mathscr{P}_2 \to \mathscr{C}$ be a weak 2-functor (a pseudofunctor) into a weak 2-category $\mathscr{C}$, or in other words, a weak $2$-bistellar datum. The datum $F$ consists of
    \begin{equation}
        \big( F(a), \ F(1_a), \ F(A^{\otimes n}), \ F(\mu), \ F(\delta), \ \kappa_{n,m}, \ \iota \big),
    \end{equation}
    where $n, m \in \mathbb{N}$, and where $F(a)$ is an object, $F(1_a)$ and $F(A^{\otimes n})$, for some $n \in \mathbb{N}$, are 1-endomorphism of $F(a)$, $F(\mu)$ is a 2-morphism $F(A^{\otimes 2}) \to F(A)$ and $\delta$, a 2-morphism that goes in the opposite direction of $\mu$. Finally, the coherences are
    \begin{equation}
        \kappa_{n,m}: F(A^{\otimes n}) \otimes F(A^{\otimes m}) \xrightarrow{\ \ \cong \ \ } F(A^{\otimes n + m}), \qquad \qquad \iota: 1_{F(a)} \xrightarrow{\ \ \cong \ \ } F(1_a),
    \end{equation}
    and they satisfy:
    \begin{enumerate}
        \item naturality for $\kappa_{m,n}$:
        \begin{equation}
            \begin{tikzcd}[column sep=large]
                F(A^{\otimes n}) \otimes F(A^{\otimes m}) \arrow[r, "F(f) \otimes F(g)"] \arrow[d, "\kappa_{n,m}"] & F(A^{\otimes p}) \otimes F(A^{\otimes q}) \arrow[d, "\kappa_{p,q}"] \\
                F(A^{\otimes n + m}) \arrow[r, "F(f \otimes g)"] & F(A^{\otimes p + q})
            \end{tikzcd},
        \end{equation}
        where $f: A^{\otimes n} \to A^{\otimes p}, \ g: A^{\otimes m} \to A^{\otimes q}$ are given by string diagrams,
        \item left unitality
        \item right unitality
        \item the hexagon axiom
        \begin{equation}
            \begin{tikzcd}
                F(A^{\otimes n}) \otimes \left( F(A^{\otimes m}) \otimes F(A^{\otimes r}) \right) \arrow[r, "{\alpha_{n, m, r}}"] \arrow[d, "{1_{F(A^{\otimes n})} \otimes \kappa_{m, r}}"'] & \left( F(A^{\otimes n}) \otimes F(A^{\otimes m}) \right) \otimes F(A^{\otimes r}) \arrow[d, "{\kappa_{n,m} \otimes 1_{F(A^{\otimes r})}}"] \\
                F(A^{\otimes n}) \otimes F(A^{\otimes m + r}) \arrow[d, "{\kappa_{n, m + r}}"'] & F(A^{\otimes n + m}) \otimes F(A^{\otimes r}) \arrow[d, "{\kappa_{n + m, r}}"] \\
                F(A^{\otimes n + m + r}) \arrow[r, Rightarrow, no head] & F(A^{\otimes n + m + r}),
            \end{tikzcd}
        \end{equation}
        where we used the abbreviation $\alpha_{n,m,r} := \alpha_{F(A^{\otimes n}), F(A^{\otimes m}), F(A^{\otimes r})}$.
    \end{enumerate}
    We do not spell out the unitality axioms because 2-bistellar data are nonunital and noncounital, so the role of the identity is not that important. If one thinks of $\tilde{A} := F(A)$ as being the underlying 1-morphism of an algebra with multiplication and comultiplication given by
    \begin{equation}
        \begin{tikzcd}
            \tilde{\mu}: F(A) \otimes F(A) \arrow[r, shift left = 1ex, "\kappa_{1,1}"] & F(A \otimes A) \arrow[l, shift left = 1ex, "\kappa_{1,1}^{-1}"] \arrow[r, shift left = 1ex, "F(\mu)"] & F(A) \arrow[l, shift left = 1ex, "F(\delta)"]: \tilde{\delta}
        \end{tikzcd}
    \end{equation}
    then this is an algebra that is associative in the appropriately weak sense, namely
    \begin{equation}
        \begin{tikzcd}[row sep = tiny]
            (\tilde{A} \otimes \tilde{A}) \otimes \tilde{A} \arrow[r, "\tilde{\mu} \otimes 1"] \arrow[dd, "\alpha_{\tilde{A},\tilde{A},\tilde{A}}"'] & \tilde{A} \otimes \tilde{A} \arrow[rd, "\tilde{\mu}"] \\
            & & \tilde{A} \\
            \tilde{A} \otimes (\tilde{A} \otimes \tilde{A}) \arrow[r, "1 \otimes \tilde{\mu}"] & \tilde{A} \otimes \tilde{A} \arrow[ru, "\tilde{\mu}"']
        \end{tikzcd},
    \end{equation}
    where $\alpha$ is the associator of the weak 2-category $C$.\footnote{Not to be confused with the associator of the walking 3-bistellar datum, which is not invertible.} Similarly, it is coassociative, Frobenius, and it also satisfies the bubble omission relations in a way that is weak in the same sense. Since we prefer to write out such relations in the string diagram (or two-dimensional graphical calculus) representation, if we were to write the relations out we would not see any difference, since the two-dimensional graphical calculus relies on strictification theorems to avoid writing out this coherence data.
\end{example}

\begin{example}\label{ex:P2withisosin2}
    Suppose that the orbifold datum $\mathscr{P}_2 \to \mathscr{C}$, in some weak $2$-category $\mathscr{C}$, factors through the maximal $(2,1)$-category $t_{(2,1)}\mathscr{C}$, so that the 2-cells of the orbifold datum are invertible. In this case, by acting with the inverse of the multiplication on the bubble move we get
    \begin{equation}
        \begin{tikzpicture}[very thick, scale=0.8 * 2/3, color=green!50!black, baseline]
            \draw (0, -1.5) -- (0, 2);
        \end{tikzpicture} \ \ \ = \ \ \ 
        \begin{tikzpicture}[very thick, scale=0.8 * 2/3, color=green!50!black, baseline]
            \tktrival{-1}{0}{0}{-0.5}{1}{1}
            \draw (0-0.2, 0-0.2) -- (0+0.2, 0+0.2);
            \draw (0+0.2, 0-0.2) -- (0-0.2, 0+0.2);
            \tktrival{1}{0}{0}{1.5}{1}{1}
            \draw (0, -1.5) -- (0, -0.5);
            \draw (0, 1.5) -- (0, 2);
        \end{tikzpicture} \ \ \ = \ \ \ 
        \begin{tikzpicture}[very thick, scale=0.8 * 2/3, color=green!50!black, baseline]
            \tktrival{1}{0}{0}{0.5}{1}{1}
            \tktrival{-1}{0}{0}{0}{1}{1}
            \tktrival{1}{0}{0}{2}{1}{1}
            \tktrival{-1}{0}{0}{-1.5}{1}{1}
            \draw (0-0.2, -1-0.2) -- (0+0.2, -1+0.2);
            \draw (0+0.2, -1-0.2) -- (0-0.2, -1+0.2);
        \end{tikzpicture} \ \ \ = \ \ \ 
        \begin{tikzpicture}[very thick, scale=0.8 * 2/3, color=green!50!black, baseline]
            \tktrival{-1}{0}{0}{-0.5}{1}{1}
            \tktrival{1}{0}{0}{1.5}{1}{1}
            \draw (0, -1.5) -- (0, -0.5);
            \draw (0, 1.5) -- (0, 2);
        \end{tikzpicture}\, ,
    \end{equation}
    where we used a crossed-out vertex to denote the inverse. This means that $\mu^{-1} = \delta$ by the uniqueness of inverses. It is a classical result that given the data of a special Frobenius algebra, there is a unique connected map $A^{\otimes n} \to A^{\otimes m}$, for each $n, m \in \mathbb{N} \setminus \{0\}$. Using only specialness, the disconnected maps would be given by acting with $\otimes$ on the connected ones, however we also have the extra condition
    \begin{equation}
        \begin{tikzpicture}[very thick, scale=0.6, color=green!50!black, baseline]
            \tktrival{1}{0}{0}{0}{1}{1}
            \tktrival{-1}{0}{0}{0}{1}{1}
        \end{tikzpicture} \ \ \ = \ \ \ 
        \begin{tikzpicture}[very thick, scale=0.6, color=green!50!black, baseline]
            \draw (-0.5, -1) -- (-0.5, 1);
            \draw (0.5, -1) -- (0.5, 1);
        \end{tikzpicture}\, ,
    \end{equation}
    whose iterated use makes any disconnected diagram equal to the unique connected one for the given input and output. Thus, 2-bistellar data in this situation are given by an object $a$, a 1-morphism $A$ and isomorphisms $c_{n, m}: A^{\otimes n} \to A^{\otimes m}$, one for each pair $n, m \in \mathbb{N} \setminus \{0\}$ such that
    \begin{equation}
        c_{m, k} \circ c_{n, m} = c_{n, k}, \qquad \qquad c_{n, n} = 1_{A^{\otimes n}}.
    \end{equation}
    This is one particular categorification of what it would mean to generalize idempotents --- the 1-dimensional bistellar data --- to two dimensions, namely providing isomorphisms for all possible $\otimes$-compositions of $A$ with itself. 

    It is a very limited categorification though, because if $\mathscr{C} = \mathsf{B} \vect_{\mathbb{k}}$, then the only possible choices of orbifold data would be $A = \mathbb{k}$ and $\mu_z (x \otimes y) = z x y$, for some invertible $z \in \mathbb{k}^{\times}$. In particular, matrix algebras, the classic examples of Frobenius algebras, cannot be included. 
\end{example}

\begin{example}
    In \cite{orbifold2d} it is shown that unital, counital, symmetric, special Frobenius algebras satisfy the $2$-bistellar moves. Moreover, $2$-condensation monads are exactly nonunital, noncounital, special Frobenius algebras, as seen in \cite{condensation}. From the perspective of molecules, though, the specialness condition is too algebraic, because it equates a non-identity 2-morphism with the identity of a 1-morphism. There is no 3-molecule which has positive and negative 2-boundaries that are of different dimensions, so this condition cannot be implemented in the language of molecules. For the same reason unitality and counitality are also not geometrically natural conditions.
    
    The closest thing we can get to implementing such conditions is to move away from molecules and consider general oriented graded posets like
    \begin{equation}
        \left(U \mathop{\amalg}_{\partial_0 U = \partial_0 V} V \right)^{\vartriangleright} \qquad \mathrm{with} \qquad U = \begin{tikzpicture}[scale=0.5, baseline]
            \node at (-2,0.2) [circle, fill, inner sep=1pt] {};
            \node at (0,0.2) [circle, fill, inner sep=1pt] {};
            \node at (2,0.2) [circle, fill, inner sep=1pt] {};

            \draw[midarrow] (-2,0.2) -- (0,0.2);
            \draw[midarrow] (0,0.2) -- (2,0.2);
            \centerarc[midarrow](0, 0.2)(-180:0:2);
            \centerarc[midarrow](0, 0.2)(180:0:2);

            \draw[double, double distance=1pt, -{Implies[]}] ({0+0}, {-0.7-0.3}) -- ({0+0}, {-0.7+0.3});
                
            \draw[double, double distance=1pt, -{Implies[]}] ({0+0}, {1.1-0.3}) -- ({0+0}, {1.1+0.3});
        \end{tikzpicture}, \quad 
        V = \begin{tikzpicture}[scale=0.5, baseline]
            \node at (0,0.2) [circle, fill, inner sep=1pt] {};
            \node at (2,0.2) [circle, fill, inner sep=1pt] {};

            \draw[midarrow] (0,0.2) -- (2,0.2);
        \end{tikzpicture}.
    \end{equation}
    If we add this particular oriented graded poset to our set of relations then the resulting strict $2$-category will be the walking nonunital, noncounital, special Frobenius algebra. But, as is exemplified by the fact that we cannot even provide a picture of it, this would be a geometrically unmotivated move.
\end{example}

On the other hand there's clearly a very close relationship between $2$-bistellar algebras and special Frobenius algebras, so it's natural to ask how closely related they are. In the presence of mild structure we have the following two lemmata:

\begin{lemma}\label{lem:2_bistellar_decomposition}
    Let $\mathscr{C}$ be a category weakly enriched in idempotent complete pre-additive categories, which we consider as a weak 2-category with extra structure. Then every $2$-bistellar datum $A: \mathscr{P}_2 \to \mathscr{C}$ decomposes into $A \cong \tilde{A} \oplus \tilde{B}$, where $\tilde{A}$ is a nonunital, noncounital, special Frobenius algebra, and $\tilde{B}$ is an algebra with trivial multiplication and comultiplication.
\end{lemma}

\begin{proof}
    Denote the $2$-bistellar algebra's data by $(a, A, \mu, \delta)$. The relations imply that the bubble morphism $\beta := \mu \delta$ is idempotent. Using idempotent completeness we can split $\beta$ so that we get an object $\tilde{A}:= \mathrm{im}(\beta)$, and $1 - \beta$, so that we get an object $\tilde{B} := \mathrm{im}(1 - \beta)$. 
    
    It's a classical fact that $A \cong \tilde{A} \oplus \tilde{B}$. Indeed, if we denote the splitting maps for $\beta$ with $p: \tilde{A} \leftrightharpoons P: j$ and the ones for $1 - \beta$ with $p'$ and $j'$ then $p$ and $p'$ exhibit $P$ as the coproduct of $\tilde{A}$ and $\tilde{B}$, because given any $X$ with maps $f: \tilde{A} \to X$ and $g: \tilde{B} \to X$ there is a unique map, $h := f j + g j'$, that satisfies $h p = f$ and $h p' = g$. The key calculation is
    \begin{equation}
        j' p = j' p' j' p j p = j' (1 - \beta) \beta p = j' 0 p = 0,
    \end{equation}
    and similarly for $j p'$. Under this decomposition $\mu$ and $\delta$ can be viewed as maps
    \begin{equation}
        (\tilde{A} \otimes \tilde{A}) \oplus (\tilde{A} \otimes \tilde{B}) \oplus (\tilde{B} \otimes \tilde{A}) \oplus (\tilde{B} \otimes \tilde{B}) \cong (\tilde{A} \oplus \tilde{B}) \otimes (\tilde{A} \oplus \tilde{B}) \xrightleftharpoons[\delta]{\ \ \mu \ \ } \tilde{A} \oplus \tilde{B}.
    \end{equation}
    Using the universal properties of the direct sum, the data of the multiplication and comultiplication is fully determined by the set of maps
    \begin{align}
        \tilde{A} \otimes \tilde{A} \xrightleftharpoons{\ \ \ \ \ \ } \tilde{A}, \qquad \qquad \tilde{A} \otimes \tilde{A} \xrightleftharpoons{\ \ \ \ \ \ } \tilde{B}, \notag\\
        \tilde{A} \otimes \tilde{B} \xrightleftharpoons{\ \ \ \ \ \ } \tilde{A}, \qquad \qquad \tilde{A} \otimes \tilde{B} \xrightleftharpoons{\ \ \ \ \ \ } \tilde{B}, \notag\\
        \tilde{B} \otimes \tilde{A} \xrightleftharpoons{\ \ \ \ \ \ } \tilde{A}, \qquad \qquad \tilde{B} \otimes \tilde{A} \xrightleftharpoons{\ \ \ \ \ \ } \tilde{B}, \notag\\
        \tilde{B} \otimes \tilde{B} \xrightleftharpoons{\ \ \ \ \ \ } \tilde{A}, \qquad \qquad \tilde{B} \otimes \tilde{B} \xrightleftharpoons{\ \ \ \ \ \ } \tilde{B}.
    \end{align}
    However, because of the definition of $\tilde{B}$, any of the above maps that involve $\tilde{B}$ in either the domain or codomain have to be the zero map. In the case of $\tilde{A} \otimes \tilde{B} \to \tilde{A}$, for example, the bubble omission property states that
    \begin{equation}
        \begin{tikzcd}
            \tilde{A} \otimes \tilde{B} \arrow[rd] \arrow[r, "1_{\tilde{A}} \otimes j' \beta p'"] & \tilde{A} \otimes \tilde{B} \arrow[d] \\
            & \tilde{A}  
        \end{tikzcd}
    \end{equation}
    commutes. But $j' \beta p' = 0$, so the multiplication map component also has to be 0. The only nontrivial components in the decomposition are
    \begin{equation}
        \tilde{A} \otimes \tilde{A} \xrightleftharpoons[\delta|_{\tilde{A}} :=  (j \otimes j) \delta p]{\ \ \mu|_{\tilde{A}} := j \mu (p \otimes p) \ \ } \tilde{A},
    \end{equation}
    and it's easily seen that they satisfy all the properties of a Frobenius algebra. Furthermore, by the definition of $\tilde{A}$, clearly $\mu|_{\tilde{A}} \delta|_{\tilde{A}} = j \beta p = 1$, satisfying specialness as well.
\end{proof}

\begin{lemma}\label{lem:2_bistellar_with_units}
    Let $\mathscr{C}$ be a weak 2-category, and $A: \mathscr{P}_2 \to \mathscr{C}$ be a 2-bistellar datum. If either the multiplication of $A$ is unital or the comultiplication of $A$ is counital then $A$ is a special Frobenius algebra.
\end{lemma}

\begin{proof}
    In the presence of a unit
    \begin{equation}
        \begin{tikzpicture}[very thick, scale=0.7 * 0.9, color=green!50!black, baseline]
            \tktrival{-1}{0}{0}{-0.75}{1}{1}
            \tktrival{1}{0}{0}{1.25}{1}{1}    
        \end{tikzpicture}
        \quad = \quad
        \begin{tikzpicture}[very thick, scale=0.7 * 2/3, color=green!50!black, baseline]
            \draw (-1,-1) -- (-1,0);
            \tktrival{-1}{0}{0.5}{-1}{1}{1}
            \tktrival{1}{-1}{-1}{0}{1}{1}
            \draw (1,0) -- (1,1);
            \tktrival{1}{-1}{-0.5}{1}{1.5}{1}
            \draw (-1, -1) node[Odot] {};
        \end{tikzpicture}
        \quad = \quad
        \begin{tikzpicture}[very thick, scale=0.7 * 2/3, color=green!50!black, baseline]
            \tktrival{1}{0}{0}{1}{1.5}{1.5}
            \draw (-0.75, -1) -- (-0.75, -0.5);
            \draw (0.75, -1) -- (0.75, -0.5);
            \draw (0, 1) -- (0, 1.5);
            \draw (-0.75, -1) node[Odot] {};
        \end{tikzpicture}
        \quad = \quad
        \begin{tikzpicture}[very thick, scale=0.7 * 2/3, color=green!50!black, baseline]
            \draw (0, -1) -- (0, 1.5);
        \end{tikzpicture}.
    \end{equation}
    For a counit we use the pictures flipped across the horizontal.
\end{proof}

\begin{remark}
    If we combine the assumptions of the two lemmata we see that either unitality of the multiplication or counitality of the comultiplication means that the trivial part $\tilde{B} = 0$, because a trivial (co)multiplication can never be (co)unital.
\end{remark}

\begin{example}
    The category $\vect_{\mathbb{k}}$ of finite dimensional vector spaces, linear maps, and with the tensor product of vector spaces, is a monoidal category, which is additive and idempotent complete. Thus, $\mathsf{B} \vect_{\mathbb{k}}$ satisfies the assumptions of \Cref{lem:2_bistellar_decomposition} and every 2-bistellar datum in it is a trivial extension of a special Frobenius algebra, while every (co)unital such datum is exactly a special Frobenius algebra.

    The most famous family of examples of special Frobenius algebras in this context are the matrix algebras $\mathsf{Mat}_{n \times n}(\mathbb{k})$ with multiplication given by matrix multiplication, unit given by the identity matrix and counit given by the trace. This data uniquely determines a comultiplication, by observing that the pairing given by multiplying and then taking the trace is non-degenerate, and using its induced copairing to define
    \begin{equation}
        \begin{tikzpicture}[very thick, scale=0.7 * 0.9, color=green!50!black, baseline]
            \tktrival{-1}{0}{0}{-0.2}{1}{1}
        \end{tikzpicture} \quad := \quad
        \begin{tikzpicture}[very thick, scale=0.7 * 0.9, color=green!50!black, baseline]
            \tktrival{1}{0}{0}{0.8}{1}{1}
            \tkarc{-1}{1}{-0.5}{-0.2}{0.5}{1}
            \tkarc{-1}{1}{0.5}{-0.7}{2}{1}
            \tkarc{1}{-1}{0}{0.8}{1}{1}
            \draw (-1.5, -0.7) -- (-1.5, 1.5);
            \draw (-1, -0.2) -- (-1, 1.5);
            \draw (0.5, -0.2) -- (0.5, -0.7);
            \draw (1, 0.8) -- (1, -1.5);
        \end{tikzpicture} \ .
    \end{equation}
\end{example}

\begin{remark}
    For the purposes of the orbifold construction there is a motivated expectation that one should do the construction iteratively. In dimension 2 that means that one should first do the 1-dimensional orbifold construction to the hom-categories of the target and only then do the 2-dimensional orbifold construction on the result. Since the walking 1-bistellar datum is the walking idempotent, in such use cases one of the assumptions of \Cref{lem:2_bistellar_decomposition} is already satisfied.
\end{remark}

\subsection{Dimension 3}\label{ssec:special_orbifold_data}

Just like in \Cref{ex:P2withisosin2}, we can consider the case of three-dimensional bistellar data factoring through the maximal $(3,2)$-subcategory $t_{(3,2)} \mathscr{C}$ of a given 3-category $\mathscr{C}$. This is a simplification in which instead of the having noninvertible (co)associators, Frobeniusators and bubbleators, we instead have invertible version of those, so that we are closer to the more classically studied situation of having associativity and Frobenius axioms that hold only up to specified isomorphisms. A similar kind of setup was previously considered in \cite{orbifoldnd} and then subsequently in \cite{orbifold3d}, and we will work up to showing that:
\begin{proposition}\label{prop:3_bistellar_from_special_orbifold}
    A special orbifold datum in a Gray category with duals $\mathscr{C}$, as defined in \cite[Def.\,4.1]{orbifold3d}, gives a representation of the walking 3-bistellar datum on the maximal $(3,2)$-subcategory of $\mathscr{C}$,
    \begin{equation}
        \left\{ \begin{gathered}
        \emph{\text{special orbifold data}} \\ \emph{\text{in }} \mathscr{C}
    \end{gathered} \right\} \xrightarrow{\qquad} \mathsf{Fun}(\mathscr{P}_3, t_{(3,2)} \mathscr{C}).
    \end{equation}
\end{proposition}

First we introduce the relevant definitions, which make substantial use of the graphical calculus for Gray categories developed in \cite{graphicalcalculus}. The three-dimensional graphical calculus pictures defining special orbifold data, which are reproduced below, will all be adapted from \cite{orbifold3d}.

The three-dimensional graphical calculus is a calculation method developed for Gray categories with or without duals. It represents the various morphisms of a Gray category by labeling a stratification of the standard cube. Objects are represented by 3-dimentional strata, 1-morphisms by 2-dimensional strata, 2-morphisms by 1-dimensional strata and 3-morphisms by 0-dimensional strata. The compositions are determined by stacking in the 3 different directions. One possible convention, which we will use below, is 
\begin{center}
    \begin{tabular}{c|c|c}
        Level & Composition & Directionality \\
        \hline
        1 & $\boxtimes$ & front to back \\
        2 & $\otimes$ & right to left \\
        3 & $\circ$ & bottom to top \\
    \end{tabular}
\end{center}
Reading a whole three-dimensional diagram then results in a specific 3-morphism in the Gray category, that may be an identity, in which case we're representing a morphism of a lower level, or a nonidentity 3-morphism. The sources and targets of the 3-morphism are presented as the boundary faces of the cube in which the diagram fits, employing the directionality convention again to discern sources from targets. When duals are present then bending the appropriate lines or planes provides the duality. In equations, we denote 1-duals by $-^\#$ and 2-duals by $-^\vee$. For a detailed review of the graphical calculus in the same conventions as above, see \cite[\S 3]{graycat}.

The following definition also provides a good example of translating well-known relations between 3-morphisms into the language of the three-dimensional graphical calculus.

\begin{definition}
    An \emph{$\mathbb{E}_1$-algebra} in a Gray category is a collection
    \begin{equation}
        \left(a, \ \ A, \ \ \mu, \ \ u, \ \ \alpha, \ \ u^{\mathrm{l}}, \ \ u^{\mathrm{r}} \right),
    \end{equation}
    where $a$ is an object, $A: a \to a$ is a 1-endomorphism, $\mu: A \boxtimes A \to A$ is a multiplication 2-morphism, $u: 1_a \to A$ is a unit 2-morphism, and finally, $\alpha$, $u^{\mathrm{l}}$, and $u^{\mathrm{r}}$ are an associator, and left and right unitors, respectively, which are invertible 3-cells, satisfying the usual relations
    \begin{equation}\label{eq:3d_graphical_pentagon}
        \begin{tikzpicture}[thick, xscale=1.77, yscale=1.5,color=blue!50!black, baseline=0.5cm, >=stealth,
		style={x={(-0.55cm,-0.4cm)},y={(1cm,-0.2cm)},z={(0cm,0.9cm)}}]
		\pgfmathsetmacro{\yy}{0.2}
		\coordinate (P) at (0.5, \yy, 0);
		\coordinate (R) at (0.625, 0.5 + \yy/2, 0);
		\coordinate (L) at (0.5, 0, 0);
		\coordinate (R1) at (0.25, 1, 0);
		\coordinate (R2) at (0.5, 1, 0);
		\coordinate (R3) at (0.75, 1, 0);
		\coordinate (N1) at ($(R1) + (R2) - (R)$);
		\coordinate (N2) at ($(R) + 2*(R2) - 2*(R)$);
		\coordinate (N3) at ($(R3) + (R2) - (R)$);
		\coordinate (N4) at ($2*(R3) - (R)$);
		\coordinate (Lo) at (0.5, 0, 2);
		\coordinate (Po) at (0.5, \yy, 2);
		\coordinate (O1) at ($(N1) + (0,0,2)$);
		\coordinate (O2) at ($(N2) + (0,0,2)$);
		\coordinate (O3) at ($(N3) + (0,0,2)$);
		\coordinate (O4) at ($(N4) + (0,0,2)$);
		\coordinate (A2) at ($(O2) - (R3) + (R)$);
		\coordinate (A3) at ($(O3) - 2*(R3) + 2*(R)$);
		\coordinate (A4) at ($(O4) - 3*(R3) + 3*(R)$);
		\coordinate (Pt) at (0.5, \yy, 1);
		\coordinate (Rt) at ($(A2) + (0,0,-1)$);
		\coordinate (Lt) at (0.5, 0, 1);
		\coordinate (M1) at ($(N1) + (0,0,1)$);
		\coordinate (M2) at ($(N2) + (0,0,1)$);
		\coordinate (M3) at ($(N3) + (0,0,1)$);
		\coordinate (M4) at ($(N4) + (0,0,1)$);
		%
		\coordinate (alphabottom) at (0.5, 0.5, 0.5);
		\coordinate (alphatop) at (0.5, 0.5, 1.5);
		%
		\fill [green!50,opacity=0.445] (P) -- (alphabottom) -- (R);
		\fill [green!60,opacity=0.445] (R) -- (alphabottom) -- ($(R3) + (0,0,1)$) -- (R3);
			 \fill [green!80!white,opacity=0.345] (P) -- (alphabottom) -- (Rt) -- (A2) -- (O1) -- (N1);
		\fill [green!50,opacity=0.345] (alphabottom) -- (Rt) -- (A2) -- (A3) -- (alphatop) -- (Pt);
		\fill [green!80,opacity=0.345] (alphatop) -- (A4) -- (A3);
		\fill[color=blue!60!black] (0.3, 0.79, 1.89) circle (0pt) node (0up) { {\scriptsize$A$} };
		%
			 \fill [green!50!white,opacity=0.445] (R) -- (alphabottom) -- (Rt) -- (A2) -- (O2) -- (N2);
		\draw[color=green!30!black, ultra thick, rounded corners=0.5mm, postaction={decorate}, decoration={markings,mark=at position .17 with {\arrow[draw=green!30!black]{>}}}] (alphabottom) -- (Rt) -- (A2);
		\fill[color=blue!60!black] (0.6, 0.5, 2) circle (0pt) node (0up) { {\scriptsize$A$} };
		\fill [green!90!magenta,opacity=0.545] (R) -- (alphabottom) -- (Pt) -- (alphatop) -- ($(R3) + (0,0,1)$)-- (R3);
			\fill [green!80!white,opacity=0.345] (R3) -- ($(R3) + (0,0,1)$) -- (alphatop) -- (A3) -- (O3) -- (N3);
		\draw[string=green!30!black, ultra thick] (alphatop) -- (A3); 	
		%
			 \fill [green!70!white,opacity=0.345] (R3) -- ($(R3) + (0,0,1)$) -- (alphatop) -- (A4) -- (O4) -- (N4);
		\fill [green!90!magenta,opacity=0.345] (L) -- (P) -- (alphabottom) -- (Pt) -- (alphatop) -- (Po) -- (Lo);
		%
		\draw[string=green!30!black, ultra thick] (P) -- (alphabottom);
		\draw[string=green!30!black, ultra thick] (R) -- (alphabottom);
		\draw[string=green!30!black, ultra thick] (R3) -- ($(R3) + (0,0,1)$);
		\draw[string=green!30!black, ultra thick] ($(R3) + (0,0,1)$) -- (alphatop);
		\draw[string=green!30!black, ultra thick] (alphatop) -- (A4);
		\draw[color=green!30!black, ultra thick, rounded corners=0.5mm, postaction={decorate}, decoration={markings,mark=at position .33 with {\arrow[draw=green!30!black]{>}}}] (alphabottom) -- (Pt) -- (alphatop);
		\draw[color=green!30!black, ultra thick, rounded corners=0.5mm] (R3) --  ($(R3) + (0,0,1)$) -- (alphatop) ;
		%
		\fill[color=green!30!black] (alphabottom) circle (1.2pt) node[left] (0up) { {\scriptsize$\alpha$} };
		\fill[color=green!30!black] (alphatop) circle (1.2pt) node[left] (0up) { {\scriptsize$\alpha$} };
		%
		\draw [black,opacity=0.4141, densely dotted, semithick] (Lt) -- (Pt) -- (0.75, 1, 1) -- (M4);
		\draw [black,opacity=0.4141, densely dotted, semithick] (0.75, 1, 1) -- (M3);
		\draw [black,opacity=0.4141, densely dotted, semithick] (Pt) -- (Rt) -- (M2);
		\draw [black,opacity=0.4141, densely dotted, semithick] (Rt) -- (M1);
		%
		\draw [black!80!white,opacity=1, thin] (R) -- (P);
		\draw [black,opacity=1, very thin] (A2) -- (A3);
		\draw [black,opacity=1, very thin] (A4) -- (A3);
		\draw [black,opacity=1, very thin] (A2) -- (O2) -- (N2) -- (R);
		\draw [black,opacity=1, very thin] (A3) -- (O3) -- (N3) -- (R3);
		\draw [black,opacity=1, very thin] (A4) -- (O4) -- (N4) -- (R3);
		\draw [black,opacity=1, very thin] (R3) -- (R);
		\draw [black,opacity=1, very thin] (L) -- (P);
		\draw [black,opacity=1, very thin] (Po) -- (Lo) -- (L);
		\draw [black,opacity=1, very thin] (R3) -- (R);
		\draw [black,opacity=1, very thin] (A2) -- (O1) -- (N1) -- (P);
		\end{tikzpicture} = 
        \begin{tikzpicture}[thick, xscale=1.77, yscale=1.1, color=green!60!black, baseline=1cm, >=stealth, 
		style={x={(-0.6cm,-0.4cm)},y={(1cm,-0.2cm)},z={(0cm,0.9cm)}}]
		\pgfmathsetmacro{\yy}{0.2}
		\coordinate (P) at (0.5, \yy, 0);
		\coordinate (R) at (0.625, 0.5 + \yy/2, 0);
		\coordinate (L) at (0.5, 0, 0);
		\coordinate (R1) at (0.25, 1, 0);
		\coordinate (R2) at (0.5, 1, 0);
		\coordinate (R3) at (0.75, 1, 0);
		\coordinate (N1) at ($(R1) + (R2) - (R)$);
		\coordinate (N2) at ($(R) + 2*(R2) - 2*(R)$);
		\coordinate (N3) at ($(R3) + (R2) - (R)$);
		\coordinate (N4) at ($2*(R3) - (R)$);
		\coordinate (1P) at ($(P) + (0,0,1)$);
		\coordinate (1L) at (0.5, 0, 1);
		\coordinate (1N1) at ($(N1) + (0,0,1)$);
		\coordinate (1N2) at ($(N2) + (0,0,1)$);
		\coordinate (1N3) at ($(N3) + (0,0,1)$);
		\coordinate (1N4) at ($(N4) + (0,0,1)$);
		\coordinate (1R) at ($(R) + (0,0,1)$);
		\coordinate (1RR) at ($(R3) + (-0.25,0,1)$);
		\coordinate (2P) at ($(P) + (0,0,2)$);
		\coordinate (2L) at (0.5, 0, 2);
		\coordinate (2N1) at ($(N1) + (0,0,2)$);
		\coordinate (2N2) at ($(N2) + (0,0,2)$);
		\coordinate (2N3) at ($(N3) + (0,0,2)$);
		\coordinate (2N4) at ($(N4) + (0,0,2)$);
		\coordinate (2RR) at ($(1RR) + (0,0,1)$);
		\coordinate (2R) at ($(2N3) - 2*(2N3) + 2*(2RR)$);
		\coordinate (3P) at ($(P) + (0,0,3)$);
		\coordinate (3L) at (0.5, 0, 3);
		\coordinate (3N1) at ($(N1) + (0,0,3)$);
		\coordinate (3N2) at ($(N2) + (0,0,3)$);
		\coordinate (3N3) at ($(N3) + (0,0,3)$);
		\coordinate (3N4) at ($(N4) + (0,0,3)$);
		\coordinate (3RR) at ($(2RR) + (-0.25,0,1)$);
		\coordinate (3R) at ($(2R) + (0,0,1)$);
		%
		\coordinate (alpha1) at (0.5, 0.75, 0.5);
		\coordinate (alpha2) at (0.5, 0.5, 1.5);
		\coordinate (alpha3) at (0.25, 0.75, 2.5);
		%
		\coordinate (psi) at (0.5, 0.75, 1.5);
		%
		\fill [green!90!magenta,opacity=0.345] (P) -- (1P) -- (alpha2) -- (2P) -- (3P) -- (3L) -- (L);
		\fill [green!80,opacity=0.345] (P) -- (1P) -- (alpha2) -- (2R) -- (alpha3) -- (3RR) -- (3N1) -- (N1);
			 \fill [green!80!white,opacity=0.345] (alpha3) -- (3RR) -- (3R);
			 \fill [green!20!white,opacity=0.345] (R) -- (alpha1) -- (1RR) -- (2RR) -- (alpha3) -- (3RR) -- (3N2) -- (N2);
		\fill[color=blue!60!black] (0.3, 0.79, 2.89) circle (0pt) node (0up) { {\scriptsize$A$} };
		\draw[string=green!30!black, ultra thick] (alpha3) -- (3RR); 	
		\fill [green!80,opacity=0.345] (R3) -- (alpha1) -- (1RR) -- (2RR) -- (alpha3) -- (3R) -- (3N3) -- (N3);
			 \fill [green!50!white,opacity=0.345] (R3) -- (alpha1) -- (R);
		\fill [green!60!blue,opacity=0.345] (alpha1) -- (1R) -- (alpha2) -- (2R) -- (alpha3) -- (2RR) -- (1RR);
		\draw[color=green!30!black, ultra thick, rounded corners=0.5mm, postaction={decorate}, decoration={markings,mark=at position .34 with {\arrow[draw=green!30!black]{>}}}] (R3) -- (alpha1) -- (1RR) -- (2RR) -- (alpha3);
			 \fill [green!70!white,opacity=0.345] (alpha2) -- (2P) -- (3P) -- (3R) -- (alpha3) -- (2R);
		\fill[color=blue!60!black] (0.6, 0.5, 3) circle (0pt) node (0up) { {\scriptsize$A$} };
		\draw[string=green!30!black, ultra thick] (alpha3) -- (3R); 		
		\fill[color=green!30!black] (alpha3) circle (1.2pt) node[left] (0up) { {\scriptsize$\alpha$} }; 	
		\draw[color=green!30!black, ultra thick, rounded corners=0.5mm, postaction={decorate}, decoration={markings,mark=at position .3 with {\arrow[draw=green!30!black]{>}}}] (alpha2) -- (2R) -- (alpha3);
		\fill[color=blue!60!black] (psi) circle (0pt) node (0up) { {\scriptsize$A$} };
		%
		\fill [green!80,opacity=0.345] (R3) -- (alpha1) -- (1R) -- (alpha2) -- (2P) -- (3P) -- (3N4) -- (N4);
		%
			 \fill [green!70!white,opacity=0.345] (P) -- (1P) -- (alpha2) -- (1R) -- (alpha1) -- (R);
		\draw[color=green!30!black, ultra thick, rounded corners=0.5mm, postaction={decorate}, decoration={markings,mark=at position .3 with {\arrow[draw=green!30!black]{>}}}] (P) -- (1P) -- (alpha2);
		\draw[color=green!30!black, ultra thick, rounded corners=0.5mm, postaction={decorate}, decoration={markings,mark=at position .3 with {\arrow[draw=green!30!black]{>}}}] (alpha1) -- (1R) -- (alpha2);
		\draw[color=green!30!black, ultra thick, rounded corners=0.5mm, postaction={decorate}, decoration={markings,mark=at position .3 with {\arrow[draw=green!30!black]{>}}}] (alpha2) -- (2P) -- (3P);
		\draw[string=green!30!black, ultra thick] (R) -- (alpha1); 		
		\draw[string=green!30!black, ultra thick] (R3) -- (alpha1); 		
		%
		\fill[color=green!30!black] (alpha1) circle (1.2pt) node[left] (0up) { {\scriptsize$\alpha$} };
		\fill[color=green!30!black] (alpha2) circle (1.2pt) node[left] (0up) { {\scriptsize$\alpha$} };
		%
		%
		\draw [black,opacity=0.4141, densely dotted, semithick] (2L) -- (2P) -- (2N4);
		\draw [black,opacity=0.4141, densely dotted, semithick] (2P) -- (2N1);
		\draw [black,opacity=0.4141, densely dotted, semithick] (2R) -- (2RR) -- (2N2);
		\draw [black,opacity=0.4141, densely dotted, semithick] (2RR) -- (2N3);
		\draw [black,opacity=0.4141, densely dotted, semithick] (1L) -- (1P) -- (1N4);
		\draw [black,opacity=0.4141, densely dotted, semithick] (1R) -- (1RR) -- (1N3);
		\draw [black,opacity=0.4141, densely dotted, semithick] (1RR) -- (1N2);
		\draw [black,opacity=0.4141, densely dotted, semithick] (1P) -- (1N1);
		%
		\draw [black,opacity=1, very thin] (3P) -- (3L) -- (L) -- (P); 
		\draw [black,opacity=1, very thin] (N1) -- (P); 
		\draw [black,opacity=1, very thin] (3RR) -- (3N1) -- (N1) -- (3N1); 
		\draw [black,opacity=1, very thin] (3RR) -- (3R); 
		\draw [black,opacity=1, very thin] (3RR) -- (3N2) -- (N2) -- (R); 
		\draw [black,opacity=1, very thin] (3R) -- (3N3) -- (N3) -- (R3); 
		\draw [black,opacity=1, very thin] (R3) -- (R); 
		\draw [black,opacity=1, very thin] (3P) -- (3RR);
		\draw [black,opacity=1, very thin] (3P) -- (3N4) -- (N4) -- (R); 
		\draw [black,opacity=1, very thin] (P) -- (R); 
		\end{tikzpicture}\, , \qquad
        \begin{tikzpicture}[thick, xscale=2.321, yscale=1.5, color=blue!50!black, baseline=0.2cm, >=stealth, 
		style={x={(-0.6cm,-0.4cm)},y={(1cm,-0.2cm)},z={(0cm,0.9cm)}}]
		\pgfmathsetmacro{\yy}{0.2}
		\coordinate (P) at (0.5, \yy, 0);
		\coordinate (R) at (0.625, 0.5 + \yy/2, 0);
		\coordinate (L) at (0.5, 0, 0);
		\coordinate (R1) at (0.25, 1, 0);
		\coordinate (R2) at (0.5, 1, 0);
		\coordinate (R3) at (0.75, 1, 0);
		\coordinate (Pt) at (0.5, \yy, 1);
		\coordinate (Rt) at (0.375, 0.5 + \yy/2, 1);
		\coordinate (Lt) at (0.5, 0, 1);
		\coordinate (R1t) at (0.25, 1, 1);
		\coordinate (R2t) at (0.5, 1, 1);
		\coordinate (R3t) at (0.75, 1, 1);
		\coordinate (alpha) at (0.5, 0.5, 0.5);
		%
		%
		\fill [green!50,opacity=0.545] (L) -- (P) -- (alpha) -- (Pt) -- (Lt);
		\fill [green!50,opacity=0.545] (Pt) -- (Rt) -- (alpha);
		\fill [green!50,opacity=0.545] (Rt) -- (R1t) -- (R1) -- (P) -- (alpha);
		%
		\fill [green!50,opacity=0.545] (R2) -- (R) -- (alpha) -- ($(alpha)+0.5*(Rt)-0.5*(alpha)$);
		\draw[green!80!black, ultra thick] (R2) -- ($(alpha)+0.5*(Rt)-0.5*(alpha)$);
		\fill[color=green!50!black] (0.5,1.03,0.28) circle (0pt) node[left] (0up) { {\scriptsize$u$} };
		\fill[color=green!30!black] ($(alpha)+0.5*(Rt)-0.5*(alpha)$) circle (1.2pt) node[right] (0up) { {\scriptsize$u^{\textrm{r}}$} };
		\draw[green!30!black, ultra thick] (alpha) -- ($(alpha)+0.5*(Rt)-0.5*(alpha)$);
		\fill [green!50,opacity=0.545] (Pt) -- (R3t) -- (R3) -- (R) -- (alpha);
		\fill [green!50,opacity=0.545] (P) -- (R) -- (alpha);
		%
		\draw[string=green!30!black, ultra thick] (P) -- (alpha);
		\draw[string=green!30!black, ultra thick] (R) -- (alpha);
		\draw[string=green!30!black, ultra thick] (alpha) -- (Pt);
		%
		\fill[color=green!30!black] (alpha) circle (1.2pt) node[left] (0up) { {\scriptsize$\alpha$} };
		\fill[color=green!30!black] (0.5,0.35,0.24) circle (0pt) node[left] (0up) { {\scriptsize$\mu$} };
		\fill[color=green!30!black] (0.5,0.72,0.21) circle (0pt) node[left] (0up) { {\scriptsize$\mu$} };
		\fill[color=green!30!black] (0.5,0.4,0.71) circle (0pt) node[left] (0up) { {\scriptsize$\mu$} };
		%
		\draw [black,opacity=1, very thin] (Pt) -- (Lt) -- (L) -- (P);
		\draw [black,opacity=1, very thin] (Pt) -- (Rt);
		\draw [black,opacity=1, very thin] (Rt) -- (R1t) -- (R1) -- (P);
		\draw [black,opacity=1, very thin] (R2) -- (R);
		\draw [black,opacity=1, very thin] (Pt) -- (R3t) -- (R3) -- (R);
		\draw [black,opacity=1, very thin] (P) -- (R);
		\end{tikzpicture} =
        \begin{tikzpicture}[thick, xscale=2.321, yscale=1.5, color=blue!50!black, baseline=0.2cm, >=stealth, 
		style={x={(-0.6cm,-0.4cm)},y={(1cm,-0.2cm)},z={(0cm,0.9cm)}}]
		\pgfmathsetmacro{\yy}{0.2}
		\coordinate (P) at (0.5, \yy, 0);
		\coordinate (R) at (0.625, 0.5 + \yy/2, 0);
		\coordinate (L) at (0.5, 0, 0);
		\coordinate (R1) at (0.25, 1, 0);
		\coordinate (R2) at (0.5, 1, 0);
		\coordinate (R3) at (0.75, 1, 0);
		\coordinate (Pt) at (0.5, \yy, 1);
		\coordinate (Rt) at (0.375, 0.5 + \yy/2, 1);
		\coordinate (Lt) at (0.5, 0, 1);
		\coordinate (R1t) at (0.25, 1, 1);
		\coordinate (R2t) at (0.5, 1, 1);
		\coordinate (R3t) at (0.75, 1, 1);
		\coordinate (alpha) at (0.5, 0.5, 0.5);
		%
		%
		\fill [green!50,opacity=0.545] (L) -- (P) -- (Pt) -- (Lt);
		\fill [green!50,opacity=0.545] (P) -- (Pt) -- (R1t) -- (R1);
		%
		\fill [green!50,opacity=0.545] (R2) -- (R) -- ($(alpha)+0.5*(R)-0.5*(alpha)$);
		\draw[green!80!black, ultra thick] (R2) -- ($(alpha)+0.5*(R)-0.5*(alpha)$);
		\fill[color=green!50!black] (0.5,0.8,0.23) circle (0pt) node[left] (0up) { {\scriptsize$u$} };
		\fill [green!50,opacity=0.545] (Pt) -- (R3t) -- (R3) -- (R) -- (P);
		%
		\draw[string=green!30!black, ultra thick] (P) -- (Pt);
		\fill[color=green!30!black] ($(alpha)+0.5*(R)-0.5*(alpha)$) circle (1.2pt) node[above] (0up) { {\scriptsize$u^{\textrm{l}}$} };
		\draw[green!30!black, ultra thick] (R) -- ($(alpha)+0.5*(R)-0.5*(alpha)$);
		%
		\fill[color=green!30!black] (0.5,0.22,0.5) circle (0pt) node[left] (0up) { {\scriptsize$\mu$} };
		%
		\draw [black,opacity=1, very thin] (Pt) -- (Lt) -- (L) -- (P);
		\draw [black,opacity=1, very thin] (Pt) -- (Rt);
		\draw [black,opacity=1, very thin] (Rt) -- (R1t) -- (R1) -- (P);
		\draw [black,opacity=1, very thin] (R2) -- (R);
		\draw [black,opacity=1, very thin] (Pt) -- (R3t) -- (R3) -- (R);
		\draw [black,opacity=1, very thin] (P) -- (R);
		\end{tikzpicture}\, ,
    \end{equation}
    which are the pentagon axiom for the associator, and the triangle axiom for the interaction of the unitors and associator, here expressed in the three-dimensional graphical calculus.
\end{definition}

\begin{remark}\label{rem:E1_coalgebra}
    We can also have the dual version of the data above, namely, we can define an $\mathbb{E}_1$-coalgebra in a Gray category $\mathscr{C}$, as an $\mathbb{E}_1$-algebra in the Gray category $\mathscr{C}^{\text{2-op}}$, which has the 2-cells reversed.
\end{remark}

\begin{definition}[{\cite[Def.\,4.1]{orbifold3d}}]\label{def:special_orbifold_data}
    A \emph{special orbifold datum} in a Gray category with duals is an $\mathbb{E}_1$-algebra $(a, A, \mu, u, \alpha, u^{\mathrm{l}}, u^{\mathrm{r}})$ that further satisfies the relations
    \begin{equation}
.
    \end{equation}
    In the above relations, we used the notation
    \begin{equation}
        \alpha' := %
,
    \end{equation}
    where the surfaces have been omitted to avoid clutter. Bending surfaces is done through the $\vee$-dualization data.
\end{definition}

\begin{remark}
    The striped surfaces in the digrams have been used to denote the $\#$-dual. In this case, since $A^\#$ and $A$ are isomorphic, they are there to indicate where dualization data has been used rather than to distinguish the actual 1-morphisms. The directions on the $\mu$-lines help with this distinction as well, because they indicate whether we're working with $\mu$ or its dual $\mu^\vee$.
\end{remark}

\begin{remark}\label{rem:coalgebra_from_special_orbifold}
    Observe that a special orbifold datum $(a, A, \mu, u, \alpha, u^{\mathrm{l}}, u^{\mathrm{r}})$ in a Gray category with duals $\mathscr{C}$ also carries an $\mathbb{E}_1$-coalgebra structure (as discussed in \Cref{rem:E1_coalgebra}), where the comultiplication is $\mu^\vee: A \to A \boxtimes A$, the counit is $u^\vee: A \to 1_a$, and the coassociator and counitors are $(\alpha^\vee)^{-1}$, $(u^{\mathrm{l} \vee})^{-1}$ and $(u^{\mathrm{r} \vee})^{-1}$, respectively. The axioms of an $\mathbb{E}_1$-coalgebra are then satisfied because they are the images of the ones for an $\mathbb{E}_1$-algebra under the functor $-^\vee$, with an additional inversion of every 3-cell.
\end{remark}

\begin{lemma}\label{lem:duality_for_special_orbifold}
    Let $(a, A, \mu, u, \alpha, u^{\mathrm{l}}, u^{\mathrm{r}})$ be a special orbifold datum in a Gray category with duals. Then there are unique isomorphisms
    \begin{equation}
        \begin{tikzpicture}[very thick, scale=0.6, color=green!50!black, baseline]
            \tktrival{1}{1}{0}{0}{1}{1}
            \tkarc{-1}{-1}{0}{0}{1}{1}
            \draw (-1, 0) -- (-1, -0.6);
            \draw[-dot-] (1, 0) -- (1, 1);
        \end{tikzpicture}\ \xrightarrow{\ \ \cong \ \ } \ 
        \begin{tikzpicture}[very thick, scale=0.8, color=green!50!black, baseline=-0.4cm]
            \tktrival{-1}{1}{0}{0}{1}{1}
        \end{tikzpicture}\, , \qquad \qquad
        \begin{tikzpicture}[very thick, xscale=-0.6, yscale=0.6, color=green!50!black, baseline]
            \tktrival{1}{1}{0}{0}{1}{1}
            \tkarc{-1}{-1}{0}{0}{1}{1}
            \draw (-1, 0) -- (-1, -0.6);
            \draw[-dot-] (1, 0) -- (1, 1);
        \end{tikzpicture} \ \xrightarrow{\ \ \cong \ \ } \ 
        \begin{tikzpicture}[very thick, scale=0.8, color=green!50!black, baseline=-0.4cm]
            \tktrival{-1}{1}{0}{0}{1}{1}
        \end{tikzpicture}\, ,
    \end{equation}
    where the dot represents the unique isomorphism between $A$ and its $\#$-dual that makes the diagram well-typed.
\end{lemma}

\begin{proof}
    Along with the usual multiplication $\mu$ and unit $u$, by \Cref{rem:coalgebra_from_special_orbifold} we also have a comultiplication $\mu^\vee$ and counit $u^\vee$. Consider the 2-morphisms
    \begin{equation}
        \begin{tikzpicture}[very thick, scale=0.8, color=green!50!black, baseline]
            \tktrival{1}{0}{0}{0.5}{1}{1}
            \draw (0, 0.5) node[Odot] {};
        \end{tikzpicture}\, , \qquad \qquad
        \begin{tikzpicture}[very thick, xscale=0.8, yscale=-0.8, color=green!50!black, baseline]
            \tktrival{1}{0}{0}{0.5}{1}{1}
            \draw (0, 0.5) node[Odot] {};
        \end{tikzpicture}\, .
    \end{equation}
    They satisfy the adjunction triangle identities
    \begin{equation}
        \begin{tikzpicture}[very thick, scale=0.6, color=green!50!black, baseline]
            \tktrival{-1}{-1}{0}{0}{1}{1}
            \tktrival{1}{1}{0}{0}{1}{1}
            \draw (0.5, -1) node[Odot] {};
            \draw (-0.5, 1) node[Odot] {};
            \draw (-1, 0) -- (-1, -1);
            \draw (1, 0) -- (1, 1);
        \end{tikzpicture}\ \xrightarrow[\ \ \cong \ \ ]{\overline{\alpha}''} \ 
        \begin{tikzpicture}[very thick, scale=0.6, color=green!50!black, baseline]
            \tktrival{-1}{0}{0}{0}{1}{1}
            \tktrival{1}{0}{0}{0}{1}{1}
            \draw (0.5, -1) node[Odot] {};
            \draw (-0.5, 1) node[Odot] {};
        \end{tikzpicture}\ \xrightarrow[\ \ \cong \ \ ]{(u^{\mathrm{l} \vee})^{-1}} \ 
        \begin{tikzpicture}[very thick, scale=0.6, color=green!50!black, baseline]
            \tktrival{1}{0}{0}{0}{1}{1}
            \draw (0.5, -1) node[Odot] {};
            \draw (0, 0) -- (0, 1);
        \end{tikzpicture}\ \xrightarrow[\ \ \cong \ \ ]{u^{\mathrm{r} \vee}} \ 
        \begin{tikzpicture}[very thick, scale=0.6, color=green!50!black, baseline]
            \draw (0, -1) -- (0, 1);
        \end{tikzpicture}\  ,
    \end{equation}
    where the symbolic 3-cells names stand for the appropriate whiskerings which make the diagram well-typed. Since adjunction data are unique up to unique isomorphism, these observations provide both the fact that $A^\#$ is uniquely isomorphic to $A$, and the evaluation and coevaluation that participate in this adjunction. With this the required 3-isomorphisms from the lemma statement are well-typed.

    For $\#$-dualization on the right, to actually provide an isomorphism we have
    \begin{equation}
        \begin{tikzpicture}[very thick, scale=0.6, color=green!50!black, baseline]
            \tktrival{1}{1}{0}{0}{1}{1}
            \tkarc{-1}{-1}{0}{0}{1}{1}
            \draw (-1, 0) -- (-1, -0.6);
            \draw[-dot-] (1, 0) -- (1, 1);
        \end{tikzpicture}\ \xrightarrow[\ \ \cong \ \ ] \ 
        \begin{tikzpicture}[very thick, scale=0.5, color=green!50!black, baseline=-0.1cm]
            \tktrival{-1}{-1}{0}{0}{1}{1}
            \tktrival{1}{1}{0}{0}{1}{1}
            \draw (0.5, -1) node[Odot] {};
            \draw (-1, 0) -- (-1, -1);
            \draw (1, 0) -- (1, 1);
        \end{tikzpicture}\ \xrightarrow[\ \ \cong \ \ ]{\overline{\alpha}''} \ 
        \begin{tikzpicture}[very thick, scale=0.5, color=green!50!black, baseline=-0.1cm]
            \tktrival{-1}{0}{0}{0}{1}{1}
            \tktrival{1}{0}{0}{0}{1}{1}
            \draw (0.5, -1) node[Odot] {};
        \end{tikzpicture}\ \xrightarrow[\ \ \cong \ \ ]{u^{\mathrm{l}}} \ 
        \begin{tikzpicture}[very thick, scale=0.8, color=green!50!black, baseline=-0.5cm]
            \tktrival{-1}{1}{0}{0}{1}{1}
        \end{tikzpicture} := \mu^\vee.
    \end{equation}
    Since the end result is $\mu^\vee$, an adjoint itself, the above isomorphism must be unique. By a similar logic, we get the result for $\#$-dualization on the left as well.
\end{proof}

\begin{proof}[Proof of \Cref{prop:3_bistellar_from_special_orbifold}]
    To provide the functor from $\mathscr{P}_3$, we just need to provide the data of where the generators are mapped to and check that they satisfy the needed relations. Using the notation of \Cref{def:special_orbifold_data} (which is also the notation of \cite[Def.\,4.1]{orbifold3d}) on the right-hand side of the equalities, and the notation for the generators in \Cref{thm:3_bistellar_datum} on the left-hand side, we set
    \begin{equation}
        a := a, \qquad \qquad A := A, \qquad \qquad \mu := \mu, \qquad \qquad \delta := \mu^\vee
    \end{equation}
    for the generators up to dimension 2, where $A = A^\#$ is self-dual, with $-^\#$ being the dualization on 1-morphisms, and $-^\vee$ being the dualization on 2-morphisms. For the invertible generator 3-cells we have
    \begin{equation}
        \alpha := \alpha, \qquad \qquad \hat{\gamma} := \alpha^\vee, \qquad \qquad \phi_L := \overline{\alpha}'', \qquad \qquad \phi_R := \alpha''^\vee, \notag
    \end{equation}
    \begin{equation}
        \hat{\beta}_I^U := (\mathrm{ev}_{\mu} \otimes 1_\mu) \circ (1_\mu \otimes \phi_I), \qquad \qquad \hat{\beta}_I^D := (1_{\mu^\vee} \otimes \mathrm{ev}_{\mu}) \circ (\phi_I \otimes 1_{\mu^\vee}),
    \end{equation}
    where $I \in \{L, R\}$, and $\mathrm{ev}_{\mu}$ is the invertible 3-cell
    \begin{equation}
        \begin{tikzpicture}[very thick, scale=0.7 * 0.9, color=green!50!black, baseline]
            \tktrival{-1}{0}{0}{-0.75}{1}{1}
            \tktrival{1}{0}{0}{1.25}{1}{1}    
        \end{tikzpicture}
        \ \ \xrightarrow[\qquad]{\mathrm{ev}_\mu} \ \ 
        \begin{tikzpicture}[very thick, scale=0.7 * 2/3, color=green!50!black, baseline]
            \draw (0, -1) -- (0, 1.5);
        \end{tikzpicture}
    \end{equation}
    that makes an appearance in \Cref{eq:condition_O8}. The conditions that special orbifold data satisfies ensure that these 3-cells are all invertible, so that they are guaranteed to land in $t_{(3,2)} \mathscr{C}$. The remaining generators are the appropriate inverses of the ones above. Notice that the bubbleators are defined through the Frobeniusators, along with $\mathrm{ev}_\mu$.

    In checking the relations it is helpful to observe that because of their pentagonal shape, there is a certain pattern of cyclic, or reflected cyclic, permutations on the vertices of the pentagon. This is not quite an action of $\mathbb{Z}_5$, because, in general, it does not respect the 3-cells. But, by construction, the 3-cells given by special orbifold data are invertible, and the plain and anti version of the generators are inverses of each other. This means that two relations that are related to one another by one of these permutations, are not independent relations for this particular setup of invertible 3-cells, specifically because of the cyclic property.

    By \cite[Lem.\,3.16]{orbifoldnd} every 1-4 move is a concatenation of a 2-3 move and a special type of bubble move, which can be seen in \cite[Eq.\,3.61]{orbifoldnd}. In our notation this move is given by
    \begin{equation}
        \begin{tikzpicture}[scale=0.5]
            \def\R{4cm}
            \def\clearance{2mm}
            \pgfmathsetlengthmacro{\shorten}
                {(1cm+\clearance)/cos(36)}

            \coordinate (P1) at (162:\R);
            \coordinate (P2) at (90:\R);
            \coordinate (P3) at (18:\R);

            \draw[triple] ($(P1)!\shorten!(P2)$) -- ($(P2)!\shorten!(P1)$);
            \draw[triple] ($(P2)!\shorten!(P3)$) -- ($(P3)!\shorten!(P2)$);
            \draw[double] ($(P1)!\shorten!(P3)$) -- ($(P3)!\shorten!(P1)$);

            \node at (P1) {\resizebox{0.7cm}{0.7cm}{\begin{tikzpicture}[very thick, color=green!50!black, yscale=-1]
                \tktrival{1}{0}{0}{0}{1}{1}
            \end{tikzpicture}}};
            \node at (P2) {\resizebox{1.2cm}{1.2cm}{\begin{tikzpicture}[very thick, color=green!50!black, yscale=-1]
                \tktrival{1}{0}{0}{3}{3}{3}
                \tktrival{-1}{1}{1.5}{0}{1}{1}
                \tktrival{1}{0}{0}{1}{1}{1}
                \tktrival{-1}{-1}{-1.5}{0}{1}{1}
                \draw (0, 1) node[Odot] {};
            \end{tikzpicture}}};
            \node at (P3) {\resizebox{0.7cm}{0.7cm}{\begin{tikzpicture}[very thick, color=green!50!black, yscale=-1]
                \tktrival{1}{0}{0}{0}{1}{1}
            \end{tikzpicture}}};
        \end{tikzpicture},
    \end{equation}
    where in the middle of the bubble we use the coevaluation for the self dual $A$, which we can also think of as the data of the unit (as seen in \Cref{lem:duality_for_special_orbifold}). While the 2-flipped version is generically an independent move, here it is implied by making use of the $\vee$-duality. There is also an immediate corollary, that the given 1-4 move contains one more bubble than the 2-3 move. Since special orbifold data satisfy the relation imposed by the special bubble move, as proved in \cite[Lem.\,3.17]{orbifoldnd}, all 1-4 moves are automatically satisfied if the 2-3 moves with fewer bubbles than them are satisfied.
    
    The previous two observations, by induction, imply that all we need to check are the relations in the marked groups that have no 1-4 move, and no bubbles (which are required to be, and in fact are, the same condition). This leaves only the marked groups 0, 3, and 13 represented by
    \begin{center}
        \begin{tikzpicture}[scale=0.5]

            \def\R{4cm}
            \def\clearance{2mm}
            \pgfmathsetlengthmacro{\shorten}
                {(1cm+\clearance)/cos(36)}
            \def\ang{20}
            \def\ext{0.9cm}

            \coordinate (P1) at (162:\R);
            \coordinate (P2) at (90:\R);
            \coordinate (P3) at (18:\R);
            \coordinate (P4) at (-126:\R);
            \coordinate (P5) at (-54:\R);
            \coordinate (L) at (132:{\R + 2*\ext});
            \coordinate (C1) at (-135:{sqrt(2)*(\R+\ext)});
            \coordinate (C2) at (-180+\ang : {(\R+\ext)/cos(\ang)});
            \coordinate (C3) at (-90-\ang : {(\R+\ext)/cos(\ang)});

            \draw[color=red!40!yellow, very thick] (C2) -- (C1) -- (C3);
            \node[draw, circle] at ($(C1) + (0.5cm,0.5cm)$) [color=red!40!yellow, very thick, inner sep=0.2ex] {0};

            \draw[triple] ($(P1)!\shorten!(P2)$) -- ($(P2)!\shorten!(P1)$);
            \draw[triple] ($(P2)!\shorten!(P3)$) -- ($(P3)!\shorten!(P2)$);
            \draw[triple] ($(P1)!\shorten!(P4)$) -- ($(P4)!\shorten!(P1)$);
            \draw[triple] ($(P4)!\shorten!(P5)$) -- ($(P5)!\shorten!(P4)$);
            \draw[triple] ($(P5)!\shorten!(P3)$) -- ($(P3)!\shorten!(P5)$);
            \draw[quadruple,->] (0,0.45cm) -- (0,-0.65cm);

            \node at (P1) {\resizebox{1cm}{1cm}{\begin{tikzpicture}[very thick, color=green!50!black, yscale=-1]
                \tktrivalcoord{2.0}{0}{3.5}{3}{1.5}{1}
                \tktrivalcoord{1.5}{1}{2.5}{3}{1.0}{2}
                \tktrivalcoord{1.0}{2}{0.5}{3}{1.5}{3}
            \end{tikzpicture}}};
            \node at (P2) {\resizebox{1cm}{1cm}{\begin{tikzpicture}[very thick, color=green!50!black, yscale=-1]
                \tktrivalcoord{2.0}{0}{3.0}{1}{1.0}{2}
                \tktrivalcoord{3.0}{1}{3.5}{3}{2.5}{3}
                \tktrivalcoord{1.0}{2}{0.5}{3}{1.5}{3}
            \end{tikzpicture}}};
            \node at (P3) {\resizebox{1cm}{1cm}{\begin{tikzpicture}[very thick, color=green!50!black, yscale=-1]
                \tktrivalcoord{2.0}{0}{2.5}{1}{0.5}{3}
                \tktrivalcoord{2.5}{1}{3.0}{2}{1.5}{3}
                \tktrivalcoord{3.0}{2}{3.5}{3}{2.5}{3}
            \end{tikzpicture}}};
            \node at (P4) {\resizebox{1cm}{1cm}{\begin{tikzpicture}[very thick, color=green!50!black, yscale=-1]
                \tktrivalcoord{2.0}{0}{3.5}{3}{1.5}{1}
                \tktrivalcoord{1.5}{1}{2.0}{2}{0.5}{3}
                \tktrivalcoord{2.0}{2}{2.5}{3}{1.5}{3}
            \end{tikzpicture}}};
            \node at (P5) {\resizebox{1cm}{1cm}{\begin{tikzpicture}[very thick, color=green!50!black, yscale=-1]
                \tktrivalcoord{2.0}{0}{2.5}{1}{0.5}{3}
                \tktrivalcoord{2.5}{1}{3.5}{3}{2.0}{2}
                \tktrivalcoord{2.0}{2}{2.5}{3}{1.5}{3}
            \end{tikzpicture}}};

        \end{tikzpicture}\, ,
        \ \ 
        \begin{tikzpicture}[scale=0.5]

            \def\R{4cm}
            \def\clearance{2mm}
            \pgfmathsetlengthmacro{\shorten}
                {(1cm+\clearance)/cos(36)}
            \def\ang{20}
            \def\ext{0.9cm}

            \coordinate (P1) at (162:\R);
            \coordinate (P2) at (90:\R);
            \coordinate (P3) at (18:\R);
            \coordinate (P4) at (-126:\R);
            \coordinate (P5) at (-54:\R);
            \coordinate (L) at (132:{\R + 2*\ext});
            \coordinate (C1) at (-135:{sqrt(2)*(\R+\ext)});
            \coordinate (C2) at (-180+\ang : {(\R+\ext)/cos(\ang)});
            \coordinate (C3) at (-90-\ang : {(\R+\ext)/cos(\ang)});

            \draw[color=red!40!yellow, very thick] (C2) -- (C1) -- (C3);
            \node[draw, circle] at ($(C1) + (0.5cm,0.5cm)$) [color=red!40!yellow, very thick, inner sep=0.2ex] {3};

            \draw[triple] ($(P1)!\shorten!(P2)$) -- ($(P2)!\shorten!(P1)$);
            \draw[triple] ($(P2)!\shorten!(P3)$) -- ($(P3)!\shorten!(P2)$);
            \draw[triple] ($(P1)!\shorten!(P4)$) -- ($(P4)!\shorten!(P1)$);
            \draw[triple] ($(P4)!\shorten!(P5)$) -- ($(P5)!\shorten!(P4)$);
            \draw[triple] ($(P5)!\shorten!(P3)$) -- ($(P3)!\shorten!(P5)$);
            \draw[quadruple,->] (0,0.45cm) -- (0,-0.65cm);

            \node at (P1) {\resizebox{1cm}{1cm}{\begin{tikzpicture}[very thick, color=green!50!black, yscale=-1]
                \tktrivalcoord{1.5}{0}{2.5}{1.0}{1.0}{2}
                \tktrivalcoord{3.5}{0}{2.5}{1.0}{3.0}{2}
                \tktrivalcoord{1.0}{2}{0.5}{3}{1.5}{3}
            \end{tikzpicture}}};
            \node at (P2) {\resizebox{1cm}{1cm}{\begin{tikzpicture}[very thick, color=green!50!black, yscale=-1]
                \tktrivalcoord{1.5}{0}{2.0}{1}{0.5}{3}
                \tktrivalcoord{2.0}{1}{2.5}{2.0}{1.5}{3}
                \tktrivalcoord{3.5}{0}{2.5}{2.0}{3.0}{3}
            \end{tikzpicture}}};
            \node at (P3) {\resizebox{1cm}{1cm}{\begin{tikzpicture}[very thick, color=green!50!black, yscale=-1]
                \tktrivalcoord{1.5}{0}{2.0}{1.0}{0.5}{3}
                \tktrivalcoord{3.5}{0}{2.0}{1.0}{2.5}{2.0}
                \tktrivalcoord{2.5}{2.0}{3.0}{3}{1.5}{3}
            \end{tikzpicture}}};
            \node at (P4) {\resizebox{1cm}{1cm}{\begin{tikzpicture}[very thick, color=green!50!black, yscale=-1]
                \tktrivalcoord{3.5}{0}{1.5}{0}{2.0}{1.0}
                \tktrivalcoord{2.0}{1.0}{3.0}{3}{1.0}{2}
                \tktrivalcoord{1.0}{2}{0.5}{3}{1.5}{3}
            \end{tikzpicture}}};
            \node at (P5) {\resizebox{1cm}{1cm}{\begin{tikzpicture}[very thick, color=green!50!black, yscale=-1]
                \tktrivalcoord{3.5}{0}{1.5}{0}{2.0}{1.0}
                \tktrivalcoord{2.0}{1.0}{2.5}{2}{0.5}{3}
                \tktrivalcoord{2.5}{2}{3.0}{3}{1.5}{3}
            \end{tikzpicture}}};

        \end{tikzpicture}\, ,
        \ \ 
        \begin{tikzpicture}[scale=0.5]

            \def\R{4cm}
            \def\clearance{2mm}
            \pgfmathsetlengthmacro{\shorten}
                {(1cm+\clearance)/cos(36)}
            \def\ang{20}
            \def\ext{0.9cm}

            \coordinate (P1) at (162:\R);
            \coordinate (P2) at (90:\R);
            \coordinate (P3) at (18:\R);
            \coordinate (P4) at (-126:\R);
            \coordinate (P5) at (-54:\R);
            \coordinate (L) at (132:{\R + 2*\ext});
            \coordinate (C1) at (-135:{sqrt(2)*(\R+\ext)});
            \coordinate (C2) at (-180+\ang : {(\R+\ext)/cos(\ang)});
            \coordinate (C3) at (-90-\ang : {(\R+\ext)/cos(\ang)});

            \draw[color=red!40!yellow, very thick] (C2) -- (C1) -- (C3);
            \node[draw, circle] at ($(C1) + (0.5cm,0.5cm)$) [color=red!40!yellow, very thick, inner sep=-0.05ex] {13};

            \draw[triple] ($(P1)!\shorten!(P2)$) -- ($(P2)!\shorten!(P1)$);
            \draw[triple] ($(P2)!\shorten!(P3)$) -- ($(P3)!\shorten!(P2)$);
            \draw[triple] ($(P1)!\shorten!(P4)$) -- ($(P4)!\shorten!(P1)$);
            \draw[triple] ($(P4)!\shorten!(P5)$) -- ($(P5)!\shorten!(P4)$);
            \draw[triple] ($(P5)!\shorten!(P3)$) -- ($(P3)!\shorten!(P5)$);
            \draw[quadruple,->] (0,0.45cm) -- (0,-0.65cm);

            \node at (P1) {\resizebox{1cm}{1cm}{\begin{tikzpicture}[very thick, color=green!50!black, yscale=1]
                \tktrivalcoord{3.0}{0}{3.5}{3}{2.5}{1.0}
                \tktrivalcoord{1.0}{0}{2.5}{1.0}{1.5}{2.0}
                \tktrivalcoord{1.5}{2.0}{2.0}{3}{0.5}{3}
            \end{tikzpicture}}};
            \node at (P2) {\resizebox{1cm}{1cm}{\begin{tikzpicture}[very thick, color=green!50!black, yscale=1]
                \tktrivalcoord{3.0}{0}{3.5}{3}{2.5}{1.5}
                \tktrivalcoord{1.0}{1}{0.5}{3}{1.5}{2.0}
                \tktrivalcoord{2.5}{1.5}{1.5}{2.0}{2.0}{3}
            \end{tikzpicture}}};
            \node at (P3) {\resizebox{1cm}{1cm}{\begin{tikzpicture}[very thick, color=green!50!black, yscale=1]
                \tktrivalcoord{1.0}{0}{0.5}{3}{1.5}{1.0}
                \tktrivalcoord{3.0}{0}{1.5}{1.0}{2.5}{2.0}
                \tktrivalcoord{2.5}{2.0}{3.5}{3}{2.0}{3}
            \end{tikzpicture}}};
            \node at (P4) {\resizebox{1cm}{1cm}{\begin{tikzpicture}[very thick, color=green!50!black, yscale=1]
                \tktrivalcoord{3.0}{0}{1.0}{0}{2.0}{1.0}
                \tktrivalcoord{2.0}{1.0}{3.5}{3}{1.5}{2}
                \tktrivalcoord{1.5}{2}{2.0}{3}{0.5}{3}
            \end{tikzpicture}}};
            \node at (P5) {\resizebox{1cm}{1cm}{\begin{tikzpicture}[very thick, color=green!50!black, yscale=1]
                \tktrivalcoord{3.0}{0}{1.0}{0}{2.0}{1.0}
                \tktrivalcoord{2.0}{1.0}{2.5}{2}{0.5}{3}
                \tktrivalcoord{2.5}{2}{3.5}{3}{2.0}{3}
            \end{tikzpicture}}};

        \end{tikzpicture}\, .
    \end{center}

    Each of these three elements is a representative for its $\mathbb{Z}_2^4$-orbit. We clearly do not have to separately check relations that are related by a flip of the 4-cell or a flip of all the 3-cells, the second because the plain versions and anti versions of the generators are inverses of each other. We also don't have to separately check relations that are related by a flip of all 2-cells since this changes the multiplication into its dual, the comultiplication, and vice versa.

    For marked group 0, this contains the pentagon axiom for the associator, which is satisfied by \Cref{eq:3d_graphical_pentagon}. A flip of all 1-cells in this group is compensated by a flip of all 3-cells, so there is nothing else to check.

    In the marked group 3 we do not have enough symmetry, so there are two independent relations to check, which are related to each other by a flip of all 1-cells. On the other hand, in marked group 13, it is also the case that a flip of all 1-cells is compensated by a flip of all 3-cells, like in the case of marked group 0, so we only need to check the one representative pictured above. These three remaining relations are all satisfied through the same argument. Namely, using \Cref{lem:duality_for_special_orbifold} we can partially $\#$-dualize the output leg(s) connected to a comultiplication, which is either on the left or the right, for the relations in marked group 3, or on both the left and right, for the relation in marked group 13.

    More specifically, \Cref{lem:duality_for_special_orbifold} implies that on the hom-1-subcategory containing any one of the given diagrams partial $\vee$-dualization and partial $\#$-dualization of the leg(s) that connect(s) to a comultiplication are uniquely isomorphic, and because of the functoriality of the dualization functors, this also extends to the 3-cells involved.
    
    Since the Frobeniusators are defined by partially $\vee$-dualizing the associator on one side or the other, as seen in the definitions of $\alpha''$ and $\overline{\alpha}''$, it must be the case that under partial $\#$-dualization we also map the Frobeniusators to the associator, and vice versa. For example, the 3-cell
    \begin{equation}
        \begin{tikzpicture}[very thick, color=green!50!black, xscale=0.7, yscale=-0.7, baseline=-1cm]
            \tktrivalcoord{2.0}{0}{2.5}{1}{0.5}{3}
            \tktrivalcoord{2.5}{1}{3.0}{2}{1.5}{3}
            \tktrivalcoord{3.0}{2}{3.5}{3}{2.5}{3}
        \end{tikzpicture}\xrightarrow[\ \ \cong \ \ ] \ \ 
        \begin{tikzpicture}[very thick, color=green!50!black, xscale=0.7, yscale=-0.7, baseline=-1cm]
            \tktrivalcoord{1.5}{0}{2.0}{1}{0.5}{3}
            \tktrivalcoord{2.0}{1}{2.5}{2.0}{1.5}{3}
            \tktrivalcoord{3.5}{1}{2.5}{2.0}{3.0}{3}
            \tkarc{-1}{-1}{3.5}{1}{1}{1}
            \draw[-dot-] (4.5, 1) -- (4.5, 3);
        \end{tikzpicture}\ \xrightarrow[\ \ \cong \ \ ]{\overline{\alpha}''} \ 
        \begin{tikzpicture}[very thick, color=green!50!black, xscale=0.7, yscale=-0.7, baseline=-1cm]
            \tktrivalcoord{1.5}{0}{2.0}{1.0}{0.5}{3}
            \tktrivalcoord{3.5}{1}{2.0}{1.0}{2.5}{2.0}
            \tktrivalcoord{2.5}{2.0}{3.0}{3}{1.5}{3}
            \tkarc{-1}{-1}{3.5}{1}{1}{1}
            \draw[-dot-] (4.5, 1) -- (4.5, 3);
        \end{tikzpicture}\ \xrightarrow[\ \ \cong \ \ ] \ 
        \begin{tikzpicture}[very thick, color=green!50!black, xscale=0.7, yscale=-0.7, baseline=-1cm]
            \tktrivalcoord{2.0}{0}{2.5}{1}{0.5}{3}
            \tktrivalcoord{2.5}{1}{3.5}{3}{2.0}{2}
            \tktrivalcoord{2.0}{2}{2.5}{3}{1.5}{3}
        \end{tikzpicture}
    \end{equation}
    has to be the appropriate whiskering of the associator. In the above, we continue to use the dot to denote the unique isomorphism between $A$ and $A^\#$, and to label 3-cells symbolically, without acknowledging whiskering.

    The partial $\#$-dualization functor thus transforms all three relations, from marked groups 3 and 13, into some cyclic or reverse cyclic permutation of the pentagon axiom for the associator, and so all relations are satisfied.
\end{proof}

\begin{remark}
    It is generically, of course, not the case that every representation of $\mathscr{P}_3$ comes from a special orbifold datum. In fact, one should not expect this for a variety of reasons:
    \begin{itemize}
        \item As explained in \cite{orbifoldnd}, the special variant of 3-orbifold data is not as general as it could be, and in fact is an easier-to-work-with version of a more general set of relations.
        \item We can now also see that the restriction to landing in $t_{(3,2)} \mathscr{C}$ eliminates many examples, similarly to the case of \Cref{ex:P2withisosin2}, because there is no reason to impose this simplification except for practical purposes.
        \item Special orbifold data are designed for the case of oriented theories, which is why above we start with a Gray category with duals. Bistellar data are closest to the framed case since they do not see anything but the underlying weak $3$-category.
    \end{itemize}
\end{remark}

\begin{remark}[The 4-dimensional case]
    Since the definition of bistellar data works in any dimension, we already know what $\mathscr{P}_4$ is, and we even have a finite generators-and-relations presentation of it, again by virtue of its definition. However, this presentation has a lot of symmetries, and in any given special case it might be possible to optimize it and give a smaller presentation, as was the case with special orbifold data in 3-dimensions in \Cref{prop:3_bistellar_from_special_orbifold}.
    
    In the most general case though, \Cref{thm:satm_in_low_d} and \Cref{thm:satm_4} already provide the generators of which there are
    \begin{equation}
        1 \ \  + \ \  1 \ \  + \ \  2 \ \  + \ \  16 \ \  + \ \  1240 \ \  = \ \  1260.
    \end{equation}
    In this presentation the list of relations would consist of $\atm{5}^{\mathrm{fa}}$, which is expected to have many elements. One would also have to worry about checking the frame-acyclicity condition explicitly, since we are no longer guaranteed that it holds automatically.
\end{remark}

\appendix

\section{Finding Simplicial Atoms with Code}\label{app:pseudocode}

This appendix gives one possible approach to a computer program that will find all simplicial atoms in low dimensions. Here we describe the  logic and provide some pseudocode for the script, while an auxiliary \texttt{python} file gives an honest implementation \cite{pythoncode}. It's important to mention that there is no currently know decision procedure to determine whether a generic oriented graded poset is a molecule or not, so the strategy will be to whittle down the set of possible options to a low enough number that can be confirmed manually.

The set of all oriented graded poset structures on a graded poset is given by $2^{|\mathrm{Hasse \ edges}|}$, which for the $n$-simplex grows like $2^{2^n}$, so we have to be efficient in our encoding and our searching for the algorithm to be able to meaningfully work. We see the vertices of the $n$-simplex as the integers $1, \dots, n$, and encode the simplices, which are subsets of such integers using bitmasks. This allows for efficient computations with the simplices. The Hasse edges are encoded as a pair \texttt{(e, s)}, where \texttt{s} is a bitmask (representing a simplex) and \texttt{e} is an integer, such that at that position the bitmask \texttt{s} has the bit set (representing a vertex from the simplex). We read this as being the Hasse edge between the simplex represented by \texttt{s} and the simplex you get by removing the specified vertex represented by \texttt{e}.

Since the automorphisms of the $n$-simplex are $\mathrm{Aut}(\Delta^n) = S_n$, the permutations on $n$-elements, the code will find $n!$-many copies that are isomorphic, but not equal because they would differ by a permutation of the vertices. To deal with this redundancy, and since the encoding of the vertices already singles out a preferred ordering of the vertices, we seed the search by fixing all Hasse diagram edges that go between levels 0 and 1. That is, we fix all the directionality for the 1-dimensional edges of our oriented graded poset to be the standard one.

To find all simplicial atoms in low dimensions the idea is to search through all possible ways of equipping the poset that represents the standard simplex with an oriented graded poset structure, so that this structure is oriented thin \cite[Def.\,2.3.10]{hadzihasanovic}. More precisely, every molecule is an oriented thin, oriented graded poset once we adjoin a minimal element to it, as a result of \cite[Prop.\,5.3.4]{hadzihasanovic}, so what we really mean is that we only enforce the oriented thinness condition whenever we wouldn't need a minimal element for it to be true. This is done by the \textsc{Propagate} function. The main search, finding all oriented thin, oriented graded posets, builds the poset from the bottom up going level by level, which is enforced by the \textsc{ConstructEdgeOrder} function. Finally, we bring it all together by doing a depth-first search on the tree of possible oriented graded poset structures.

\vspace{0.15cm}
\noindent
\begin{minipage}[h]{0.46\linewidth}
    \begin{algorithmic}[1]
        \ttfamily
        \Function{GenerateOrientedThin}{$n$}
            \State $E$ $\gets$ \Call{ConstructEdgeOrder}{$n$}
            \State
            \Function{DFS}{$i$, $P$}
                \While{$i < |E|$ and $E_i \in P$}
                    \State $i \gets i + 1$
                \EndWhile
                \State

                \If{$i = |E|$}
                    \Yield{$P$}
                    \Return{}
                \EndIf
                \State

                \For{$\varepsilon \in \{-, +\}$}
                    \State $P'$ $\gets$ \Call{copy}{$P$}
                    \State $P'(E_i) \gets \varepsilon$
                    \State queue $\gets$ $\{E_i\}$
                    \If{\textbf{not} \Call{Propagate}{\tikzmarknode{call-propagate}{$P'$}, queue}}
                        \State \textbf{continue}
                    \EndIf
                    \Yield{\textbf{from} \Call{DFS}{$i$, $P'$}}
                \EndFor
            \EndFunction
            \State
            \State $P_0$ $\gets$ \Call{SeedInitialEdges}{\tikzmarknode{call-seed}{$n$}}
            \Yield{\textbf{from} \Call{DFS}{$0$, $P_0$}}
        \EndFunction
    \end{algorithmic}    
\end{minipage}\ \ \ \ \ \ 
\begin{minipage}[h]{0.5\linewidth}
    {\footnotesize
    \begin{algorithmic}[1]
        \ttfamily
        \Function{\tikzmarknode{def-propagate}{P}ropagate}{$P$, queue}
            \While{queue is nonempty}
                \State get edge $f$ from queue
                \If{$f$ is one top edge in diamond $\{f, g, f', g'\} \subset P$}
                    \If{all 4 edge signs known}
                        \State \Return{$P(f) P(g) P(f') P(g')$ is $-$}
                    \ElsIf{only $P(g)$ unknown}
                        \State $P(g) \gets - P(f) P(f') P(g')$
                        \State add $g$ to queue
                    \EndIf
                \EndIf
            \EndWhile
            \State \Return{True}
        \EndFunction
    \end{algorithmic}}
    \vspace{1cm}
    {\footnotesize
    \begin{algorithmic}[1]
        \ttfamily
        \Function{\tikzmarknode{def-seed}{S}eedInitialEdges}{$n$}
            \State $P$ $\gets$ $\emptyset$
            \For{$i \in [0, \dots, n]$}
                \For{$j \in [i + 1, n]$}
                    \State $P \big( \{i\} \subseteq \{i, j\} \big) \gets +$
                    \State $P \big( \{j\} \subseteq \{i, j\} \big) \gets -$ 
                \EndFor
            \EndFor
            \State \Return{$P$}
        \EndFunction
    \end{algorithmic}}
\end{minipage}
\begin{tikzpicture}[remember picture,overlay]
    \draw[->, >=Stealth, thick, black, shorten >=0.3cm, shorten <=2.2cm] ({def-propagate.west}) to[out=180,in=60,looseness=1.25] (call-propagate.east);

    \draw[->, >=Stealth, thick, black, shorten >=0.3cm, shorten <=2.2cm] (def-seed.west) to[out=180,in=40,looseness=1.25] (call-seed.east);
\end{tikzpicture}
\vspace{0.3cm}

Being oriented thin, though, is not enough for an oriented graded poset to be a molecule, so the next thing we check is globularity. We do this by checking that the $k$-iterated boundary is equal to the one-shot $k$-boundary, for every $0 \leq k < n$. This step eliminates certain candidates that have the correct `local' behavior, but have global deficiencies by having simplices that locally look like sources or targets, but appear outside the iterated $k$-boundary. This is done with a straightforward implementation of \Cref{def:boundary}.

Finally, there is a computationally more intensive check that we can do last, which is to recursively compute the slices of a molecule, and see whether they admit a topological sort of their maximal elements. If there is no topological sort on the maximal elements of $P$ then the function \textsc{TopologicalOrder} returns \texttt{None}, and we determine that $P$ is not a molecule.

\vspace{0.15cm}
\noindent
\begin{minipage}[t]{0.43\linewidth}
    \begin{algorithmic}[1]
        \ttfamily
        \Function{HasLayerings}{$P$}
            \State queue $\gets \{P\}$
            \While{queue is nonempty}
                \State get $R$ from queue
                \State $\mathsf{Slc} \gets $ \Call{Slices}{\tikzmarknode{call-slices}{$R$}}
                \If{$\mathsf{Slc}$ is None}
                    \State \Return{False}
                \EndIf
                \If{$\dim R > 1$}
                    \State extend queue with $\mathsf{Slc}$
                \EndIf
            \EndWhile
            \State \Return{True}
        \EndFunction
    \end{algorithmic}
\end{minipage}\ \ \ \ \ \ \ 
\begin{minipage}[t]{0.5\linewidth}
    {\footnotesize
    \begin{algorithmic}[1]
        \ttfamily
        \Function{\tikzmarknode{def-slices}{S}lices}{$P$}
            \State $\mathsf{O} \gets$ \Call{TopologicalOrder}{$\mathsf{M}P$}
            \If{$\mathsf{O}$ is None}
                \State \Return{None}
            \EndIf
            \State $\mathsf{Slc} \gets \{\!\!\!$ \Call{Boundary}{$-$, $\dim P - 1$, $P$} $\!\!\!\}$
            \For{$x \in \mathsf{O}$}
                \State $S \gets \left( \mathsf{Slc}[-1] \setminus \downarrow \Delta^- x \right) \cup \downarrow \Delta^+ x$
                \State add $S$ to $\mathsf{Slc}$
            \EndFor
            \State \Return{$\mathsf{Slc}$}
        \EndFunction
    \end{algorithmic}}
\end{minipage}
\begin{tikzpicture}[remember picture,overlay]
    \draw[->, >=Stealth, thick, black, shorten >=0.3cm, shorten <=2.2cm] ({def-slices.west}) to[out=180,in=0,looseness=1.25] (call-slices.east);
\end{tikzpicture}
\vspace{0.3cm}

After these checks we can organize our candidates into orbits of the $\mathbb{Z}_2^n$ action, which reverses arrows, and then inspect the list. In dimensions $\leq 4$ there are few enough candidates that we can double-check that all of them are indeed molecules.

\section{Table of Simplicial 4-Atoms}\label{app:4_atoms_list}

Using the above logic one can generate all simplicial $n$-atoms, for $n \leq 4$, which are guaranteed to be frame-acyclic by \Cref{lem:fa_less5_automatic}. The set of $n$-atoms has a $\mathbb{Z}_2^n$ action where each factor flips the sources and targets in some dimension $k$, with $0 < k \leq n$. Omitting the cases with a unique element, $n=0$ and $n=1$, the numbers are
\begin{center}
    \begin{tabular}{c|c|c}
        Dimension $n$ & $\#$ of simplicial $n$-atoms & $\#$ of $\mathbb{Z}_2^n$-orbits \\
        \hline
        2 & 2 & 1 \\
        3 & 16 & 3 \\
        4 & 1240 & 80 \\
    \end{tabular}
\end{center}

Thus, below we will draw one representative of each $\mathbb{Z}_2^4$-orbit. For any such representative, the diagram shows the five different 2-dimensional slices of a simplicial 4-atom at the vertices of a pentagon. These 2-dimensional slices are drawn with their string diagram representation, instead of as their polytope representation, to make their categorical reading easier. The 3-cells are denoted by the edges of the pentagon, which are unlabeled in the pictures to avoid clutter, but they carry the unique 3-cell which is a whiskering of the generators, and whose source and target are the fixed 2-cells. Finally, the 4-cell is either of the 2-3 or 1-4 kind, and sits in the middle of the pentagon.

While not necessarily a symmetry there is a further pattern that becomes clear from the pictures of the elements below, which is a cyclic, or reflected cyclic, permutation of the vertices of the pentagon. The reason it is not a symmetry, is because it does not respect all 3-cell directionality. The extra marking of a numbered, yellow, lower left corner on certain elements denotes the beginning of a block in which elements are related by such a permutation. Such a block ends when the following block begins. This is very relevant in the case of $\mathscr{P}_3$ where if the 3-generators happen to be such that the plain and anti versions are inverses of each other (as is the case for the special orbifold data of \Cref{ssec:special_orbifold_data} and \cite{orbifoldnd,orbifold3d}), then every block will only contribute one independent $\mathbb{Z}_2^4$-orbit of relations. There are 18 of these marked blocks.

\begin{theorem}\label{thm:satm_4}
    The $\mathbb{Z}_2^4$-orbits of simplicial 4-atoms are represented by the following list of elements:
    %
    \begin{center}
}};

        \end{tikzpicture}
    \end{center}
\end{theorem}

\printbibliography

@article{metayer,
  title={Cofibrant objects among higher-dimensional categories},
  author={François M{\'e}tayer},
  journal={Homology, Homotopy and Applications},
  year={2008},
  volume={10},
  pages={181-203},
  url={https://intlpress.com/JDetail/1805807090826362882},
  eprint={hal-00259476},
  archivePrefix = {HAL},
}

@book{ara,
   title={Polygraphs: From Rewriting to Higher Categories},
   ISBN={9781009498982},
   DOI={10.1017/9781009498968},
   publisher={Cambridge University Press},
   author={Ara, Dimitri and Burroni, Albert and Guiraud, Yves and Malbos, Philippe and Métayer, François and Mimram, Samuel},
   year={2025},
   month=Mar,
   eprint={2312.00429},
   archivePrefix={arXiv},
   primaryClass={math.CT}
}

@article{hadzihasanovic,
      title={Combinatorics of higher-categorical diagrams},
      author={Amar Hadzihasanovic},
      year={2024},
      eprint={2404.07273},
      archivePrefix={arXiv},
      primaryClass={math.CT},
      note = {[preprint]}
}

@article{condensation,
    author = "Gaiotto, Davide and Johnson-Freyd, Theo",
    title = "{Condensations in higher categories}",
    eprint = "1905.09566",
    archivePrefix = "arXiv",
    primaryClass = "math.CT",
    month = "5",
    year = "2019",
    note = "[preprint]"
}

@inbook{ising,
   title={Topological Defects},
   ISBN={9780323957069},
   DOI={10.1016/b978-0-323-95703-8.00098-7},
   booktitle={Encyclopedia of Mathematical Physics},
   author={Carqueville, Nils and Del Zotto, Michele and Runkel, Ingo},
   year={2025},
   pages={621–647},
   eprint={2311.02449},
   archivePrefix={arXiv},
   primaryClass={math-ph},
   editor = {Richard Szabo and Martin Bojowald},
   edition = {Second Edition},
}

@article{orbifold2d,
   title={Orbifold completion of defect bicategories},
   volume={7},
   ISSN={1664-073X},
   DOI={10.4171/qt/76},
   number={2},
   journal={Quantum Topology},
   publisher={European Mathematical Society - EMS - Publishing House GmbH},
   author={Carqueville, Nils and Runkel, Ingo},
   year={2016},
   month=Feb,
   pages={203–279},
   eprint={1210.6363},
   archivePrefix={arXiv},
   primaryClass={math.QA}
}

@article{dkr,
    author = "Davydov, Alexei and Kong, Liang and Runkel, Ingo",
    title = "{Field theories with defects and the centre functor}",
    eprint = "1107.0495",
    archivePrefix = "arXiv",
    primaryClass = "math.QA",
    reportNumber = "ZMP-HH-11-8",
    pages = "71--128",
    month = "7",
    year = "2011"
}

@article{pachner,
    title = {Shellings of simplicial balls and P.L. manifolds with boundary},
    journal = {Discrete Mathematics},
    volume = {81},
    number = {1},
    pages = {37-47},
    year = {1990},
    issn = {0012-365X},
    doi = {10.1016/0012-365X(90)90178-K},
    author = {Udo Pachner}
}

@article{graycat,
    author = "Carqueville, Nils and Meusburger, Catherine and Schaumann, Gregor",
    title = "{3-dimensional defect TQFTs and their tricategories}",
    eprint = "1603.01171",
    archivePrefix = "arXiv",
    primaryClass = "math.QA",
    doi = "10.1016/j.aim.2020.107024",
    journal = "Advances in Mathematics",
    volume = "364",
    pages = "107024",
    year = "2020"
}

@article{cobordismhypothesis,
    author = "Baez, J. C. and Dolan, J.",
    title = "{Higher dimensional algebra and topological quantum field theory}",
    eprint = "q-alg/9503002",
    archivePrefix = "arXiv",
    doi = "10.1063/1.531236",
    journal = "Journal of Mathematical Physics",
    volume = "36",
    pages = "6073--6105",
    year = "1995"
}

@article{orbifoldnd,
    author = "Carqueville, Nils and Runkel, Ingo and Schaumann, Gregor",
    title = "{Orbifolds of n-dimensional defect TQFTs}",
    eprint = "1705.06085",
    archivePrefix = "arXiv",
    primaryClass = "math.QA",
    doi = "10.2140/gt.2019.23.781",
    journal = "Geometry \& Topology",
    volume = "23",
    pages = "781--864",
    year = "2019"
}

@article{orbifold3d,
    author = {Carqueville, Nils and M{\"u}ller, Lukas},
    title = "{Orbifold Completion of 3-Categories}",
    eprint = "2307.06485",
    archivePrefix = "arXiv",
    primaryClass = "math.QA",
    doi = "10.1007/s00220-025-05434-y",
    journal = "Communications in Mathematical Physics",
    volume = "407",
    number = "1",
    pages = "8",
    year = "2026"
}

@inproceedings{ffrs,
    author = "{Fr\"ohlich}, {J\"urg} and {Fuchs}, {J\"urgen} and {Runkel}, {Ingo} and {Schweigert}, {Christoph}",
    title = {Defect Lines, Dualities and Generalised Orbifolds},
    booktitle = {16th International Congress on Mathematical Physics},
    eprint = "0909.5013",
    archivePrefix = "arXiv",
    primaryClass = "math-ph",
    reportNumber = "KCL-MTH-09-10, ZMP-HH-09-20",
    doi = "10.1142/9789814304634_0056",
    pages = "608--613",
    year = "2010"
}

@article{cobordismhypothesissingularities,
    title={On the Classification of Topological Field Theories},
    author={Lurie, Jacob},
    journal={Current Developments in Mathematics},
    volume={2008},
    number={1},
    pages={129--280},
    year={2008},
    publisher={International Press of Boston},
    eprint={0905.0465},
    archivePrefix={arXiv},
    primaryClass={math.CT}
}

@incollection{orbifoldsurvey,
    title = {Orbifolds of Topological Quantum Field Theories},
    editor = {Richard Szabo and Martin Bojowald},
    booktitle = {Encyclopedia of Mathematical Physics},
    edition = {Second Edition},
    pages = {618-634},
    year = {2025},
    doi = {10.1016/B978-0-323-95703-8.00008-2},
    author = {Nils Carqueville},
    ISBN={9780323957069},
    eprint={2307.16674},
    archivePrefix={arXiv},
    primaryClass={math-ph}
}

@article{graphicalcalculus,
   title={Gray categories with duals and their diagrams},
   volume={450},
   ISSN={0001-8708},
   url={http://dx.doi.org/10.1016/j.aim.2024.109740},
   DOI={10.1016/j.aim.2024.109740},
   journal={Advances in Mathematics},
   author={Barrett, John W. and Meusburger, Catherine and Schaumann, Gregor},
   year={2024},
   month={7},
   pages={109740}
}

@online{pythoncode,
  author       = {Ivanov, Aleksandar},
  title        = {simplicial-atoms},
  shorthand    = {sim},
  date         = {2026-11-11},
  url          = {https://github.com/ackiivanov/simplicial-atoms},
  urldate      = {2026-11-13},
  license      = {MIT License},
  repository   = {https://github.com/ackiivanov/simplicial-atoms},
  note         = {GitHub repository}
}

@article{sphericalfusion,
    author = "Carqueville, Nils and Runkel, Ingo and Schaumann, Gregor",
    title = "Orbifolds of Reshetikhin-Turaev TQFTs",
    eprint = "1809.01483",
    archivePrefix = "arXiv",
    primaryClass = "math.QA",
    journal = "Theory and Applications of Categories",
    volume = "35",
    pages = "513--561",
    year = "2020"
}

\authoraddress{University of Vienna, Faculty of Physics, Mathematical Physics Group, Boltzmanngasse 5, 1090 Vienna, Austria}{aleksandar.ivanov@univie.ac.at}

\end{document}